\documentclass[11pt]{article}
\usepackage[margin=1in]{geometry}
\usepackage{amsmath,amssymb,amsthm,mathtools}
\usepackage[T1]{fontenc}
\usepackage{xcolor}
\usepackage{hyperref}
\hypersetup{colorlinks=true,linkcolor=red,citecolor=blue,urlcolor=blue}

\newtheorem{theorem}{Theorem}[section]
\newtheorem{proposition}[theorem]{Proposition}
\newtheorem{lemma}[theorem]{Lemma}
\newtheorem{corollary}[theorem]{Corollary}
\newtheorem{conjecture}[theorem]{Conjecture}

\theoremstyle{definition}
\newtheorem{definition}[theorem]{Definition}
\newtheorem{remark}[theorem]{Remark}

\newcommand{\Q}{\mathbb Q}
\newcommand{\C}{\mathbb C}
\newcommand{\Z}{\mathbb Z}
\newcommand{\PP}{\mathbb P}
\newcommand{\ft}{\mathfrak t}

\newcommand{\ev}{\mathrm{ev}}
\newcommand{\Mbar}{\overline{\mathcal M}}
\newcommand{\Frac}{\mathrm{Frac}}

\newcommand{\wts}{\mathrm{wts}}
\newcommand{\Gh}{\widehat\Gamma}
\newcommand{\Ghs}{\widehat\Gamma^{*}}
\DeclareMathOperator{\Sym}{Sym}
\DeclareMathOperator{\Aut}{Aut}
\DeclareMathOperator{\val}{val}
\DeclareMathOperator*{\Res}{Res}

\title{Convergence of the big J-function of a GKM manifold}
\author{Todor Milanov}
\date{}

\begin{document}
\maketitle

\begin{abstract}
The genus-$0$ equivariant Gromov--Witten (GW) invariants of a GKM manifold can be
computed recursively by a recursion discovered by Brown and studied further by Coates--Corti--Iritani--Tseng. In this paper, we prove a slightly modified version of this result by realizing the J-function as a solution to a fixed-point problem. As an application, 
we prove that both the equivariant and the non-equivariant quantum cohomologies of a GKM manifold are convergent. We also show that the equivariant central charges have no poles in the spectral parameter, compute their monodromy and their asymptotics at $z=0$, and formulate an equivariant mirror symmetry conjecture for GKM manifolds.
\end{abstract}

\tableofcontents

\section{Introduction and summary}

\subsection{Overview}\label{sec:overview}

The higher-genus reconstruction in Gromov--Witten (GW) theory was
conjectured by Givental in \cite{Givental2001a} and proved by Teleman
in \cite{Teleman2012}. While this solves the problem of computing the
higher-genus GW invariants in terms of genus-$0$ ones, the problem of
computing genus-$0$ GW invariants is still quite challenging. There
are several interesting conjectures, such as the mirror symmetry
conjecture and the refined Dubrovin conjecture where having an
efficient method for computing the genus-$0$ GW invariants is
crucial. Nevertheless, there was a considerable progress in the case when 
the target manifold is a Goresky--Kottwitz--MacPherson (GKM) manifold -- see Section
\ref{sec:gkm} for some background. Using fixed-point localization, Jeff Brown \cite{Brown2014} (see also \cite{CCIT2015}) was able to find a  
recursion that in principle can be used to compute all genus-$0$ GW invariants in terms of
the combinatorics of the GKM graph. In this paper, we would like to report on an interesting twist in carrying out the fixed-point localization computation that was discovered with the AI assistance. Namely, the so-called {\em Abel's identity} (see \cite{Abel1826}) allows us to obtain a closed formula for the graph sum in the genus-$0$ fixed point localization formula. Of course, our result can be derived easily by using the results in \cite{Brown2014} or \cite{CCIT2015}. Nevertheless, we decided to give a self-contained proof because the Abel's identity might be useful in other more sophisticated computations and also our results are in more general settings. The main feature of our formula is that it realizes the big J-function as a solution to a fixed-point problem. This allows us to give an elementary proof of the convergence of the big J-function, i.e. the convergence of the equivariant quantum cohomology. Moreover, by passing to the non-equivariant limit, we obtain also the convergence of the non-equivariant J-function too. 

The analytic aspects of quantum cohomology have a long history.
Iritani~\cite[\S1]{Iritani2007} proved that the big quantum cohomology of
a smooth projective toric variety is convergent and generically
semisimple, and announced that the big quantum cohomology, as well as the
higher genus potentials, converge for any projective manifold admitting a
Hamiltonian torus action with only isolated fixed points and isolated
one-dimensional orbits, that is, for GKM manifolds. More recently,
convergence and semisimplicity of the big quantum cohomology of flag
varieties $G/P$ were established in~\cite{LLYZ2026}. Theorem~\ref{thm:J-convergence}
gives a convergence result for the $J$-function in the GKM setting.
The convergence for GKM manifolds was announced by Iritani in an abstract for the Fall 2005 meeting of the Mathematical
Society of Japan and in lecture notes for a meeting for graduate students in Kinosaki, but a proof has not been
published. Iritani's method is more straightforward than ours: it estimates directly the graph sums in the fixed-point
localization formulas, and it yields the convergence of all genus-$g$ descendant potentials without relying on
Givental's higher-genus reconstruction (private communication). Our approach is different: the convergence is
obtained from the realization of the big $J$-function as the solution of a fixed-point problem.

Finally, motivated by the gamma conjectures in quantum cohomology (see \cite{GGI2016} and \cite{CDG2024}), we compare the poles of the equivariant J-function with respect to spectral parameter $z$ with the gamma class of the manifold $X$. Namely, the J-function is a sum of elementary fractions whose poles can be canceled by the poles of the equivariant gamma class of the target manifold. We prove that the ratio of the restrictions of the $J$-function and the gamma class to the fixed points satisfy an interesting relation. As an application, we show that the equivariant central charges, i.e., the pairings of the $J$-function with the gamma-integral classes of equivariant vector bundles, have no poles in $z\in\C^\times$, we compute their monodromy, and we determine their asymptotics as $z\to0$. In the rest of this introduction we will focus in formulating our
results.

\subsection{Notation and conventions}\label{sec:notation}

\paragraph{Tori conventions and notation.}
Through out the paper $T\cong (\C^*)^\ell$ will be a fixed algebraic torus. Put  
$\ft:=\mathrm{Lie}(T)$ for the Lie algebra of $T$ and $\ft^*:=\mathrm{Hom}_\C(\ft,\C)$ for the linear forms on $\ft$. Recall that a character of $T$ is just a group homomorphism $\chi:T\to \C^*$ and 
let $\C_\chi$ be the corresponding one-dimensional representation of $T$ .
The set of all characters of $T$ is denoted by $X^*(T)$. 
Taking the differential at $1\in T$ defines an embedding 
\[
X^*(T)\to \ft^*,\qquad \chi\mapsto d\chi
\]
Its image $\ft^*_{\mathbb Z}$ is a free abelian group of rank $\dim T$ spanning $\ft^*$. 
We usually identify $X^*(T)\cong \ft^*_{\mathbb Z}$ and refer to both of them as the {\em character lattice} of $T$. Nevertheless, in order to make a distinction, 
we call the elements of $X^*(T)$ characters, those of 
$\ft^*_\Z$ -- \emph{weights}, and for $t\in T$ and a weight $m=d\chi$ we define $t^m:=\chi(t)$. 
Note that, since $d(\chi\chi')=d\chi+d\chi'$, we have $t^{m}t^{m'}=t^{m+m'}$.
Similarly, we embed the lattice of cocharacters $X_*(T):=\mathrm{Hom}(\C^*,T)$ into $\ft$ via $\sigma\mapsto d\sigma(1)$, where
$d\sigma:\C=\mathrm{Lie}(\C^*)\to\ft$ is the differential of $\sigma$ at $1$, and we denote by $\ft_\Z$ the image of $X_*(T)$. If
$\chi$ is a character and $\sigma$ a cocharacter, then $\chi(\sigma(s))=s^{\langle\chi,\sigma\rangle}$ for some integer $\langle\chi,\sigma\rangle$, and
$d\chi\big(d\sigma(1)\big)=\langle\chi,\sigma\rangle$. Therefore the two lattices $\ft_\Z$ and $\ft^*_\Z$ are dual to each other.

Let $K^0_T({\rm pt})$ be the representation ring of virtual $T$-modules.
Recall that every $T$-module $V$ decomposes as a direct sum of its weight subspaces:  
\[
V=\oplus_{\mu \in \ft^*_\Z}  V_\mu,\qquad
V_\mu:=\{v\in V\mid t v=t^\mu v\ \forall t\in T\},
\]
where $t\mapsto t^\mu$ is the character corresponding to the weight $\mu$. 
Using the weight space decomposition, we define the equivariant \emph{Euler class}
\[
e^T: K^0_T({\rm pt})\to \C[\ft]:=\Sym(\ft^*),\qquad 
e^T(V):= \prod_\mu \mu^{\dim V_\mu},
\]
and the \emph{character map} 
\[
\operatorname{ch}_T:K^0_T({\rm pt})\to\C[\![\ft]\!],\qquad 
[V]\mapsto\sum_{\mu}(\dim V_\mu)\,e^\mu,\qquad 
e^\mu:=\exp(\mu),
\]
where $\C[\![\ft]\!]$ is the completion of $\C[\ft ]$ at $0$. 
The Euler class satisfies $e^T(V+W)=e^T(V)e^T(W)$ while 
$\operatorname{ch}_T$ is an injective ring homomorphism.

Finally, we define a ring isomorphism $\operatorname{Sym}(\ft^*_\Z)\cong H^*_T(\mathrm{pt})$ such that a weight $\mu=d\chi \in \mathfrak t^*_{\mathbb Z}$ is mapped to $c_1(L_\chi)$ where the line bundle $L_\chi := (ET\times \mathbb C)/ T$ on $BT$ is defined through the {\em right} action of $T$ on $ET\times \mathbb C$, i.e. $(v,z)\cdot t = (vt,t^{-\mu} z)$. Note that when $T=\C^*$ and $\chi=\operatorname{id}$, we have $BT=\C P^\infty$ and $L_\chi=\mathcal O(-1)$. Under this identification, the character map $\operatorname{ch}_T: K^0_T({\rm pt})\to \widehat H^*_T({\rm pt})$ can be interpreted as an equivariant version of the Chern character map.

\paragraph{Equivariant GW theory.} We refer to \cite{Givental1996} and 
\cite{GraberPandharipande1999} for more 
details. Let $X$ be a smooth projective variety equipped with the action of $T$.  
Let $\overline{\mathcal M}_{0,n+1}(X,\beta)$ be the moduli space of
stable maps $f:(C,z_0,z_1,\dots,z_n)\to X$ where $C$ is a genus-$0$
nodal curve equipped with $n+1$ marked points and $f$ is a holomorphic
map of degree $\beta$, i.e. $f_*[C]=\beta$. The $T$-action on $X$
induces a $T$-action on the moduli space and the virtual fundamental
cycle (see \cite{BehrendFantechi1997}) has an equivariant counterpart
(see \cite[\S1]{GraberPandharipande1999} and \cite{EdidinGraham1998}
for equivariant Chow groups) 
\[
[\overline{\mathcal M}_{0,n+1}(X,\beta)]^{\rm vir}\in 
H^T_*(\overline{\mathcal M}_{0,n+1}(X,\beta))
\]
which allows us to define the equivariant pushforward map
\[
\ev^T_{0*}: H_T^*(\overline{\mathcal M}_{0,n+1}(X,\beta))\to H_T^*(X),
\qquad
\ev^T_{0*}(\alpha):= \mathrm{PD}_X^{-1}\,\ev_{0*} (\alpha \cap [\overline{\mathcal M}_{0,n+1}(X,\beta)]^{\rm vir}),
\]
where $\ev_0: \overline{\mathcal M}_{0,n+1}(X,\beta)\to X$ is the evaluation map at the $0$th marked point.
Here the virtual class is regarded as an equivariant homology class and $\mathrm{PD}_X$ denotes equivariant Poincar\'e
duality on $X$; the precise conventions are recalled in \S\ref{sec:eqhom}.
The $T$-equivariant big $J$-function of $X$ is by definition the following formal power series:
\begin{equation}\label{eq:Jdef}
J(\tau,z):=-z+\tau+\sum_{\beta\in {\rm NE}(X)}\sum_{n=0}^\infty
\frac{q^\beta}{n!}\,
\ev^T_{0*}\Big(\frac{1}{-z-\psi_0} \cup 
\ev_1^*\tau \cup\dots \cup \ev_n^*\tau
\Big),
\end{equation}
where ${\rm NE}(X)\subset H_2(X,\Z)$ is the monoid of effective curve classes and 
the RHS should be expanded as a power series in $z^{-1}$. Here and below, the correlators corresponding to the unstable
moduli spaces, i.e. $(\beta,n)=(0,0)$ and $(0,1)$, are defined to be $0$.
Equivalently, by separating the $\beta=0$ terms from the sum, we can write
\[
J(\tau,z)=-z\,e^{-\tau/z}+\sum_{\beta\ne0}q^\beta J_\beta(\tau,z)
\]
with
\[
J_\beta(\tau,z)=\sum_n\frac1{n!}\ev^T_{0*}\Big(
\frac{1}{-z-\psi}\cup \ev_1^*\tau\cup \dots\cup \ev_n^*\tau\Big),
\]
where on the RHS we suppressed the dependence on the degree $\beta$. Let us fix our interpretation 
of the Novikov variables. Let $L_1,L_2,\dots,L_r$ be ample line bundles such that $c_1(L_i)$ form a $\Z$-basis of the lattice $H^2(X;\Z)\cap H^{1,1}(X)$. Then $\langle c_1(L_i),\beta\rangle $ is a non-negative integer for all $\beta\in \operatorname{NE}(X)$. The Novikov variables $q^\beta$ are identified with elements of the ring $\C[[q]]:=\C[[q_1,\dots,q_r]]$ via 
$q^\beta := q_1^{\langle c_1(L_1),\beta\rangle} \dots q_r^{\langle c_1(L_r),\beta\rangle} $.

\subsection{GKM manifolds and the big \texorpdfstring{$J$}{J}-function}\label{sec:intro-gkm}

Let $X$ be a GKM manifold (Definition~\ref{def:gkm}). By definition a $T$-curve in $X$ is a closed one-dimensional orbit of $T$. Every $T$-curve $C$ is a smooth submanifold of $X$ isomorphic to $\PP^1$ that carries precisely two fixed points -- see Proposition \ref{prop:gkmcurves}. Therefore, it is natural to define a graph, called the {\em GKM graph} of $X$, whose vertices are the fixed points $X^T$ and whose edges are the $T$-curves. For each vertex $p\in X^T$, the tangent space $T_pX$ is a $T$-module and we put $e_p:=e^T(T_pX)$, i.e. the product of the weights of $T_pX$. If $C$ is a $T$-curve connecting two fixed points $p$ and $p'\in X^T$ and $k\geq 1$ is an integer number, then we define a rational function $E_{k}(C)\in \operatorname{Frac}(\C[\ft])$ as follows.  Let us fix a parametrization of $C$, i.e. an isomorphism
\[
\phi:\PP^1\to C,\qquad \phi([1,0])=p,\qquad \phi([0,1])=p'
\]
such that the induced $T$-action on $\PP^1$ is given by $t\cdot [z_0,z_1]=[z_0, t^\epsilon z_1]$ for some weight $\epsilon \in \ft^*_\Z$. Since $d\phi$ induces an isomorphism between $T_{[1,0]}\PP^1$ and $T_pC$ we get that $\epsilon$ is the weight of the $T$-action on $T_pC$. Furthermore, a $T$-equivariant line bundle on $\PP^1$ must have the form 
\[
L(a,u):= ((\C^2\setminus{0})\times \C)/\C^*,\qquad
t\cdot[z_0,z_1,x]:=[z_0,t^\epsilon z_1, t^u x],\qquad
t\in T,
\]
where $a\in \Z$ is an integer, $u\in \ft^*_\Z$ is a weight, and the $\C^*$-action in the quotient defining $L(a,u)$ is $\lambda\cdot (z_0,z_1,x)=(\lambda z_0,\lambda z_1, \lambda^a x)$. The number $a$ is called the degree, $\epsilon$ -- the base weight at $p$,  and $u$ the fiber weight at $p$. It is a result of Kumar \cite{Kumar2003} (see also Proposition \ref{prop:eqBG}) that every $T$-equivariant vector bundle on $\PP^1$ decomposes as a direct sum of equivariant line bundles, i.e line bundles of the form $L(a,u)$. The restriction of the tangent bundle $TX|_C$ is an equivariant $T$-bundle on $C$ and hence $\phi^*TX$ decomposes as a direct sum of equivariant line bundles $L(a_i,u_i)$ ($1\leq i\leq \operatorname{dim} X$). Since one of these line bundles is $T_pC=L(2,\epsilon)$ we may assume that $a_1=2$ and $u_1=\epsilon$. Then the edge factor $E_{k}(C)$ is defined as follows:
\begin{equation}\label{eq:edge}
E_{k}(C)=\prod_{\substack{m=0\\ m\ne k}}^{2k}\Big(1-\frac mk\Big)\varepsilon\cdot
\prod_{\substack{i\ge2\\ a_i\ge0}}\ \prod_{m=0}^{ka_i}\Big(u_i-\frac mk\varepsilon\Big)\cdot
\prod_{\substack{i\ge2\\ a_i<0}}\ \prod_{m=ka_i+1}^{-1}\Big(u_i-\frac mk\varepsilon\Big)^{-1}.
\end{equation}
Note that the base and fiber weights of $L(a_i,u_i)$ at the other fixed point $p'$ are respectively $\epsilon'=-\epsilon$ and $u_i'=u_i - a_i \epsilon$. It is straightforward to check that the definition of the edge factor $E_{k}(C)$ is independent of the choice of a fixed point $p\in C^T$. In fact, as it will be proved later on, $E_{k}(C)$ is the equivariant Euler class of the moving part, i.e. the product of the non-zero weights, of the virtual $T$-representation $H^0(\PP^1,f^*TX)-H^1(\PP^1,f^*TX)$ where $f:\PP^1\to X$ is the stable map obtained from $\phi$ and the standard $k$-fold covering map $\PP^1\to \PP^1$, $[z_0,z_1]\mapsto [z_0^k,z_1^k]$. 

Let us denote by $B$ the set of all triples $(p,C,k)$ consisting of a fixed point $p\in X^T$, 
a $T$-curve $C$ passing through $p$, and an integer $k\geq 1$. The pair $F=(p,C)$ will be sometimes referred to as a \emph{GKM-flag}. Put
\[
z_b:=\frac{\varepsilon}{k},\qquad
R_b(\tau):=
\frac{e_p}{k\,E_{k}(C)}\; q^{k[C]}\; J_{p'}(\tau,z_b), \qquad b=(p,C,k)\in B,
\]
where $\varepsilon:=e^T(T_pC)$ is the weight of $T_pC$, $[C]\in \operatorname{NE}(X)$ is the homology class of $C$, $p'\in X^T$ is the fixed point connected to $p$ via $C$, and $J_p(\tau,z):=J(\tau,z)|_{p}$ is the restriction to $p$. 
\begin{theorem}\label{thm:J}
Let $u_p$ ($p\in X^T$) be the unique solution of the fixed-point problem
\begin{equation}\label{eq:up}
u_p=\tau|_p+\sum_{b\in B(p)}
\tfrac{R_{b}(\tau)}{z_b}\, 
e^{u_p/z_b} ,
\end{equation}
where $B(p):=\{b\in B\ |\ p\in b\}$. Then
\begin{equation}\label{eq:Jp}
J_p(\tau,z)=-z\,e^{-u_p/z}\Big(1+\sum_{b\in B(p)}
\tfrac{R_{b}(\tau)}{z_b(z-z_b)}\,
e^{u_p/z_b}\Big).
\end{equation}
Moreover, $u_p$ are canonical coordinates for the equivariant quantum cohomology.
\end{theorem}

The solution $u_p$ in Theorem \ref{thm:J} is constructed  by using successive iterations: start with $u_p^{(0)}:=0$, substitute it for $u_p$ in the RHS and call the result $u_p^{(1)}$. By repeating this process we obtain a sequence $u_p^{(n)}$ of formal power series in $q$ and $\tau$ which converges in the formal $I$-adic topology defined by the ideal $I=(q,\tau)$. The formula for $J_p(\tau,z)$ implies that 
\begin{equation}\label{eq:J-res}
R_{b}(\tau)=-\operatorname{Res}_{z=z_b} J_p(\tau,z) dz.
\end{equation}
Since $R_{b}(\tau)$ involves $J_{p'}$ with $p'\neq p$ and has a factor of $q^{k[C]}$, Theorem \ref{thm:J} computes the big $J$-function order by order in $q$, for all $\tau$ at once. 

\subsection{Convergence}\label{sec:intro-convergence}

Formulas \eqref{eq:up}--\eqref{eq:J-res} can be interpreted as a fixed-point problem for 
$u_p$ ($p\in X^T$) and $R_b$ $(b\in B)$. We will prove in Section \ref{sec:analytic} that if we specialize 
all equivariant parameters in $H^*_T({\rm pt})=\C[\ft]$ to a generic $\lambda\in \ft$, then 
the fixed-point problem can be interpreted analytically and the formal solution in Theorem \ref{thm:J} is in fact convergent.  
\begin{theorem}\label{thm:J-convergence}
Let $X$ be a GKM manifold. 
For every admissible $\lambda\in \ft$ (see Definition \ref{def:admissible}; the admissible points form an open dense
subset of $\ft$) there is a real number $r=r(\lambda)>0$
such that for every $p\in X^T$, the specialization of $J_p(\tau,z)$ to $\lambda$ 
converges for $|q^{[C]}|<r$ (all $T$-curves $C$) and $|\tau|_p|<r$ (all $p\in X^T$),
and its sum extends to a meromorphic function of
$z\in\C^\times$ whose poles are simple and contained in 
$\{z_{b}\}_{b\in B(p)}$. Moreover, $J_p(\tau,z)+z$ is holomorphic at $z=\infty$.
\end{theorem}
By taking the non-equivariant limit and using that the quantum differential equations have a regular singularity 
at the large radius limit $q=0$, we obtain the analyticity of the non-equivariant J-function. 
\begin{corollary}\label{cor:none-J}
Let $X$ be a GKM manifold. Then there exists a positive number $r>0$ such that the non-equivariant big J-function $J(\tau,z)$ of $X$ is convergent and it extends to an analytic function for all 
\[
|q|<r, \qquad |\tau|<r,\qquad  z\in \PP^1\setminus\{0\}.
\]
\end{corollary}

\subsection{Equivariant central charges and mirror symmetry}\label{sec:intro-central}

Partially motivated by the gamma conjectures of Galkin--Golyshev--Iritani \cite{GGI2016} and 
the refined Dubrovin conjecture in \cite{CDG2024}, 
we would like to investigate the equivariant central charges of $X$. Put $n:=\dim_\C X$. Since $T$ is a torus and $X$ is
smooth, the ample line bundles $L_1,\dots,L_r$ used to define the Novikov variables admit $T$-equivariant structures
\cite{KKLV1989}; we fix such structures and put $P_i:=c_1^T(L_i)\in H^2_T(X)$. Let $K^0_T(X)$ be the Grothendieck group of
$T$-equivariant vector bundles on $X$, let $\operatorname{ch}_T:K^0_T(X)\to\widehat H^*_T(X)$ be the equivariant Chern character,
and let $\Gh_X:=\prod_x\Gamma(1+x)$ be the equivariant Gamma class, the product running over the equivariant Chern roots
$x$ of $TX$. Let $\theta$ be the \emph{degree operator}, $\theta(\alpha):=(\tfrac n2-\deg\alpha)\,\alpha$ for $\alpha\in H^*_T(X)$
homogeneous of complex degree $\deg\alpha$; here the equivariant parameters have degree $1$, so $\theta$ is not
$H^*_T({\rm pt})$-linear. Let $\rho:=c_1^T(TX)\cup$. Following Iritani \cite{Iritani2009}, for $E\in K^0_T(X)$ we put
\[
\Psi(E):=(2\pi)^{\frac{1-n}{2}}\,\Gh_X\cup q_1^{-P_1}\cdots q_r^{-P_r}\cup(2\pi i)^{\deg}\operatorname{ch}_T(E),
\]
\[
Z(E)(\tau,z):=\frac{1}{\sqrt{2\pi}}\int_X\frac{J(\tau,z)}{-z}\cup z^{\theta}z^{\rho}\,\Psi(E),
\]
where $q_i^{-P_i}:=e^{-P_i\log q_i}$, $z^\rho:=e^{\rho\log z}$, and $(2\pi i)^{\deg}$ multiplies the degree-$k$ component by
$(2\pi i)^k$. We call $Z(E)$ the \emph{equivariant central charge} of $E$. The operator $z^\theta$ multiplies a homogeneous class
of degree $k$ by $z^{n/2-k}$; in particular it rescales every weight $u$ to $u/z$. Therefore, after specializing the
equivariant parameters to an admissible $\lambda\in\ft$ and applying the localization formula, we get
\begin{multline}\label{eq:Z-loc}
Z(E)(\tau,z)=(2\pi)^{-\frac n2}\,z^{\frac n2}\sum_{p\in X^T}\frac{J_p(\tau,z)}{-z\,e_p}\,
\prod_{u\in\wts(T_pX)}\Gamma\Big(1+\frac uz\Big)\\
\times e^{c_1^T(TX)|_p\log z/z}\;\prod_{i=1}^r e^{-P_i|_p\log q_i/z}\;
\operatorname{ch}_T(E_p)\Big(\frac{2\pi i}{z}\Big),
\end{multline}
where $z^{n/2}:=e^{\frac n2\log z}$, $E_p$ is the fiber of $E$ at $p$, and for a virtual $T$-module $V$ with
$\operatorname{ch}_T(V)=\sum_w m_w(V)\,e^w$ we write $\operatorname{ch}_T(V)(x):=\sum_w m_w(V)\,e^{x\,w(\lambda)}$. We take
\eqref{eq:Z-loc} as the definition of the specialized central charge. It depends on the choices of $\log z$ and of
$\log q_i$, i.e., $Z(E)$ is a function on the universal cover of the space of parameters $(q,z)$ with $q_i\neq0$ and
$z\neq0$. Recall that by Theorem~\ref{thm:J} the canonical coordinates $u_p$ and the residues $R_b$ are holomorphic in
$(q,\tau)$ (see also Proposition~\ref{prop:contraction}). Put
\[
g_p(0):=1-\sum_{b\in B(p)}\frac{R_b(\tau)}{z_b^2}\,e^{u_p/z_b},\qquad p\in X^T,
\]
i.e., the value at $z=0$ of the bracket in \eqref{eq:Jp}.
\begin{corollary}\label{cor:central-charge}
Let $X$ be a GKM manifold and $\lambda$ admissible. Then there exists a positive number $r_0=r_0(\lambda)\leq1$ such that for every
$E\in K^0_T(X)$ the central charge $Z(E)(\tau,z)$ is a holomorphic function on the universal cover of
\[
0<|q_i|<r_0\ (1\leq i\leq r),\qquad |\tau|_p|<r_0\ (p\in X^T),\qquad z\in\C^\times .
\]
Moreover, the following statements hold.
\begin{enumerate}
\item[(a)] The central charge has no poles in $z\in\C^\times$: although the individual terms in \eqref{eq:Z-loc} have
poles at the points $z_b$ ($b\in B$), these poles cancel in the sum.
\item[(b)] $Z(E)$ is multivalued in $z$ and in $q_i$, with the following monodromy. The analytic continuation along the loop
$z\mapsto e^{2\pi i s}z$, $0\leq s\leq1$, transforms $Z(E)$ into $(-1)^nZ(E\otimes K_X^{-1})$, where $K_X^{-1}=\det TX$ carries
its induced $T$-equivariant structure. The analytic continuation along the loop $q_i\mapsto e^{2\pi i s}q_i$, $0\leq s\leq1$,
transforms $Z(E)$ into $Z(E\otimes L_i^{-1})$.
\item[(c)] Let $S\subset\C^\times$ be a closed sector that contains none of the rays $\arg z=\arg u(\lambda)$ and
$\arg z=\arg u(\lambda)+\pi$, where $u$ runs over the weights of $T_pX$ for all $p\in X^T$, and fix a branch of $\log z$
on $S$. For each weight $u$ of $T_pX$ put $\ell_u:=\operatorname{Log}(u/z)+\log z$, where $\operatorname{Log}$ is the principal
branch; $\ell_u$ is a logarithm of $u$ independent of $z\in S$. Then, for fixed $(q,\tau)$, as $z\to0$ in $S$,
\[
Z(E)(\tau,z)=\sum_{p\in X^T}\big(\kappa_p\,g_p(0)+O(z)\big)\sum_{w}m_w(E_p)\,e^{-F_{p,w}/z},
\]
where the $O(z)$ term depends on $p$ but not on $E$,
\[
\kappa_p:=\frac{1}{e_p}\prod_{u\in\wts(T_pX)}e^{\ell_u/2}\quad(\text{so that }\kappa_p^2=e_p^{-1}),
\]
and
\[
F_{p,w}:=u_p(\tau)+\sum_{i=1}^rP_i|_p\log q_i-\sum_{u\in\wts(T_pX)}u\,(\ell_u-1)-2\pi i\,w .
\]
\end{enumerate}
\end{corollary}

Corollary~\ref{cor:central-charge} is consistent with the following equivariant mirror symmetry conjecture, which is
motivated by equivariant mirror symmetry for toric manifolds and partial flag manifolds: the equivariant central charges
of a GKM manifold should be oscillatory integrals. Let $T^\vee:=\mathrm{Hom}(\ft_\Z,\C^*)$ be the dual torus, so that $\mathrm{Lie}(T^\vee)=\mathrm{Hom}(\ft_\Z,\C)=\ft^*$. The exponential map $\ft^*\to T^\vee$, $\mu\mapsto\big(\xi\mapsto e^{\mu(\xi)}\big)$, identifies $T^\vee$ with $\ft^*/2\pi i\,\ft^*_\Z$; in particular, for a map $g$ to $T^\vee$, $\log g\in\ft^*$ is defined up to $2\pi i\,\ft^*_\Z$. For a holomorphic
map $g:Y\to T^\vee$ and $\lambda\in\ft$, the function $g^{\lambda/z}:=e^{\langle\lambda,\log g\rangle/z}$ is multivalued: changing the branch
of $\log g$ by $2\pi i\mu$, $\mu\in\ft^*_\Z=X^*(T)$, multiplies it by $e^{2\pi i\mu(\lambda)/z}$.
\begin{conjecture}\label{conj:mirror}
Let $X$ be a GKM manifold of dimension $n$. Then there exist
\begin{itemize}
\item a family of complex manifolds $\pi:Y\to\mathcal B$ with fibers $Y_s=\pi^{-1}(s)$ of dimension $n=\dim Y-\dim\mathcal B$, and a
local biholomorphism (the mirror map) from $\mathcal B$ to a neighborhood of the large radius limit in the space of
parameters $(q,\tau)$;
\item a holomorphic function $f:Y\to\C$ and a holomorphic map $g:Y\to T^\vee$;
\item a primitive form $\omega$ for the family of functions $F:=f-\langle\lambda,\log g\rangle$ in the sense of K.~Saito
(see \cite[Section 4.6]{MilanovSaito});
\end{itemize}
such that for every admissible $\lambda\in\ft$ and every $E\in K^0_T(X)$
\[
Z(E)(\tau,z)=\frac{1}{(2\pi z)^{n/2}}\int_{C(E)}e^{-f/z}\,g^{\lambda/z}\,\omega,
\]
where $C(E)\in H_n\big(Y_s,\operatorname{Re}(F/z)\gg0;\mathcal L\big)$ is an integration cycle with integer coefficients and
$\mathcal L$ is the local system on $Y\times\ft\times\C^\times$ that determines the values of $g^{\lambda/z}$, i.e., the branch of
$\log g$ along the cycle; the condition $\operatorname{Re}(F/z)\gg0$ is understood with this branch. Moreover, the map $E\mapsto
C(E)$ is an isomorphism of $\Z[X^*(T)]$-modules onto the lattice of such cycles, where a character $\mu$ acts on cycles by
the deck transformation of $\mathcal L$ that shifts $\log g$ by $2\pi i\mu$.
\end{conjecture}
The evidence for Conjecture~\ref{conj:mirror} is discussed in Section~\ref{sec:mirror-conj}.

\section*{Statement on the use of AI tools}

This work was carried out with substantial assistance from the AI system 
Claude (Anthropic; model Claude Opus 5.5, used from
28~September to 13~November 2026). The project grew out of an exercise. In order to understand the
semisimplicity of the quantum cohomology of $G/P$, the author asked Claude
to revisit Givental's fixed-point localization computation~\cite{Givental2001a}
and to look for a closed formula for the canonical coordinates. Claude
found a way to resum the localization sums using an identity of
Abel~\cite{Abel1826}. This unexpected outcome led the author to revisit
Givental's localization argument, which resulted in
Theorem~\ref{thm:J}. Although Theorem~\ref{thm:J} also admits an
elementary proof based on Brown's result~\cite{Brown2014}, the formulation
obtained in this way suggested Theorem~\ref{thm:J-convergence}, the convergence of
the $J$-function. Claude also assisted in finding the estimates needed in
the proof of Theorem~\ref{thm:J-convergence}. In addition, Claude was used to
search the literature for Section~\ref{sec:gkm}, to draft parts of the
Appendix, and to proofread the manuscript.

The author has checked all arguments, computations and references in the
paper, including those that originated with or were drafted by the AI, and
takes full responsibility for its content. 

\subsection*{Acknowledgments}
I am thankful to Chow Chi Hong for useful discussions on partial flag varieties, and to Hiroshi Iritani for the
communication about his unpublished work on the convergence of Gromov--Witten potentials of GKM manifolds. 
This work is supported by the World Premier International Research Center Initiative (WPI Initiative), MEXT, Japan and by JSPS Kakenhi Grant Number JP22K03265.

\clearpage
\section{GKM manifolds}\label{sec:gkm}

In this section we give the necessary background on GKM manifolds. 
For the sake of completeness, we would like to recall also some of the key proofs. 
\begin{definition}\label{def:gkm}
A \emph{GKM manifold} is a smooth projective complex variety $X$ with an algebraic action of an algebraic torus $T$ such
that
\begin{itemize}
\item[(G1)] the fixed point set $X^T$ is finite;
\item[(G2)] for every $p\in X^T$ the weights of the $T$-module $T_pX$ are pairwise linearly independent.
\end{itemize}
\end{definition}
From now on we will assume that $X$ is a GKM manifold. The main examples are given by toric and partial flag manifolds $G/P$ with maximal torus $T\subset P$. According to Bia{\l}ynicki-Birula
\cite{BialynickiBirula1973}, $X$ has a cell decomposition with $N:=|X^T|$ even-dimensional cells, so $H^{\mathrm{odd}}(X)=0$ and
$H^*_T(X)$ is a free $\C[\ft]$-module of rank $N$. Furthermore, restriction to $X^T$ is injective, and the Atiyah--Bott localization 
formula
\[
\int_Xa=\sum_{p\in X^T}\frac{a|_p}{e_p}\qquad(a\in H^*_T(X))
\]
holds (see \cite{AtiyahBott1984,BerlineVergne1982}). Note that $e_p\neq 0$ thanks to condition (G1). 

\subsection{The Atiyah--Bott--Lefschetz formula}\label{sec:ABL}

Let $M$ be a compact complex manifold with a holomorphic action of a
complex torus $T$, and let $E\to M$ be a 
$T$-equivariant holomorphic vector bundle. Then $T$ acts on the
cohomology groups $H^i(M,E)$ by $(t\cdot s)(x)=t\,s(t^{-1}x)$. The
\emph{equivariant holomorphic Euler characteristic} is 
\[
\chi_T(M,E):=\sum_{i\ge0}(-1)^i\operatorname{ch}_T H^i(M,E)\in
\operatorname{ch}_T\big(K^0_T({\rm pt})\big)=\mathbb Z\big[e^\mu\ |\ \mu\in\ft^*_\Z\big].
\]
For a fixed point $p\in M^T$, let $w_1(p),\dots,w_n(p)$ be the weights of $T_pM$.

\begin{theorem}[Atiyah--Bott--Lefschetz]\label{thm:ABL}
If $M^T$ is finite, then in the fraction field of $\mathbb Z[e^\mu\ |\ \mu\in\ft^*_\Z]$,
\[
\chi_T(M,E)=\sum_{p\in M^T}\frac{\operatorname{ch}_T(E_p)}{\prod_{j=1}^n\big(1-e^{-w_j(p)}\big)} .
\]
\end{theorem}

This is the holomorphic Lefschetz formula of Atiyah and Bott \cite{AtiyahBott1968}, applied to the Dolbeault
complex of $E$.

\subsection{Algebraic torus actions on \texorpdfstring{$\PP^1$}{P1}}
Let us first recall the following well known result about algebraic tori actions.
\begin{proposition}\label{prop:fixed-locus}
Suppose that $X$ is a smooth projective variety equipped with an algebraic $T$-action and let 
$K\subseteq T$ be an algebraic subtorus. 
Then the fixed locus $X^K$ is a smooth submanifold of $X$ and the tangent space $T_pX^K$ at a point $p\in X^K$ is the sum of the weight subspaces $(T_pX)_\mu$ of $T_pX$ such that $\mu|_{\mathrm{Lie}(K)}=0$. 
\end{proposition}
\begin{proof}
Indeed, let $K_c\subset K$ be the maximal compact subgroup; if $K\cong(\C^*)^r$, then $K_c=(S^1)^r$. We have $X^K=X^{K_c}$:
for $x\in X^{K_c}$ the orbit map $K\to X$, $k\mapsto k\cdot x$, is algebraic and constant on $K_c$, which is Zariski
dense in $K$. Since $K_c$ is compact, its action can be linearized by averaging (Bochner). For every $p\in X^{K_c}$
there are a $K_c$-invariant neighbourhood $U$ of $p$ and a $K_c$-equivariant biholomorphism from $U$ onto an open
neighbourhood of $0$ in $T_pX$, where $k\in K_c$ acts on $T_pX$ by $dk_p$. Under this map, $X^{K_c}\cap U$ corresponds to an open
subset of the linear subspace $(T_pX)^{K_c}$. Hence $X^{K_c}$ is a complex submanifold with $T_p(X^{K_c})=(T_pX)^{K_c}$. Since it
is also a closed subvariety, it is smooth as a variety. Finally, $(T_pX)^{K_c}$ is the sum of the $T$-weight spaces $(T_pX)_\mu$, $\mu=d\chi_\mu$, with $\chi_\mu$ trivial on $K_c$. By Zariski
density these are the $\chi_\mu$ trivial on $K$, and since $K$ is connected these are the $\mu$ with $\mu|_{\mathrm{Lie}(K)}=0$.
\end{proof}

\begin{proposition}\label{prop:P1actions}
Let $T$ be an algebraic torus.
\begin{enumerate}
\item[(a)] For distinct points $p,p'\in\PP^1$ and a non-trivial character $\chi\in X^*(T)$, choose homogeneous coordinates
with $p=[1,0]$ and $p'=[0,1]$. Then
\[
t\cdot[z_0,z_1]=[z_0,\chi(t)z_1]
\]
defines a non-trivial algebraic action of $T$ on $\PP^1$.
Its fixed points are $p$ and $p'$, and the tangent weights are $d\chi$ at $p$ and $-d\chi$ at $p'$.
\item[(b)] Every non-trivial algebraic action of $T$ on $\PP^1$ arises as in (a). The triple $(p,p',\chi)$ is determined by
the action up to the relation $(p,p',\chi)\sim(p',p,\chi^{-1})$.
\end{enumerate}
Thus non-trivial algebraic $T$-actions on $\PP^1$ correspond bijectively to triples $(p,p',\chi)$ with $p\ne p'$ and
$\chi\ne1$, modulo $(p,p',\chi)\sim(p',p,\chi^{-1})$. Two such actions are $T$-equivariantly isomorphic if and only if their
characters agree up to inversion.
\end{proposition}

\begin{proof}
(a) The action is the composite of $\chi:T\to\C^*$ with the action $s\cdot[z_0,z_1]=[z_0,sz_1]$ of $\C^*$, so it is algebraic. It
is non-trivial since $\chi\ne1$. Homogeneous coordinates with $p=[1,0]$ and $p'=[0,1]$ are unique up to
$(z_0,z_1)\mapsto(az_0,bz_1)$, $a,b\in\C^*$, and these rescalings commute with the action. A point $[z_0,z_1]$ is fixed if and only
if $z_0z_1(\chi(t)-1)=0$ for all $t$, i.e.\ if and only if it is $p$ or $p'$, since $\chi\ne1$.

The chart $z\mapsto[1,z]$ is $T$-equivariant for the action of $T$ on $\C$ through $\chi$. The chart $w\mapsto[w,1]$ is
$T$-equivariant through $\chi^{-1}$, since $t\cdot[w,1]=[\chi(t)^{-1}w,1]$. The differential at $0$ of a $T$-equivariant chart is an
isomorphism of $T$-modules. Hence $T_p\PP^1=\C_\chi$ and $T_{p'}\PP^1=\C_{\chi^{-1}}$, with weights $d\chi$ and $-d\chi$.

(b) An algebraic action of $T$ on $\PP^1$ is a homomorphism of algebraic groups $\rho:T\to\operatorname{Aut}(\PP^1)=\mathrm{PGL}(2,\C)$.
Its image is a torus, and it is non-trivial, so it is one-dimensional, hence a maximal torus of $\mathrm{PGL}(2,\C)$.
Since all maximal tori are conjugate, there are homogeneous coordinates in which $\rho(T)$ is contained in the diagonal torus
\[
D=\{[\operatorname{diag}(a,b)]\}.
\]
The map $D\to\C^*$, $[\operatorname{diag}(a,b)]\mapsto b/a$, is an isomorphism, and composing it with $\rho$ gives a non-trivial character
$\chi$. In these coordinates $t\cdot[z_0,z_1]=[z_0,\chi(t)z_1]$, so the action arises as in (a) from
$(p,p',\chi)=([1,0],[0,1],\chi)$.

By (a), the action determines the unordered pair $\{p,p'\}$ as its fixed locus, and it determines $\chi$ as the tangent
weight at $p$. Interchanging $p$ and $p'$ replaces $\chi$ by $\chi^{-1}$. Conversely, $(p,p',\chi)$ and $(p',p,\chi^{-1})$ define the same
action, because exchanging $z_0$ and $z_1$ turns $[z_0,\chi(t)z_1]$ into $[\chi(t)z_0,z_1]=[z_0,\chi(t)^{-1}z_1]$.

For the last statement, note that the tangent weights $\{d\chi,-d\chi\}$ are invariants of the equivariant isomorphism
class. Conversely, given the same $\chi$ for two actions, a projective transformation sending the first triple
$(p_1,p_1',\chi)$ to the second $(p_2,p_2',\chi)$ is $T$-equivariant.
\end{proof}

\subsection{Equivariant vector bundles on \texorpdfstring{$\PP^1$}{P1}}

Our next goal is to classify the $T$-equivariant vector bundles on
$\PP^1$ (see \cite{Kumar2003}). We will give a proof following an
argument of Bolibruch.  Let $\chi\in X^*(T)$ be a non-trivial
character and let $T$ 
act on $\PP^1$ by $t\cdot[z_0,z_1]=[z_0,\chi(t)z_1]$, with fixed points
$[1,0]$ and $[0,1]$. Let $E$ be a 
$T$-equivariant vector bundle of rank $r$ on $\PP^1$, i.e.\ $T$ acts algebraically on the total space of $E$, linearly on the
fibres and compatibly with the projection. If $E$ is given as a holomorphic bundle, it is algebraic by GAGA.

Let $U_0=\{z_0\ne0\}$ and $U_1=\{z_1\ne0\}$ be the standard affine charts, with coordinates $z=z_1/z_0$ on $U_0$ and $w=z_0/z_1$ on
$U_1$. Then $t\cdot z=\chi(t)z$ and $t\cdot w=\chi(t)^{-1}w$. The module $\Gamma(U_0,E)$ is a finitely generated projective module over the
principal ideal domain $\C[z]$, hence free; similarly $\Gamma(U_1,E)$ is free over $\C[w]$. Both charts are $T$-invariant, so each
restriction $E|_{U_a}$ is a $T$-equivariant vector bundle on $U_a\cong\mathbb A^1$. We treat $U_0$; the case of $U_1$ is the same, with
$\chi$ replaced by $\chi^{-1}$.

\paragraph{The cocycle of a frame.} Let $e=(e_1,\dots,e_r)$ be a frame of $E|_{U_0}$, written as a row vector. For $t\in T$, the
map $t:E_z\to E_{t\cdot z}$ is linear, so
\begin{equation}\label{eq:cocycle-def}
t\cdot e(z)=e(t\cdot z)\,A(t,z)
\end{equation}
for a unique matrix $A(t,z)\in\mathrm{GL}_r$. Since the action is algebraic, $A\in\mathrm{GL}_r(\C[T]\otimes\C[z])$. With $\C[T]=\C[\ft^*_{\mathbb Z}]$
and $t^m:=\chi_m(t)$ for $m=d\chi_m\in\ft^*_{\mathbb Z}$, this means
\[
A(t,z)=\sum_{m\in\ft^*_{\mathbb Z}}a_m(z)\,t^m,\qquad a_m\in\mathrm{Mat}_r(\C[z]),
\]
a finite sum. Applying \eqref{eq:cocycle-def} twice, and using that $t_2$ acts linearly on the fibres, gives
\[
(t_2t_1)\cdot e(z)=t_2\cdot\big(e(t_1\cdot z)A(t_1,z)\big)=e(t_2t_1\cdot z)\,A(t_2,t_1\cdot z)\,A(t_1,z).
\]
Hence $A$ satisfies the cocycle conditions
\begin{equation}\label{eq:cocycle}
A(t_2t_1,z)=A(t_2,t_1\cdot z)\,A(t_1,z),\qquad A(1,z)=1.
\end{equation}
If $e'=e\,g$ with $g\in\mathrm{GL}_r(\C[z])$ is another frame, then
\begin{equation}\label{eq:gauge}
A'(t,z)=g(t\cdot z)^{-1}A(t,z)\,g(z).
\end{equation}

\paragraph{Cocycles and gradings.} The group $T$ acts on $P:=\Gamma(U_0,E)$ by $(t\cdot s)(z):=t\cdot s(t^{-1}\cdot z)$, and on $\C[z]$ by
$(t\cdot\varphi)(z):=\varphi(t^{-1}\cdot z)$. These actions satisfy $t\cdot(\varphi s)=(t\cdot\varphi)(t\cdot s)$. The coordinate $z$ is homogeneous of
weight $-d\chi$, so $\C[z]$ is graded by $\ft^*_{\mathbb Z}$ with $\deg z^n=-n\,d\chi$. Writing $s=e\,c$ with $c\in\C[z]^r$, we get
\[
t\cdot s=e\,A(t,t^{-1}\cdot z)\,c(t^{-1}\cdot z).
\]
The right-hand side is a finite sum $\sum_\mu t^\mu\,e\,\pi_\mu(c)$. This defines operators $\pi_\mu$ on $P$, and we put
$P_\mu:=\pi_\mu(P)$.

\begin{lemma}\label{lem:cocycle-grading}
\begin{enumerate}
\item[(a)] For every frame $e$, the operators $\pi_\mu$ are orthogonal idempotents with $\sum_\mu\pi_\mu=\mathrm{id}$, and $z\,P_\mu\subseteq P_{\mu-d\chi}$.
Thus $P=\bigoplus_\mu P_\mu$ is a $\ft^*_{\mathbb Z}$-graded $\C[z]$-module, and $P_\mu=\{s\in P\mid t\cdot s=t^\mu s\ \text{for all }t\in T\}$.
\item[(b)] There is a frame $e$ consisting of homogeneous elements, $e_j\in P_{\sigma_j}$. In such a frame
\[
A(t,z)=D(t):=\operatorname{diag}\big(t^{\sigma_1},\dots,t^{\sigma_r}\big).
\]
\item[(c)] Consequently every cocycle has the form
\[
A(t,z)=g(\chi(t)z)\,D(t)\,g(z)^{-1},\qquad g\in\mathrm{GL}_r(\C[z]),
\]
where the multiset $\{\sigma_1,\dots,\sigma_r\}$ is the multiset of weights of the $T$-module $E_0$, the fibre at $z=0$. In
particular $E|_{U_0}\cong U_0\times E_0$ equivariantly, with $T$ acting by $t\cdot(z,v)=(\chi(t)z,\,t\cdot v)$.
\end{enumerate}
\end{lemma}

\begin{proof}
(a) Expanding $(t_2t_1)\cdot s=t_2\cdot(t_1\cdot s)$ in $t_1,t_2$ gives $\sum_\mu(t_2t_1)^\mu\pi_\mu=\sum_{\mu,\nu}t_2^\nu t_1^\mu\,\pi_\nu\pi_\mu$, so
$\pi_\nu\pi_\mu=\delta_{\mu\nu}\pi_\mu$, since distinct characters are linearly independent. Setting $t=1$ gives $\sum_\mu\pi_\mu=\mathrm{id}$.
For $s\in P_\mu$ we have $t\cdot(zs)=(t\cdot z)(t\cdot s)=t^{-d\chi}z\,t^\mu s$, so $zs\in P_{\mu-d\chi}$. The description of $P_\mu$ follows from
the expansion of $t\cdot s$.

(b) \emph{Graded Nakayama.}
We claim that if $N\subseteq P$ is a graded submodule such that
$N+zP=P$, then $N=P$. 
Choose a linear form $\ell:\ft^*_{\mathbb Z}\to\mathbb R$ with
$\ell(d\chi)<0$, so that $\ell(\deg z)>0$. Since $P$ is generated by
finitely many homogeneous elements, of degree say $\mu_i$, the linear
form $\ell$ takes values in finitely many arithmetic progressions
$\ell(\mu_i)+\ell(\deg z)\,\mathbb Z_{\ge0}$ on the weights of $P$. In
particular, these values are bounded below. Now the claim can be
proved as follows. Put $M:=P/N$. Then $M$ is
graded and $M=zM$. If $M\ne0$, choose a 
homogeneous $0\ne m\in M$ with $\ell(\deg m)$ minimal. Since $m\in
zM$, $m$ is a sum of elements of the form $zm'$ with $m'$ 
homogeneous and $\ell(\deg m')=\ell(\deg m)-\ell(\deg z)<\ell(\deg
m)$. By minimality each such $m'$ is $0$, so $m=0$, a 
contradiction. Hence $N=P$.

\emph{The frame.} The fibre $E_0=P/zP$ is graded, with $(E_0)_\mu$ the $\mu$-weight space of the $T$-module $E_0$. Choose a basis
of $E_0$ of weight vectors, of weights $\sigma_1,\dots,\sigma_r$, and lift them to homogeneous $e_j\in P_{\sigma_j}$. The submodule $N$
generated by the $e_j$ satisfies $N+zP=P$, so $N=P$. 
Thus the map $\C[z]^r\to P$, $c\mapsto e\,c$, is surjective, and we have a short exact sequence of $\C[z]$-modules
\[
0\longrightarrow Q\longrightarrow\C[z]^r\xrightarrow{\;c\,\mapsto\,e\,c\;}P\longrightarrow0,\qquad Q:=\{c\in\C[z]^r\mid e\,c=0\}.
\]
Since $\C(z)$ is a localization of $\C[z]$, it is flat (Atiyah--Macdonald, Prop.~3.3), so
\[
0\longrightarrow Q\otimes\C(z)\longrightarrow\C(z)^r\longrightarrow P\otimes\C(z)\longrightarrow0
\]
is exact. Since $P$ is free of rank $r$, the last two terms are $r$-dimensional vector spaces over $\C(z)$. So the
surjection between them is an isomorphism, and $Q\otimes\C(z)=0$. Thus every $q\in Q$ satisfies $sq=0$ for some $0\ne s\in\C[z]$. Since
$Q\subseteq\C[z]^r$ is torsion-free, $Q=0$. Hence $e$ is a frame.

Finally, homogeneity $t\cdot e_j=t^{\sigma_j}e_j$ means $t\cdot e_j(t^{-1}\cdot z)=t^{\sigma_j}e_j(z)$, i.e.\ $t\cdot e_j(z)=e_j(t\cdot z)\,t^{\sigma_j}$. Comparing with
\eqref{eq:cocycle-def} gives $A=D$.

(c) Let $e$ be any frame and $e'=e\,g$ a homogeneous frame as in (b). By \eqref{eq:gauge}, $D(t)=g(t\cdot z)^{-1}A(t,z)\,g(z)$. At $z=0$
we get $A(t,0)=g(0)D(t)g(0)^{-1}$, and $A(t,0)$ is the matrix of the action of $T$ on $E_0$ in the basis $e(0)$. So the
$\sigma_j$ are the weights of $E_0$, and $D$ is unique up to permutation. In the homogeneous frame, $e'(z)c\mapsto(z,e'(0)c)$ is the
asserted equivariant isomorphism $E|_{U_0}\cong U_0\times E_0$.
\end{proof}

Over $U_1$ the same lemma applies with $z$, $\chi$ and $E_0$ replaced by $w$, $\chi^{-1}$ and $E_\infty$, the fibre at $[0,1]$. Thus $E$ has
homogeneous frames $e^0$ over $U_0$ and $e^1$ over $U_1$, with $T$ acting on them by $D_0(t)$ and $D_\infty(t)$, the diagonal
matrices of the weights of $E_0$ and $E_\infty$.

\paragraph{Homogeneous matrices.} For $\alpha=(\alpha_1,\dots,\alpha_r)\in(\ft^*_{\mathbb Z})^r$ put
\[
t^\alpha:=\operatorname{diag}(t^{\alpha_1},\dots,t^{\alpha_r}).
\]
A matrix
$F\in\mathrm{Mat}_r(\C[z,z^{-1}])$ is \emph{homogeneous of type $(\alpha,\beta)$} if
\[
F(\chi(t)z)=t^\alpha F(z)\,t^{-\beta}\qquad(t\in T).
\]
Writing $F_{ij}=\sum_nc_nz^n$, this says $c_n\big(t^{n\,d\chi}-t^{\alpha_i-\beta_j}\big)=0$ for all $n$. Since $d\chi\ne0$, each entry $F_{ij}$ is either $0$ or a
monomial $c\,z^n$ with $n\,d\chi=\alpha_i-\beta_j$. The following properties are immediate:
\begin{itemize}
\item a product of matrices of types $(\alpha,\beta)$ and $(\beta,\gamma)$ has type $(\alpha,\gamma)$;
\item the inverse of an invertible matrix of type $(\alpha,\beta)$ has type $(\beta,\alpha)$;
\item for $K=\operatorname{diag}(k_1,\dots,k_r)\in\mathrm{Mat}_r(\mathbb Z)$ the matrix $z^K$ has type $(\alpha,\alpha-K\,d\chi)$ for every $\alpha$, where
$(\alpha-K\,d\chi)_j=\alpha_j-k_j\,d\chi$;
\item a permutation matrix $P$ has type $(\alpha,\alpha P)$, where $(\alpha P)_j=\alpha_{\pi(j)}$ if $Pe_j=e_{\pi(j)}$.
\end{itemize}
A constant matrix of type $(\alpha,\beta)$ has non-zero entries only where $\alpha_i=\beta_j$. In particular, if $F\in\mathrm{Mat}_r(\C[z])$ has type
$(\alpha,\beta)$, then the $j$-th column of $F(0)$ lies in the span of the $e_i$ with $\alpha_i=\beta_j$.

\paragraph{The transition matrix.} Let $E$ be a $T$-equivariant vector bundle of rank $r$ on $\PP^1$. By
Lemma~\ref{lem:cocycle-grading}, there are homogeneous frames $e^0$ over $U_0$ and $e^1$ over $U_1$, with
\[
t\cdot e^0(x)=e^0(t\cdot x)\,t^{\sigma},\qquad t\cdot e^1(x)=e^1(t\cdot x)\,t^{\tau},
\]
where $\sigma$ and $\tau$ list the weights of $E_0$ and $E_\infty$. On $U_0\cap U_1$ write $e^1=e^0M$ with $M\in\mathrm{GL}_r(\C[z,z^{-1}])$. Then
\[
t\cdot e^1(x)=t\cdot e^0(x)\,M(z(x))=e^0(t\cdot x)\,t^\sigma M(z(x)).
\]
On the other hand,
\[
t\cdot e^1(x)=e^1(t\cdot x)\,t^\tau=e^0(t\cdot x)\,M(\chi(t)z(x))\,t^\tau .
\]
Hence $M$ is homogeneous of type $(\sigma,\tau)$. Replacing $e^0$ by $e^0g_0$ and $e^1$ by $e^1g_1$, with $g_0\in\mathrm{GL}_r(\C[z])$ of type
$(\sigma,\sigma')$ and $g_1\in\mathrm{GL}_r(\C[z^{-1}])$ of type $(\tau,\tau')$, gives homogeneous frames again, and replaces $M$ by $g_0^{-1}Mg_1$.

\paragraph{The line bundles $L(k,\sigma)$.} For $k\in\mathbb Z$ and $\sigma\in\ft^*_{\mathbb Z}$ let 
\[
L(k,\sigma):=((\C^2\setminus{\{0\}} )\times \C)/\C^*\to \mathbb P^1,
\]
where $\C^*$ acts on $\C$ via the character $t\mapsto t^k$, that is, $L(k,\sigma)$ is the line bundle $\mathcal O(k)$. Let us define frames
\[
e^0([z_0,z_1])=[1,z_1/z_0,1],\qquad 
e^1(w)=[z_0/z_1,1,1],\qquad [z_0,z_1]\in \mathbb P^1.
\]
We have 
\[
e^1([z_0,z_1])=[z_0/z_1,1,1]=[1/z,1,1]=[1,z,z^k]=e^0([z_0,z_1]) z^k,
\]
where $z=z_1/z_0$. The transition function is $z^k.$ Let us define a $T$-action on $L(k,\sigma)$ such that $e^0$ is homogeneous of weight $\sigma$, i.e.\
$t\cdot e^0(x)=e^0(t\cdot x)\,t^\sigma$. Then
\[
t\cdot e^1(x)=z(x)^k\,e^0(t\cdot x)\,t^\sigma=z(x)^kz(t\cdot x)^{-k}\,e^1(t\cdot x)\,t^\sigma=e^1(t\cdot x)\,t^{\sigma-k\,d\chi}.
\]
So this action is algebraic on all of $\mathcal O(k)$, and the weights of $L(k,\sigma)$ at $[1,0]$ and $[0,1]$ are $\sigma$ and $\sigma-k\,d\chi$.

\begin{proposition}[Equivariant Birkhoff--Grothendieck theorem]\label{prop:eqBG}
Let $T$ act on $\PP^1$ by $t\cdot[z_0,z_1]=[z_0,\chi(t)z_1]$ with $\chi\ne1$, and let $E$ be a $T$-equivariant vector bundle of rank $r$ on $\PP^1$.
Then there are $(k_j,\sigma_j)\in\mathbb Z\times\ft^*_{\mathbb Z}$, $j=1,\dots,r$, and a $T$-equivariant isomorphism
\[
E\cong\bigoplus_{j=1}^rL(k_j,\sigma_j).
\]
The multiset $\{(k_j,\sigma_j)\}$ is uniquely determined by $E$.
\end{proposition}

\begin{proof}
\emph{Existence.} We follow the proof of the Birkhoff factorization in \cite[Prop.~3.5]{Milanov2017} (originally due to Bolibruch)
keeping track of homogeneity. Let $M$ be the transition matrix of homogeneous frames, of type $(\sigma,\tau)$. We look for a
factorization
\begin{equation}\label{eq:BF}
M=L\,z^K\,W
\end{equation}
with $K=\operatorname{diag}(k_1,\dots,k_r)$, $k_1\ge\dots\ge k_r$, and:
\begin{itemize}
\item $L\in\mathrm{Mat}_r(\C[z])$ homogeneous of type $(\sigma,w)$, for some $w\in(\ft^*_{\mathbb Z})^r$, with $\det L=c\,z^v$, $c\ne0$, $v\ge0$;
\item $W\in\mathrm{GL}_r(\C[z^{-1}])$ homogeneous of type $(w-K\,d\chi,\tau)$.
\end{itemize}
Such a factorization exists. Indeed, $\det M$ is a unit of $\C[z,z^{-1}]$, i.e.\ $\det M=c\,z^m$, and for $N\gg0$ we may take
$L=z^NM$ (of type $(\sigma,\tau-N\,d\chi)$), $K=-N\cdot1$ and $W=1$.

We show that if $v>0$, there is another factorization \eqref{eq:BF} of the same kind with $v$ replaced by $v-1$. Since $\det L(0)=0$,
the columns $m_1,\dots,m_r$ of $L(0)$ are linearly dependent. Let $i<r$ be maximal such that $m_1,\dots,m_i$ are linearly
independent. Then
\[
m_{i+1}=\sum_{a\le i}s_am_a,\qquad s_a\in\C .
\]
Since $L$ has type $(\sigma,w)$, we have $t^\sigma m_a = m_a t^{w_a}$ for all $1\leq a\leq i+1$. Therefore, multiplying the above identity by $t^\sigma$ yields
\[
m_{i+1} t^{w_{i+1}}=\sum_{a\le i}s_am_a \, t^{w_a}.
\]
Since the characters of $T$ are linearly independent functions, by comparing the coefficients in front of $t^{w_{i+1}}$, we get that $m_{i+1}$ is a linear combination only of those $m_a$ for which $w_a=w_{i+1}$.  Therefore, we may assume that $s_a=0$ unless $w_a=w_{i+1}$.

Put $S:=\sum_{a\le i}s_aE_{a,i+1}$. Then $S^2=0$, and $S$ has type $(w,w)$. The matrix $L':=L(1-S)$ has type $(\sigma,w)$, and its $(i+1)$-st column
vanishes at $z=0$. Hence $L'':=L'z^{-E_{i+1,i+1}}\in\mathrm{Mat}_r(\C[z])$, of type $(\sigma,w+d\chi\,e_{i+1})$, and $\det L''=c\,z^{v-1}$. Since $(1-S)^{-1}=1+S$ and
$z^{-K}E_{a,i+1}z^K=z^{k_{i+1}-k_a}E_{a,i+1}$, we get
\[
M=L'(1+S)z^KW=L''\,z^{K+E_{i+1,i+1}}\,W',
\]
where
\[
W':=\Big(1+\sum_{a\le i}s_a\,z^{k_{i+1}-k_a}E_{a,i+1}\Big)W .
\]
Here $k_{i+1}-k_a\le0$ for $a\le i$, so $W'\in\mathrm{GL}_r(\C[z^{-1}])$. The factor in front of $W$ has type $(w-K\,d\chi,w-K\,d\chi)$, because
$s_a\ne0$ only if $w_a=w_{i+1}$. So $W'$ has type $(w-K\,d\chi,\tau)$. With $w':=w+d\chi\,e_{i+1}$ and $K':=K+E_{i+1,i+1}$ we have
$w'-K'd\chi=w-K\,d\chi$. Finally, conjugating by a permutation matrix $P$, replacing $(L'',z^{K'},W')$ by $(L''P,P^{-1}z^{K'}P,P^{-1}W')$,
restores the order $k_1\ge\dots\ge k_r$ and preserves all homogeneity conditions.

After finitely many steps $v=0$. Then $\det L$ is a non-zero constant, so $L\in\mathrm{GL}_r(\C[z])$. Hence
\[
f^0:=e^0L,\qquad f^1:=e^1W^{-1}=e^0MW^{-1}=e^0L\,z^K=f^0z^K
\]
are homogeneous frames over $U_0$ and $U_1$. The weights of $f^0$ are $w$, those of $f^1$ are $w-K\,d\chi$, and $f^1_j=z^{k_j}f^0_j$. So the $j$-th frame vectors span a
$T$-equivariant line subbundle isomorphic to $L(k_j,w_j)$, and $E\cong\bigoplus_jL(k_j,w_j)$.

\emph{Uniqueness.} Suppose $E\cong\bigoplus_jL(k_j,\sigma_j)\cong\bigoplus_jL(k'_j,\sigma'_j)$, with the summands ordered so that $k_1\ge\dots\ge k_r$ and
$k'_1\ge\dots\ge k'_r$. The frames $z_0^{k_j}$ and $z_1^{k_j}$ of the summands give homogeneous frames $f^0$ over $U_0$ and $f^1$ over $U_1$, with
\[
f^1=f^0z^K,
\]
where $K=\operatorname{diag}(k_1,\dots,k_r)$ and $f^0$ has weights $\sigma=(\sigma_1,\dots,\sigma_r)$. Similarly the second decomposition gives
$f'^0$ and $f'^1=f'^0z^{K'}$, where $f'^0$ has weights $\sigma'$. Write $f'^0=f^0G$ and $f'^1=f^1H$ with $G\in\mathrm{GL}_r(\C[z])$ and $H\in\mathrm{GL}_r(\C[z^{-1}])$.
As shown above, $G$ is homogeneous of type $(\sigma,\sigma')$. From $f'^0z^{K'}=f^1H=f^0z^KH$ we get
\[
G=z^K\,H\,z^{-K'},\qquad\text{i.e.}\qquad G_{ij}(z)=z^{k_i-k'_j}H_{ij}(z^{-1}).
\]
The left-hand side is a polynomial in $z$, and the right-hand side involves only powers $z^n$ with $n\le k_i-k'_j$. Hence
\begin{equation}\label{eq:Gvanish}
G_{ij}=0\quad\text{if }k_i<k'_j,\qquad G_{ij}\in\C\quad\text{if }k_i=k'_j .
\end{equation}

We will first argue that $K=K'$. Suppose not, and let $s\ge0$ be such that $k_i=k'_i$ for $i\le s$ and $k_{s+1}\ne k'_{s+1}$. Exchanging the two
decompositions, which replaces $G$ by $G^{-1}=z^{K'}H^{-1}z^{-K}$, we may assume $k_{s+1}<k'_{s+1}$. For $i\ge s+1$ and $j\le s+1$ we have
\[
k_i\le k_{s+1}<k'_{s+1}\le k'_j ,
\]
so $G_{ij}=0$ by \eqref{eq:Gvanish}. Thus the first $s+1$ columns of $G$ lie in the span of $e_1,\dots,e_s$. They are therefore
linearly dependent over $\C(z)$, so $\det G=0$, contradicting $G\in\mathrm{GL}_r(\C[z])$. Hence $K=K'$.

It remains only to prove that the two sets $\{(k_j,\sigma_j)\}_{1\leq
  j\leq r}$ and $\{(k'_j,\sigma'_j)\}_{1\leq
  j\leq r}$ coincide. For $n\in\mathbb Z$ let $I_n:=\{j\mid k_j=n\}$. Since $K=K'$, \eqref{eq:Gvanish} shows that $G$ is block upper triangular
with respect to the blocks $I_n$, ordered by decreasing $n$, and that its diagonal blocks $G_n:=(G_{ij})_{i,j\in I_n}$ are constant.
Since $\det G=\prod_n\det G_n$ is a non-zero constant, each $G_n$ is invertible. Homogeneity of type $(\sigma,\sigma')$ gives, for the
constant block $G_n$,
\[
t^{\sigma|_{I_n}}\,G_n=G_n\,t^{\sigma'|_{I_n}}\qquad(t\in T).
\]
So $G_n$ is an isomorphism between the $T$-modules with weights $\{\sigma'_j\}_{j\in I_n}$ and $\{\sigma_i\}_{i\in I_n}$. Hence these multisets
agree for every $n$, and the multisets $\{(k_j,\sigma_j)\}$ and $\{(k'_j,\sigma'_j)\}$ coincide.

\end{proof}

\begin{corollary}\label{cor:ABLP1} 
Suppose that $T$ acts on $\PP^1$ with tangent weight $\epsilon$ at
$p:=[1,0]$ and $-\epsilon$ at $p'=[0,1]$. Let $E$ be a $T$-equivariant
vector bundle on $\PP^1$ and $E\cong \bigoplus L(a_i,u_i)$ be the
$T$-equivariant decomposition into line bundles. Then the weights of the $T$-modules
$H^0(\PP^1,E)$ and $H^1(\PP^1,E) $ are given
by 
\begin{align*}
  \wts \ H^0(\PP^1,E) &  = \bigcup_{i:a_i\geq 0}
  \{u_i-m\epsilon \ |\  0\leq m\leq a_i\}, \\
  \wts \ H^1(\PP^1,E) & = \bigcup_{i:a_i\leq -2}
  \{u_i-m\epsilon \ | \ a_i+1\leq m\leq -1\}.
\end{align*}
\end{corollary}
\begin{proof}
It is sufficient to prove the corollary only in the case when 
$E=L(a,u)$ is a $T$-equivariant line bundle on $\PP^1$. 
By the holomorphic Lefschetz formula \cite{AtiyahBott1968},
\[
\chi_T(\PP^1,E)=\frac{e^u}{1-e^{-\epsilon}}+\frac{e^{u-a\epsilon}}{1-e^{\epsilon}}
=\begin{cases}\sum_{m=0}^{a}e^{u-m\epsilon}, & a\ge0,\\[2pt] -\sum_{m=a+1}^{-1}e^{u-m\epsilon}, & a<0.\end{cases}
\]
Therefore, $H^0(\PP^1, E)$ has the weights $u-m\epsilon$, $0\le m\le
a$, and $H^1(\PP^1,E)$ the weights $u-m\epsilon$, $a<m<0$. 
\end{proof}

\subsection{\texorpdfstring{$T$}{T}-curves and edge factors}
The $T$-curves in a GKM-manifold can be classified as follows.
\begin{proposition}[$T$-curves]\label{prop:gkmcurves}
Let $X$ be a GKM manifold and $p\in X^T$.
\begin{enumerate}
\item[(a)] For each weight $\chi$ of $T_pX$ there is a unique irreducible $T$-invariant curve $C\ni p$ with $T_pC$ the
$\chi$-weight line. It is smooth, $C\cong\PP^1$, and $C^T=\{p,p'\}$ with $p'\ne p$.
\item[(b)] Every irreducible $T$-invariant curve through $p$ arises in this way; there are finitely many.
\item[(c)] If $\varepsilon_C(p)$ denotes the weight of $T_pC$, then $\varepsilon_C(p')=-\varepsilon_C(p)$.
\end{enumerate}
\end{proposition}

\begin{proof}
(a) Since $p$ is an isolated fixed point, recalling Proposition \ref{prop:fixed-locus}, we get $(T_pX)^T=T_p(X^T)=0$. Therefore, since $0$ is not a weight of the $T$-module $T_pX$, we have $d\chi\ne0$, that is $\chi\ne1$. Let $K:=(\ker\chi)^0$, a subtorus of codimension one, and let $Y$ be the connected component through $p$ of $X^K$.
A weight of $T_pX$ vanishing on $\mathrm{Lie}(K)=\ker d\chi$ must be proportional to $d\chi$ and hence equal to $d\chi$ thanks to condition (G2). Hence $T_pY$ coincides with the one-dimensional $\chi$-weight subspace of $T_pX$. In particular, $Y$ is a smooth connected curve. It is projective, being closed in $X$. It is $T$-invariant: $T$ commutes with $K$, so it preserves $X^K$, and since $T$ is connected it preserves the connected component $Y$ of $X^K$.

The action of $T$ on $Y$ is non-trivial, since otherwise $Y\subseteq X^T$, contradicting (G1). Choose $y\in Y$ that is not
fixed, and a cocharacter $\lambda:\C^*\to T$ with $\langle\chi,\lambda\rangle\ne0$. Since $\lambda(\C^*)$ is a one-dimensional subtorus not
contained in $K$, we have $T=K\cdot\lambda(\C^*)$, and $K$ acts trivially on $Y$. Hence the orbit map
\[
f:\C^*\to Y,\qquad s\mapsto\lambda(s)\cdot y,
\]
is non-constant. It is a morphism of algebraic varieties, being the composite of $s\mapsto(\lambda(s),y)$ with the action
morphism $T\times X\to X$, which is algebraic by assumption. Recalling the valuative criterion, we get that $f$ extends to a morphism $\PP^1\to Y$. Alternatively, the extension can be constructed also as follows. Embed $Y\subseteq\PP^N$, and choose the homogeneous coordinates so that
$z_0(f(1))\ne0$. Then $f^*(z_i/z_0)$, $1\le i\le N$, are rational functions of $s$. Write $f^*(z_i/z_0)=f_i(s)/f_0(s)$ with
$f_0,\dots,f_N\in\C[s]$ without common factor. They have no common zero in $\C$, so $[f_0(s):\dots:f_N(s)]$ is a morphism
$\C\to\PP^N$. It agrees with $f$ on a dense open subset of $\C^*$, hence extends $f$ across $s=0$. The same argument in the
coordinate $s^{-1}$ applies at $s=\infty$. The image of the extension lies in $Y$, since $Y$ is closed. We get a non-constant
morphism $f:\PP^1\to Y$. Since $Y$ is an irreducible curve, $f$ is surjective, so the function field of $Y$ embeds into
$\C(s)$. By L\"uroth's theorem it is a purely transcendental extension of $\C$, i.e. $Y$ is a smooth projective rational
curve, and hence $Y\cong\PP^1$.

By Proposition~\ref{prop:P1actions}(b), the action of $T$ on $Y$ is given by a triple $(p_0,p_0',\chi')$. Its fixed points are $p_0$ and
$p_0'$. One of them is $p$, so after replacing the triple by $(p_0',p_0,(\chi')^{-1})$ if necessary, we may assume $p_0=p$. Put
$p':=p_0'$. Then $Y^T=\{p,p'\}$ with $p'\ne p$. By Proposition~\ref{prop:P1actions}(a), the tangent weight of $Y$ at $p$ is $d\chi'$. Since
$T_pY$ is the $\chi$-weight space, $\chi'=\chi$. Put $C:=Y$. Uniqueness of $C$ follows from (b).

(b) Let $C'\ni p$ be an irreducible closed $T$-invariant curve.

\emph{Step 1: the stabilizer.} By (G1), $T$ does not fix $C'$ pointwise. Choose $y\in C'$ with $T\cdot y\ne\{y\}$, and let
$T_y\subseteq T$ be its stabilizer. It is a closed complex Lie subgroup, and $T_y\ne T$. Put $S:=T/T_y$, a connected complex
manifold of dimension $m:=\dim T-\dim T_y$. If $m=0$, then $T_y$ would be open in $T$, hence equal to $T$, since $T$ is
connected. So $m\ge1$. The orbit map induces an injective holomorphic map
\[
a:S\to C',\qquad sT_y\mapsto s\cdot y .
\]

\emph{Step 2: $m=1$.} 
Let us recall the following result in dimension theory. If $X$ is a complex space and $x\in X$, then Grauert--Remmert defined $\dim_x X$ as the smallest number $k$ of functions $f_1,\dots,f_k\in \mathcal O_{X,x}$ such that $x$ is an isolated point for the zero locus $f_1=\dots=f_k=0$ in $X$ (see \cite{GR1984}, Ch.~5). It is known that $\dim_xX$ coincides with the Krull dimension of $\mathcal O_{X,x}$ and hence with $\dim X$ when $X$ is a complex manifold. 
Fix $t\in S$. Since $\dim_{t\cdot y}C'=1$, there is $g\in\mathcal O_{C',t\cdot y}$ such that $t\cdot y$ is an isolated point of the zero set of $g$. We claim that $t$ is an isolated zero of $a^*g\in\mathcal O_{S,t}$.
Otherwise, there are $t_n\to t$ with $t_n\ne t$ and $g(a(t_n))=0$. Then $a(t_n)\to a(t)$ by continuity, and $a(t_n)\ne a(t)$ by injectivity, which contradicts the choice of $g$. Hence $m=\dim S = \dim_tS\le1$ and therefore $m=1$, that is \ $\dim T_y=\dim T-1$.

\emph{Step 3: $T_y$ fixes $C'$.} Since $T$ is abelian, $T_{s\cdot y}=T_y$ for all $s\in T$. So $T_y$ fixes $T\cdot y$ pointwise. The
fixed locus
\[
(C')^{T_y}=\bigcap_{h\in T_y}\{x\in C'\mid h\cdot x=x\}
\]
is Zariski closed in $C'$. Indeed, each $h$ acts by an automorphism of the variety $X$, so each set in the intersection is the preimage of the diagonal of $X\times X$ under the algebraic morphism $C'\to X\times X$, $x\mapsto(x,h\cdot x)$. Algebraic maps are continuous in the Zariski topology and the diagonal is Zariski closed, so each set in the intersection is Zariski closed. Finally, since the Zariski topology is a topology, the intersection of any set of Zariski closed subsets is Zariski closed. On the other hand, a proper Zariski closed subset of the irreducible curve $C'$ must be finite while $(C')^{T_y}$ is infinite because it contains $T\cdot y=a(S)$ and the latter is infinite because $a$
is injective and $\dim S=1$. Therefore $(C')^{T_y}$ is not proper, i.e. it coincides with $C'$. We get that $T_y$ fixes $C'$ pointwise. Put $H:=(T_y)^0$, a connected closed subgroup of $T$ with $\dim H=\dim T-1$, which fixes $C'$ pointwise.

\emph{Step 4: the weight.} Identify $T=(\C^*)^n$, $\ft=\C^n$,  so that $T_c=(S^1)^n=\exp(i\mathbb R^n)$. Since $T_c$ is compact and $p$ is a $T_c$-fixed point, the $T_c$-action can be linearized locally near $p$, that is, there exists a $T_c$-equivariant neighborhood $U$ of $p$ and a holomorphic embedding $\psi:U\to T_pX$ such that $\psi(tx)=dt_p \psi(x)$ for all $x\in U$ and $t\in T_c$. 
We claim that $\psi$ is locally a $T$-equivariant linearization, that is $\psi(tx)=dt_p \psi(x)$ for all $(t,x)\in T\times U$ sufficiently close to $(1,p)$. To this end, choose a ball $B\subset\ft$
around $0$ and a neighbourhood $U'\subseteq U$ of $p$ with $\exp(B)\cdot U'\subseteq U$. For $x\in U'$ the holomorphic map
\[
F_x:B\to T_pX,\qquad F_x(\xi):=\psi(\exp(\xi)\cdot x)-d(\exp\xi)_p\,\psi(x),
\]
vanishes on $B\cap i\mathbb R^n$. Writing $\xi=\sigma+i\tau$ with $\sigma,\tau\in\mathbb R^n$, the Cauchy--Riemann equations give
$\partial^\alpha F_x/\partial\xi^\alpha=(-i)^{|\alpha|}\partial^\alpha F_x/\partial\tau^\alpha$ for every multi-index $\alpha$. The right-hand side
vanishes at $0$, since $F_x=0$ on $B\cap i\mathbb R^n$. Hence the Taylor series of $F_x$ at $0$ is zero, and $F_x\equiv0$ on the ball $B$.
Thus $\psi(h\cdot x)=dh_p\,\psi(x)$ for $h\in\exp(B)$ and $x\in U'$.

Since $C'$ is a curve, $p$ is not an isolated point of $C'$. Choose $x\in C'\cap U'$ with $x\ne p$, and put $v:=\psi(x)\ne0$. For
$h\in H\cap\exp(B)$ we have $h\cdot x=x$, hence $dh_p\,v=v$. The group $H$ is connected, so it is generated by $H\cap\exp(B)$.
Therefore $v\in(T_pX)^H$, and $(T_pX)^H\ne0$.

Write $T_pX=\bigoplus_\mu V_\mu$ as the sum of its weight spaces, where $h=\exp\xi$ acts on $V_\mu$ by $e^{\mu(\xi)}$. Since $H$ is
connected and generated by $\exp(\mathrm{Lie}(H))$, we get
\[
(T_pX)^H=\bigoplus_{\mu|_{\mathrm{Lie}(H)}=0}V_\mu .
\]
So some weight $\chi$ vanishes on $\mathrm{Lie}(H)$. We have $\chi\ne0$, because $(T_pX)^T=0$. Since $\mathrm{Lie}(H)$ is a hyperplane,
$\ker\chi=\mathrm{Lie}(H)$. A second weight with this property would have the same kernel, hence be proportional to $\chi$,
which is excluded by (G2). So $\chi$ is unique.

\emph{Step 5: $C'$ is the curve of (a).} Let $K:=(\ker\chi)^0$, as in (a). Then $\mathrm{Lie}(K)=\ker\chi=\mathrm{Lie}(H)$, and $K$ and $H$
are connected, so $K=H$. Hence $C'\subseteq X^K$. Since $C'$ is connected and contains $p$, it lies in the component $Y$ of $X^K$
through $p$. By (a), $Y$ is an irreducible curve. Since $C'\subseteq Y$ are both irreducible curves, $C'=Y$.
Thus every irreducible $T$-invariant curve through $p$ is the curve of (a) for a weight of $T_pX$. Since $X^T$ is finite
and each $T_pX$ has finitely many weights, there are finitely many $T$-invariant curves.

(c) By (a), the action of $T$ on $C$ is given by the triple $(p,p',\chi)$ of Proposition~\ref{prop:P1actions}. By Proposition~\ref{prop:P1actions}(a),
\[
\varepsilon_C(p)=d\chi,\qquad \varepsilon_C(p')=-d\chi=-\varepsilon_C(p).
\]
\end{proof}

\begin{lemma}[Edge factors]\label{lem:edgegen}
Let $C$ be a $T$-curve with $C^T=\{p,p'\}$, $\varepsilon:=\varepsilon_C(p)$, and let $f:\PP^1\to C$ be the $k$-fold cover totally
ramified over $p,p'$, equivariant for a finite cover of $T$ acting with tangent weight $\varepsilon/k$ at $f^{-1}(p)$.
\begin{enumerate}
\item[(a)] $TX|_C\cong\bigoplus_{i=1}^{\dim X}L_i$ equivariantly, where $L_i$ has fibre weight $u_i$ at $p$ (the weights of
$T_pX$) and degree $a_i$; one may take $u_1=\varepsilon$, $L_1\cong TC$, $a_1=2$. The weights of $T_{p'}X$ are $u_i-a_i\varepsilon$,
and $c_1(C)=\sum_ia_i=2+\sum_{i\ge2}a_i$.
\item[(b)] In $\chi_T(\PP^1,f^*TX)=H^0-H^1$ the weight $0$ occurs exactly once and the Euler class of the moving part is
\[
E_k(p,C)=\prod_{\substack{m=0\\ m\ne k}}^{2k}\Big(1-\frac mk\Big)\varepsilon\cdot
\prod_{\substack{i\ge2\\ a_i\ge0}}\ \prod_{m=0}^{ka_i}\Big(u_i-\frac mk\varepsilon\Big)\cdot
\prod_{\substack{i\ge2\\ a_i<0}}\ \prod_{m=ka_i+1}^{-1}\Big(u_i-\frac mk\varepsilon\Big)^{-1}.
\]
\item[(c)] $E_k(p,C)$ depends only on the $T$-modules $T_pX$, $T_{p'}X$; in particular $E_k(p,C)=E_k(p',C)$.
\end{enumerate}
\end{lemma}

\begin{proof}
(a) By Proposition~\ref{prop:eqBG}, $TX|_C\cong\bigoplus L(a_i,u_i)$. By (G2) the weight $\varepsilon$ occurs once at $p$; its summand
has weight $(1-a_1)\varepsilon$ at $p'$, which must be the unique weight of $T_{p'}X$ proportional to $\varepsilon$, namely $-\varepsilon$.
So $a_1=2$ and, by uniqueness in Proposition~\ref{prop:eqBG}, $L_1\cong TC$. 

(b) $f^*L_i$ has degree $ka_i$, fiber weight $u_i$ at
$f^{-1}(p)$, and base weight $\varepsilon/k$. Therefore, $f^*L_i=L(ka_i,u_i)$ and the formula for $E_k(p,C)$ follows from Corollary~\ref{cor:ABLP1} with $\epsilon=\varepsilon/k$ and the fact that a weight $u_i-\frac mk\varepsilon$ vanishes, and hence it should be excluded from the moving part, only for $i=1$, $m=k$, by (G2). 

(c) This is a direct consequence of the Lefschetz formula 
\[
\chi_T(\PP^1,f^*TX)=\frac{\mathrm{ch}_T(T_pX)}{1-e^{-\varepsilon/k}}+\frac{\mathrm{ch}_T(T_{p'}X)}{1-e^{\varepsilon/k}}.
\]
\end{proof}

\section{Localization for genus-zero stable maps}\label{sec:loc}

\subsection{Equivariant homology and the virtual class}\label{sec:eqhom}

The moduli space $\overline{\mathcal M}_{0,n+1}(X,\beta)$ is a proper Deligne--Mumford stack with a projective coarse moduli space
$\overline M_{0,n+1}(X,\beta)$, on which $T$ acts. Since we work with rational coefficients, we push all classes forward to the coarse
moduli space. Graber and Pandharipande \cite{GraberPandharipande1999} construct the equivariant virtual fundamental class in the
equivariant Chow group $A^T_*(\overline{\mathcal M}_{0,n+1}(X,\beta))_\Q$ of Edidin--Graham \cite{EdidinGraham1998}, by running the
construction of Behrend--Fantechi \cite{BehrendFantechi1997} for the $T$-equivariant perfect obstruction theory $R\pi_*f^*TX$. We regard
it as an element of the equivariant homology $H^T_*(\overline M_{0,n+1}(X,\beta);\Q)$ by applying the equivariant cycle map. Here, for a
compact $T$-space $M$,
\[
H^T_i(M):=H_{i+2\ell(N-1)}\big(U\times_TM;\Q\big),\qquad U=(\C^N\setminus0)^\ell,\quad N\gg0,
\]
which is independent of $N$ for $N$ large (see \cite{EdidinGraham1998,Brion1997,AndersonFulton2024}). Since $U/T=(\PP^{N-1})^\ell$ and $M$
are compact, $U\times_TM$ is compact and ordinary homology suffices. Note that $H^T_*(M)$ is not the homology of the Borel construction
$ET\times_TM$: it is non-zero in negative degrees (e.g.\ $H^T_{-2k}(\mathrm{pt})\cong H^{2k}_T(\mathrm{pt})$), and it is a module over
$H_T^*(M)$ via cap product.
Equivariant cohomology $H_T^*(M):=H^*(ET\times_TM;\Q)$ is the usual one, defined via the Borel construction; in particular
$H_T^*(\mathrm{pt})=\Sym(\ft^*)$, where $\ft=\mathrm{Lie}(T)$ is the Lie algebra of $T$. It agrees with $H^*(U\times_TM;\Q)$ in degrees $<2N-1$, and cap product on the approximations
$U\times_TM$, compatible with increasing $N$, makes $H^T_*(M)$ a module over $H_T^*(M)$. The classes that we cap with the virtual
class ($\psi$-classes and pull-backs $\ev_i^*\tau$) are Borel classes; with rational coefficients the Borel cohomology of the stack
agrees with that of its coarse moduli space. Finally, for the smooth projective $X$, Poincar\'e duality on the compact complex
manifolds $U\times_TX$ gives the equivariant Poincar\'e duality isomorphism $\mathrm{PD}_X:H_T^{2\dim X-j}(X)\xrightarrow{\sim}H^T_j(X)$,
$a\mapsto a\cap[X]$, which is used in the definition of $\ev^T_{0*}$.

\subsection{Fixed loci and decorated graphs}\label{sec:graphs}

Fix $n\ge-1$ and $\beta\ne0$. A stable map $f:(C;x_0,\dots,x_n)\to X$
is \emph{$T$-fixed} if $t\cdot f\cong f$ for all $t\in T$. In
particular, if $f$ is a fixed point, then there exists a $T$-action on
$C$ such that $f$ is $T$-equivariant. Recalling Propositions
\ref{prop:P1actions} and \ref{prop:gkmcurves}, we get the following
description of the $T$-fixed stable maps.
\begin{lemma}\label{lem:fixedmaps}
If $f$ is $T$-fixed, every irreducible component of $C$ is either contracted to a point of $X^T$ or mapped onto a
$T$-curve by $t\mapsto t^k$, totally ramified over its two fixed points; all nodes and marked points map to $X^T$.
\end{lemma}

To a $T$-fixed $f$ we attach the \emph{decorated graph} $\Gamma$: vertices $v$ are the connected components $Z_v$ of
$f^{-1}(X^T)$, labelled by $p(v):=f(Z_v)\in X^T$; edges $e$ are the
non-contracted components $C_e$, labelled by the pair $(C(e),k_e)$ of
the $T$-curve $C(e):=f(C_e)$ and the degree $k_e$ of the covering
$C_e\to C(e)$  (see Lemma~\ref{lem:fixedmaps}); flags
are pairs $F=(v,e)$ with $Z_v\cap C_e\ne\emptyset$; and $\mu:\{0,\dots,n\}\to V(\Gamma)$ records the marked points. Write $\val(v)$ for the
number of flags at $v$ and $n(v):=\val(v)+|\mu^{-1}(v)|$. For a flag $F=(v,e)$ put
\[
\omega_F:=\frac{\varepsilon_{C(e)}(p(v))}{k_e}\quad(\text{the weight of }T_{y_F}C_e),\qquad E_e:=E_{k_e}(C(e)).
\]
The decorated graph $\Gamma$ has the following properties. 
\begin{lemma}[Properties of $\Gamma$]\label{lem:graphprops}
\begin{enumerate}
\item[(a)] $\Gamma$ is a tree.
\item[(b)] (Compatibility.) If $e$ has flags $F=(v,e)$ and $F'=(v',e)$, then $C(e)^T=\{p(v),p(v')\}$ and $\omega_{F'}=-\omega_F$.
\item[(c)] (Degree.) $\sum_ek_e[C(e)]=\beta$.
\item[(d)] (Vertices.) $Z_v$ is a connected nodal curve of genus $0$, contracted by $f$, iff $n(v)\ge3$; it is then stable
with the $n(v)$ special points $\{y_F\}_{F\ni v}\cup\{x_i\}_{i\in\mu^{-1}(v)}$. Otherwise $Z_v$ is a single point, namely
\begin{itemize}
\item an unmarked end of an edge if $n(v)=1$;
\item a node $C_{e'}\cap C_{e''}$ if $n(v)=2$, $\val(v)=2$;
\item a marked end of an edge if $n(v)=2$, $\val(v)=1$.
\end{itemize}
\end{enumerate}
\end{lemma}

\begin{proof}
(a) $C$ has arithmetic genus $0$, and $\Gamma$ is the dual graph of the decomposition of $C$ into the $Z_v$ and the
$C_e$. (b) is Proposition~\ref{prop:gkmcurves}(c). (c) follows from $(f_e)_*[C_e]=k_e[C(e)]$.
(d) If $Z_v$ contains a contracted component, it is a tree of such components, and stability of $f$ makes it a
stable pointed curve. If not, $Z_v$ is a point. At most two branches pass through a node, marked points are
distinct and smooth, and $\beta\ne0$ forces $\val(v)\ge1$, which gives the three cases.
\end{proof}

An \emph{isomorphism} of decorated graphs is a bijection of vertices, edges and flags preserving incidence, all
labels, and the marking $\mu$. Write $\Aut\Gamma$ for the automorphism group. Conversely, every decorated tree
satisfying (b) and (c) arises from some $T$-fixed map. Let $M_\Gamma\subseteq\Mbar_{0,n+1}(X,\beta)$ be the locus of $T$-fixed
maps of type $\Gamma$. Then $\Mbar_{0,n+1}(X,\beta)^T=\bigsqcup_\Gamma M_\Gamma$, over isomorphism classes, and gluing gives
\[
M_\Gamma\cong\Big[\prod_{v:\,n(v)\ge3}\Mbar_{0,n(v)}\Big/A_\Gamma\Big],\qquad
1\to\prod_e\mathbb Z/k_e\to A_\Gamma\to\Aut\Gamma\to1 .
\]
Here $\mathbb Z/k_e$ are the deck transformations of $f_e$. In particular, integration over $M_\Gamma$ is
\[
\int_{M_\Gamma}=\frac{1}{|\Aut\Gamma|\prod_ek_e}\int_{\prod_v\Mbar_{0,n(v)}} .
\]
For $n(v)\ge3$, the cotangent line of $Z_v$ at $y_F:=Z_v\cap C_e$ restricts to the tautological class $\psi_F$ on $\Mbar_{0,n(v)}$.

\subsection{The virtual normal bundle}

The deformation theory of $f$ is governed by the tangent--obstruction sequence
\[
0\to\Aut(C,x)\to H^0(C,f^*TX)\to T^1_f\to\mathrm{Def}(C,x)\to H^1(C,f^*TX)\to T^2_f\to0,
\]
where $\Aut(C,x)=\operatorname{Ext}^0(\Omega_C(\textstyle\sum x_i),\mathcal O_C)$ and $\mathrm{Def}(C,x)=\operatorname{Ext}^1(\Omega_C(\sum x_i),\mathcal O_C)$
are the infinitesimal automorphisms and deformations of the pointed curve. The sequence is equivariant for the
action of $T$ (through a finite cover). The fixed parts give the obstruction theory of $M_\Gamma$, and the
\emph{virtual normal bundle} $N_\Gamma$ is the moving part of $T^1_f-T^2_f$. Hence, in the representation ring,
\begin{equation}\label{eq:normal}
N_\Gamma|_f=\chi(C,f^*TX)^{\mathrm{mov}}-\Aut(C,x)^{\mathrm{mov}}+\mathrm{Def}(C,x)^{\mathrm{mov}} .
\end{equation}

\begin{lemma}[The three terms]\label{lem:normal}
Let $f$ have type $\Gamma$. Then:
\begin{enumerate}
\item[(a)] $\chi(C,f^*TX)^{\mathrm{mov}}=\sum_e\chi(C_e,f_e^*TX)^{\mathrm{mov}}-\sum_v(\val(v)-1)\,T_{p(v)}X$.
\item[(b)] $\Aut(C,x)^{\mathrm{mov}}=\bigoplus_{F=(v,e):\ n(v)=1}T_{y_F}C_e$ where $y_F\in C_e$ is an unmarked end of an edge $e$, i.e. $Z_v=\{y_F\}$ and $y_F$ is not a marked point. 
The equivariant Euler class is the product of the weights $\omega_F:=e^T(T_{y_F}C_e)$. 
\item[(c)] $\mathrm{Def}(C,x)^{\mathrm{mov}}=\bigoplus_{F=(v,e):\ n(v)\ge3}T_{y_F}Z_v\otimes T_{y_F}C_e\ \oplus\bigoplus_{v:\ n(v)=2,\ \val(v)=2}T_{y_v}C_{e'}\otimes T_{y_v}C_{e''}$,
where in the second sum $Z_v=\{y_v\}$ is the node joining the edges $e',e''$ at $v$. Over $M_\Gamma$ the first summands
have equivariant Euler class $\omega_F-\psi_F$, the second ones $\omega_{F'}+\omega_{F''}$.
\item[(d)] The fixed parts of $T^1_f$ and $T^2_f$ are $T_fM_\Gamma$ and $0$. Hence $M_\Gamma$ carries its fundamental class as
virtual class.
\end{enumerate}
Consequently
\[
\frac{1}{e^T(N_\Gamma)}=\prod_e\frac1{E_e}\;\prod_ve_{p(v)}^{\val(v)-1}\;\prod_{\substack{F=(v,e)\\ n(v)=1}}\omega_F\;
\prod_{\substack{F=(v,e)\\ n(v)\ge3}}\frac{1}{\omega_F-\psi_F}\;\prod_{\substack{v:\ n(v)=2\\ \val(v)=2}}\frac1{\omega_{F'_v}+\omega_{F''_v}} ,
\]
where $F'_v,F''_v$ are the two flags at such a $v$.
\end{lemma}

\begin{proof}
(a) \emph{The normalization (gluing) sequence.} $C$ is obtained from $\bigsqcup_vZ_v\sqcup\bigsqcup_eC_e$ by identifying, for every
flag $F=(v,e)$, the point $y_F\in C_e$ with the point $y_F\in Z_v$. A function on the disjoint union descends to $C$ iff
its two values at each $y_F$ agree, so
\[
0\to\mathcal O_C\to\bigoplus_v\mathcal O_{Z_v}\oplus\bigoplus_e\mathcal O_{C_e}\to\bigoplus_{F}\C_{y_F}\to0
\]
is exact. (When $Z_v$ is a point, $\mathcal O_{Z_v}=\C$ and the sequence glues the edges at $v$ through this point.)
Tensor with the locally free sheaf $f^*TX$. Since $f$ is constant on the connected genus-$0$ curve or point $Z_v$,
$H^0(Z_v,f^*TX)=T_{p(v)}X$ and $H^1(Z_v,f^*TX)=0$. The long exact sequence gives
\begin{align*}
\chi(C,f^*TX)&=\sum_e\chi(C_e,f_e^*TX)+\sum_vT_{p(v)}X-\sum_FT_{p(v_F)}X\\
&=\sum_e\chi(C_e,f_e^*TX)-\sum_v(\val(v)-1)T_{p(v)}X .
\end{align*}
No weight of $T_{p(v)}X$ is zero, by (G1), so these terms are entirely moving. The Euler class of
$\chi(C_e,f_e^*TX)^{\mathrm{mov}}$ is $E_e$ by Lemma~\ref{lem:edgegen}.

(b) Infinitesimal automorphisms of a nodal pointed curve are vector fields on its normalization vanishing at the
preimages of the nodes and of the marked points. On a stable component $Z_v$ there are none. On $C_e\cong\PP^1$,
with coordinate $t$, $t(y_F)=0$, $t(y_{F'})=\infty$, we have $H^0(T\PP^1)=\langle\partial_t,t\partial_t,t^2\partial_t\rangle$ with weights $\omega_F$,
$0$, $-\omega_F=\omega_{F'}$. The point $y_F$ is special (a node or a marked point) iff $n(v)\ge2$. Vanishing at $y_F$ removes
$\partial_t$, and vanishing at $y_{F'}$ removes $t^2\partial_t$. The field $t\partial_t$ is fixed.

(c) There is an exact sequence
\[
0\to\bigoplus_{C'}H^1\big(C',T_{C'}(-\text{special points})\big)\to\mathrm{Def}(C,x)\to\bigoplus_{\text{nodes }y}T_yC'\otimes T_yC''\to0 .
\]
The first sum runs over the components of the normalization, and the second term is the space of first-order
smoothings of the nodes. For an edge, $C_e\cong\PP^1$ carries at most two special points, so the $H^1$ vanishes. The
$H^1$ terms of the components of $Z_v$, together with the smoothings of the nodes inside $Z_v$, form
$T\Mbar_{0,n(v)}$, on which $T$ acts trivially. These are fixed. The
remaining nodes are the $y_F=C_e\cap Z_v$ with $n(v)\ge3$, where
$T_{y_F}C_e\otimes T_{y_F}Z_v$ has weight $\omega_F$, and the nodes $y_v=C_{e'}\cap
C_{e''}$ for which $T_{y_v}C_{e'}\otimes T_{y_v}C_{e''}$ has weight
$\omega_{F'}+\omega_{F''}$. The $T$-lines corresponding to both types
of nodes are moving, that is $\omega_F\ne0$ and
$\omega_{F'}+\omega_{F''}\ne0$. For the former we recall condition
(G1) and for the latter, if the two edges cover different
$T$-curves, their tangent weights at $p(v)$ are non-proportional by (G2); if they cover the same $T$-curve $C$, then
$\omega_{F'}+\omega_{F''}=\varepsilon_C(p(v))(1/k_{e'}+1/k_{e''})\ne0$.

(d) By (a)--(c) and Lemma~\ref{lem:edgegen}(b), the fixed part of $\chi(C,f^*TX)-\Aut(C,x)+\mathrm{Def}(C,x)$ consists of one
weight-$0$ line in $H^0(C_e,f_e^*TX)$ for each edge, cancelled by $t\partial_t\in\Aut(C,x)$, and of $\bigoplus_vT\Mbar_{0,n(v)}=T_fM_\Gamma$.
In particular the fixed part of $H^1(C,f^*TX)$ vanishes, so $(T^2_f)^{\mathrm{fix}}=0$.

The formula for $1/e^T(N_\Gamma)$ follows from \eqref{eq:normal}.
\end{proof}

\subsection{The localization formula}

\begin{proposition}[Localization; Kontsevich,
  Graber--Pandharipande]\label{prop:loc} 
Let $X$ be a GKM variety, $n\ge0$, $\phi_0,\dots,\phi_n\in H^*_T(X)$, and let $g_0,\dots,g_n$ be polynomials. Then
\[
\int_{[\Mbar_{0,n+1}(X,\beta)]^{\mathrm{vir}}}\prod_{i=0}^n\ev_i^*\phi_i\,g_i(\psi_i)
=\sum_\Gamma\frac{1}{|\Aut\Gamma|}\prod_e\frac{1}{k_e\,E_e}\prod_ve_{p(v)}^{\val(v)-1}
\int_{\Mbar_{0,n(v)}}\frac{\prod_{i\in\mu^{-1}(v)}\phi_i|_{p(v)}\,g_i(\psi_i)}{\prod_{F\ni v}(\omega_F-\psi_F)} ,
\]
summed over isomorphism classes of decorated graphs (\S\ref{sec:graphs}). The integrals at vertices with
$n(v)\le2$ are defined by the conventions
\begin{align*}
n(v)=1:&\qquad \int\frac1{\omega_F-\psi}:=\omega_F,\\
n(v)=2,\ \val(v)=2:&\qquad \int\frac{1}{(\omega_{F'}-\psi')(\omega_{F''}-\psi'')}:=\frac1{\omega_{F'}+\omega_{F''}},\\
n(v)=2,\ \val(v)=1,\ \mu^{-1}(v)=\{i\}:&\qquad \int\frac{\phi_i|_{p(v)}\,g_i(\psi_i)}{\omega_F-\psi}:=\phi_i|_{p(v)}\,g_i(-\omega_F).
\end{align*}
The same formula, with no insertions, holds for $\Mbar_{0,0}(X,\beta)$.
\end{proposition}

\begin{proof}
The virtual localization formula of Graber--Pandharipande gives
\[
\int_{[\Mbar_{0,n+1}(X,\beta)]^{\mathrm{vir}}}\alpha=\sum_\Gamma\int_{[M_\Gamma]^{\mathrm{vir}}}\frac{\alpha|_{M_\Gamma}}{e^T(N_\Gamma)} ,
\]
and $[M_\Gamma]^{\mathrm{vir}}=[M_\Gamma]$ by Lemma~\ref{lem:normal}(d). (If $X$ is convex, $\Mbar_{0,n+1}(X,\beta)$ is smooth of the
expected dimension, the virtual class is the fundamental class, and this is the Atiyah--Bott formula.) Write
$\int_{M_\Gamma}=\frac{1}{|\Aut\Gamma|\prod_ek_e}\int_{\prod_{n(v)\ge3}\Mbar_{0,n(v)}}$, and take $1/e^T(N_\Gamma)$ from
Lemma~\ref{lem:normal}. On $M_\Gamma$, $\ev_i^*\phi_i$ restricts to $\phi_i|_{p(\mu(i))}$. If $n(\mu(i))\ge3$, $\psi_i$ restricts to $\psi_i$ on
$\Mbar_{0,n(\mu(i))}$. If $\mu(i)=v$ with $n(v)=2$, $\val(v)=1$, then $x_i=y_F$ lies on $C_e$, and $\psi_i$ restricts to the weight of
$T^*_{y_F}C_e$, namely $-\omega_F$.

Now group the factors by vertices. A stable vertex contributes $e_{p(v)}^{\val(v)-1}$, its insertions, and
$\prod_{F\ni v}(\omega_F-\psi_F)^{-1}$, integrated over $\Mbar_{0,n(v)}$. A vertex with $n(v)=1$ contributes $e^0=1$ and the factor
$\omega_F$ from $\Aut(C,x)^{\mathrm{mov}}$. A node vertex ($n(v)=2$, $\val(v)=2$) contributes $e_{p(v)}$ and $1/(\omega_{F'}+\omega_{F''})$. A
marked end ($n(v)=2$, $\val(v)=1$) contributes $e^0=1$ and $\phi_i|_{p(v)}g_i(-\omega_F)$. These are exactly the three
conventions. The edges contribute $1/(k_eE_e)$.
\end{proof}

\section{Root decomposition and canonical coordinates for GKM varieties}\label{sec:rootgen}

In this section $X$ is a GKM variety as in \S\ref{sec:loc}. Organizing the localization sum around the root vertex and
resumming it with Abel's identity, we obtain an exact formula for the restriction of the big equivariant $J$-function
to a fixed point in terms of its own pole data. It identifies Givental's canonical coordinates, and it is the solved
form of Brown's characterization of the Lagrangian cone.

\subsection{Coefficients and the big \texorpdfstring{$J$}{J}-function}\label{sec:bigJ}

Let $S:=\C[\ft]$ and $K:=\Frac(S)$. Fix an ample class on $X$, so that
for every integer $d$ there are only finitely many effective classes $\beta$
of degree $<d$. Choose a basis $(\phi_a)$ of $H^*_T(X)\otimes_SK$, write $\tau=\sum_at_a\phi_a$ with formal variables $t_a$, and let
\[
\mathcal R:=K(z)[[q,t]]
\]
be the ring of power series in the Novikov variables $q^\beta$ ($\beta$ effective) and the $t_a$, with coefficients in the field
$K(z)$ of rational functions in $z$. Let $I\subset\mathcal R$ be the ideal generated by the $q^\beta$ ($\beta\ne0$) and the $t_a$. Then $\mathcal R$ is
complete and separated in the $I$-adic topology, i.e. every Cauchy
sequence is convergent and it has a unique limit. We regard $K[[q,t]]\subset\mathcal R$ as the elements not depending on $z$.

\paragraph{Big \texorpdfstring{$J$}{J}-function.} Put
\[
J(\tau,z):=-z+\tau+\sum_{\substack{(n,\beta)\ne(0,0),(1,0)}}\frac{q^\beta}{n!}\,
\ev^T_{0*}\Big(\frac{\prod_{i=1}^{n}\ev_i^*\tau}{-z-\psi_0}\cap[\Mbar_{0,n+1}(X,\beta)]^{\mathrm{vir}}\Big),
\]
and let $J_p(\tau,z)$ be its restriction to $p\in X^T$. Let
$\phi_p=i_{p*}(1)$ where $i_p:\{p\}\to X$ is the inclusion map. Since
$i^*_pi_{p*}(\alpha)=e_p\cup \alpha$, we have 
$\phi_p|_{p'}=\delta_{pp'}e_p$. By the projection formula, the coefficient of $q^\beta t^\alpha$ in $J_p$ is an integral over
$[\Mbar_{0,n+1}(X,\beta)]^{\mathrm{vir}}$ with insertion $\ev_0^*\phi_p$ at the point $0$. For $\beta\ne0$, Proposition~\ref{prop:loc} applies with
$\phi_0=\phi_p$, $g_0(\psi)=1/(-z-\psi)$, and $\phi_i=\tau$, $g_i=1$ for $i\ge1$. The marked points $1,\dots,n$ are distributed over the
vertices by $\mu$, and a marked point at $v$ contributes $\tau|_{p(v)}$. For $\beta=0$, where $\Mbar_{0,n+1}(X,0)=X\times\Mbar_{0,n+1}$, the
term is
\begin{equation}\label{eq:vint}
  \frac{(\tau|_p)^n}{n!}\int_{\Mbar_{0,n+1}}\frac{1}{-z-\psi_0}=
  \frac{1}{n!}(\tau|_p)^n\;(-z)^{-n+1},\qquad  n\ge2.
\end{equation}
This is the single-vertex graph, with
the vertex integral of Lemma~\ref{lem:vertex} below.

\begin{lemma}[Poles]\label{lem:poles}
Each coefficient of $J_p(\tau,z)$ is a rational function of $z$, with poles only at $z=0$ and at the points
$\varepsilon_{C'}(p)/j$, for $T$-curves $C'\ni p$ and $j\ge1$. In particular $J_p\in\mathcal R$.
\end{lemma}

\begin{proof}
For fixed $(n,\beta)$ only finitely many graphs occur. In each, $z$ enters only through the vertex $v_0=\mu(0)$. If
$n(v_0)\ge3$ its integral is a polynomial in $z^{-1}$. If $n(v_0)=2$ and $v_0$ is the end of an edge, it is $1/(\omega_F-z)$ with
$\omega_F=\varepsilon_{C'}(p)/j$. If $v_0$ carries no edge, then $\beta=0$, and the contribution is either the term $\tau|_p$
($n=1$) or has the form \eqref{eq:vint} ($n\geq2$); in both cases it is a polynomial in $z$ and $z^{-1}$. This proves that $J_p$ is
a formal power series in $q^\beta$ and $\tau$ whose coefficients are rational functions in $z$ with the claimed pole
structure.
\end{proof}

\paragraph{Branch data.} Recall the set $B$ of all triples $(p,C,k)$ where 
$F=(p,C)$ is a GKM-flag and $k\geq 1$ is an integer. For $b=(p,C,k)\in B$, put 
\[
z_b:=\frac{\varepsilon_C(p)}{k}=-\frac{\varepsilon_C(p')}{k},\qquad
R_{b}(\tau):=q^{k[C]}\,\frac{e_p}{k\,E_k(C)}\;J_{p'}\big(\tau,z_b\big)\in K[[q,t]],
\]
\[
x_b:=R_{b}/z_b,\qquad y_b:=1/z_b,
\]
where $\varepsilon_C(p)$ is the weight of $T_pC$ and
$E_k(C):=E_k(p,C)$ is the edge factor -- see Lemma \ref{lem:edgegen}.
The value $J_{p'}(\tau,z_b)$ makes sense for the following reason. By Lemma~\ref{lem:poles} the poles of $J_{p'}$ are at $0$ and at $\varepsilon_{C'}(p')/j$. Equality
$\varepsilon_{C'}(p')/j=z_b=-\varepsilon_C(p')/k$ would make
$\varepsilon_{C'}(p')$ and $\varepsilon_C(p')$ proportional, hence
$C'=C$ by (G2), and then $-1/j=1/k$, which is impossible. This proves
that $J_{p'}(\tau,z_b)\in K[[q,t]]$. Since $R_{b}\in
q^{k[C]}K[[q,t]]$ and there are finitely many $T$-curves, the family
$\{x_b\}_{b\in B}$ tends to $0$ $I$-adically, i.e. for every $n$ all
$x_b\in I^n$ except for finitely many $b\in B$.

\subsection{The vertex integrals}

\begin{lemma}[Vertex integrals]\label{lem:vertex}
For $n,m\ge0$ and non-zero $\omega_1,\dots,\omega_n$,
\[
\int_{\Mbar_{0,1+n+m}}\frac{1}{z-\psi_0}\prod_{i=1}^n\frac{1}{\omega_i-\psi_i}
=\frac{1}{z\,\prod_i\omega_i}\Big(\frac1z+\sum_{i=1}^n\frac1{\omega_i}\Big)^{n+m-2},
\]
where the $m$ remaining points carry no $\psi$-class. For $1+n+m\le2$ the right side reproduces the unstable
conventions of Proposition~\ref{prop:loc} and the $\beta=0$ terms:
\begin{itemize}
\item $(n,m)=(0,0)$ gives $z$;
\item $(1,0)$ gives $1/(z+\omega_1)$;
\item $(0,1)$ gives $1$, the term $\tau$ in $J$.
\end{itemize}
\end{lemma}

\begin{proof}
For $N\ge3$ and $d_1+\dots+d_N=N-3$ one has $\int_{\Mbar_{0,N}}\prod_j\psi_j^{d_j}=\binom{N-3}{d_1,\dots,d_N}$. This follows by induction on $N$ from
the string equation. Hence $\int_{\Mbar_{0,N}}\prod_{j=1}^N(1-a_j\psi_j)^{-1}=(\sum_ja_j)^{N-3}$. Write $1/(z-\psi_0)=z^{-1}(1-\psi_0/z)^{-1}$
and $1/(\omega-\psi)=\omega^{-1}(1-\psi/\omega)^{-1}$, and put $a_j=0$ at the remaining points. The three unstable cases are checked
directly.
\end{proof}

\subsection{The Lagrange--B\"urmann formula and Abel's identity}
Let us recall the following notation:
for $f(x)\in A[[x]]$, a formal power series with coefficients in some
algebra $A$, we denote by $[x^k]f(x)\in A$ the coefficient of $f(x)$
in front of $x^k$. 
\begin{lemma}[Lagrange--B\"urmann]\label{lem:LB}
Let $A$ be a commutative $\Q$-algebra, $\varphi(s)=\sum_{j\ge0}\varphi_js^j\in A[[s]]$ and $H(s)\in A[[s]]$.
\begin{enumerate}
\item[(a)] There is a unique $u(t)\in tA[[t]]$ with $u(t)=t\,\varphi(u(t))$. Its coefficients are polynomials with integer
coefficients in the $\varphi_j$, so its formation commutes with ring homomorphisms $A\to A'$.
\item[(b)] For every $N\ge1$,
\[
[t^N]\,H(u(t))=\frac1N\,[s^{N-1}]\big(H'(s)\,\varphi(s)^N\big).
\]
\end{enumerate}
\end{lemma}

\begin{proof}
(a) Write $u=\sum_{N\ge1}u_Nt^N$. The coefficient of $t^N$ in $t\varphi(u(t))$ only involves $u_1,\dots,u_{N-1}$. So the equation
determines $u_N=[t^{N-1}]\varphi(u(t))$ recursively and uniquely, as an integer polynomial in $\varphi_0,\dots,\varphi_{N-1}$.

(b) \emph{Formal residues.} For a commutative $\Q$-algebra $A$ let $\operatorname{Res}_s:A((s))\to A$ take the coefficient of $s^{-1}$. We use two
facts.
\begin{itemize}
\item[(R1)] $\operatorname{Res}_sF'(s)=0$ for every $F\in A((s))$.
\item[(R2)] Let $\psi\in sA[[s]]$ with $\psi'(0)\in A^\times$. Then for every $G\in A((t))$ the series $G(\psi(s))$ is defined in $A((s))$, and
$\operatorname{Res}_sG(\psi(s))\psi'(s)=\operatorname{Res}_tG(t)$.
\end{itemize}
(R1) is clear. For (R2) write $\psi=s\,v(s)$ with $v\in A[[s]]^\times$. Then $\psi^m=s^mv^m\in s^mA[[s]]$ for all $m\in\mathbb Z$, so
$G(\psi)=\sum_{m\ge m_0}g_m\psi^m$ converges $s$-adically. Moreover only the finitely many terms with $m_0\le m\le-1$ can
contribute to the coefficient of $s^{-1}$ in $G(\psi)\psi'$. By linearity it suffices to take $G=t^m$. For $m\ne-1$,
$\psi^m\psi'=(\psi^{m+1})'/(m+1)$ has residue $0$ by (R1). For $m=-1$, $\psi'/\psi=s^{-1}+v'/v$ with $v'/v\in A[[s]]$, so the
residue is $1$.

\emph{The case $\varphi_0\in A^\times$.} Put $\psi(s):=s/\varphi(s)\in sA[[s]]$, so $\psi'(0)=\varphi_0^{-1}\in A^\times$. Then $\psi$ has a compositional
inverse. Since $\varphi(u(t))$ has the unit constant term $\varphi_0$, the equation $u=t\varphi(u)$ gives $\psi(u(t))=t$, hence $u=\psi^{-1}$
and $u(\psi(s))=s$. For $N\ge1$, using (R2) with $G(t)=t^{-N}(H\circ u)'(t)$ and the chain rule,
\begin{align*}
N\,[t^N]H(u(t))&=[t^{N-1}](H\circ u)'(t)=\operatorname{Res}_t\,t^{-N}(H\circ u)'(t)
=\operatorname{Res}_s\,\psi(s)^{-N}(H\circ u)'(\psi(s))\,\psi'(s)\\
&=\operatorname{Res}_s\,\psi(s)^{-N}\frac{d}{ds}H\big(u(\psi(s))\big)=\operatorname{Res}_s\,s^{-N}\varphi(s)^NH'(s)
=[s^{N-1}]\big(H'(s)\varphi(s)^N\big).
\end{align*}

\emph{The general case.} Write $H=\sum_jh_js^j$. By (a), both sides of (b) are polynomials with rational coefficients in
$\varphi_0,\dots,\varphi_N,h_1,\dots,h_N$, the same polynomials for every $A$. Let $A_0:=\Q[\varphi_0,\varphi_1,\dots,h_0,h_1,\dots]$ be the
polynomial ring in these symbols, with the universal $\varphi$ and $H$. Since $A_0$ is a domain, $A_0\to A_0[\varphi_0^{-1}]$ is injective,
and (b) holds in $A_0[\varphi_0^{-1}]$ by the previous case. Hence the two polynomials agree in $A_0$. The general case follows
by the ring homomorphism $A_0\to A$ sending the symbols to the given coefficients, which is compatible with the
formation of $u$ by (a).
\end{proof}

\begin{lemma}[Abel's identity]\label{lem:abel}
Let $\mathcal A$ be a commutative $\Q$-algebra, complete and separated for the topology of an ideal $I$. Let $B$ be a
countable set and $(x_b,y_b)_{b\in B}$ a family with $x_b\in I$, $y_b\in\mathcal A$, and $x_b\to0$, i.e.\ for every $n$ we have $x_b\in I^n$
for all but finitely many $b$. For a tuple $(b_1,\dots,b_N)\in B^N$ write $S=S(b_1,\dots,b_N):=y_{b_1}+\dots+y_{b_N}$.
\begin{enumerate}
\item[(a)] There is a unique $u\in I$ with $u=\sum_bx_b\,e^{y_bu}$.
\item[(b)] For every $a\in\mathcal A$,
\[
\sum_{N\ge0}\frac1{N!}\sum_{(b_1,\dots,b_N)\in B^N}\Big(\prod_{i=1}^Nx_{b_i}\Big)\,a\,(a+S)^{N-1}=e^{au},
\]
where the term $N=0$ is $1$, and for $N\ge1$, $a(a+S)^{N-1}$ is a polynomial in $a$.
\end{enumerate}
All series converge $I$-adically.
\end{lemma}

\begin{proof}
\emph{Convergence.} Put $\varphi(s):=\sum_bx_be^{y_bs}\in\mathcal A[[s]]$. Its coefficients $\varphi_j=\sum_bx_by_b^j/j!$ converge and lie in $I$. In the
$N$-th term of (b), each product $\prod_ix_{b_i}$ lies in $I^N$, and for each $n$ only finitely many tuples give a product not in
$I^n$. So the inner sum converges, the $N$-th term lies in $I^N$, and the outer sum converges. Since $u\in I$, the series
$e^{au}$ and $e^{y_bu}$ converge.

\emph{Uniqueness in (a).} If $u,u'\in I$ are solutions, then $u-u'=\sum_j\varphi_j(u^j-u'^j)=(u-u')c$ with $c\in I$, because $\varphi_j\in I$ and
$u^j-u'^j$ is divisible by $u-u'$. Hence $u-u'\in\bigcap_nI^n=0$.

\emph{Evaluation at $t=1$.} Let $\mathcal A\langle t\rangle:=\{c(t)\in
\mathcal A[[t]]:\ [t^N] c(t)\in I^N\}$. This is a subring of $\mathcal A[[t]]$, and
$\mathrm{ev}_1(\sum c_Nt^N):=\sum c_N$ is a ring homomorphism $\mathcal A\langle t\rangle\to\mathcal A$, since the Cauchy product of two $I$-adically convergent
series with these bounds converges to the product. The ring
$\mathcal{A}\langle t\rangle$ has the following decreasing filtration
\[
\mathcal{A}\langle t\rangle=
\mathcal{A}^0\langle t\rangle\supseteq
\mathcal{A}^1\langle t\rangle\supseteq
\mathcal{A}^2\langle t\rangle \supseteq \dots ,\qquad
\mathcal{A}^k\langle t\rangle:=
\Big\{c(t)\in \mathcal{A}\langle t\rangle:\ [t^N]c(t) \in I^k\ \text{for all } N\Big\}
\]
satisfying $\mathcal{A}^i\langle t\rangle \mathcal{A}^j\langle
t\rangle\subseteq \mathcal{A}^{i+j}\langle t\rangle$, that is this is
a ring filtration. 
Then $\mathrm{ev}_1(\mathcal{A}^k\langle t\rangle)\subseteq I^k$, and $\mathrm{ev}_1$ commutes with sums $\sum_nG_n$ with $G_n\in\mathcal{A}^{k_n}\langle t\rangle$, $k_n\to\infty$.

\emph{Existence in (a).} Let $u(t)\in t\mathcal A[[t]]$ be the solution of $u=t\varphi(u)$ from Lemma~\ref{lem:LB}(a). By Lemma~\ref{lem:LB}(b) with
$H(s)=s$, $[t^N]u(t)=\frac1N[s^{N-1}]\varphi(s)^N\in I^N$, so $u(t)\in\mathcal{A}^1\langle t\rangle$ and $u(t)^j\in\mathcal{A}^j\langle t\rangle$. The element $t\varphi_j$ lies in $\mathcal A\langle t\rangle$,
and $t\varphi_j\,u(t)^j\in\mathcal{A}^{j+1}\langle t\rangle$. Applying $\mathrm{ev}_1$ to $u(t)=\sum_j(t\varphi_j)\,u(t)^j$ gives, for $u:=\mathrm{ev}_1(u(t))\in I$, the equation
$u=\sum_j\varphi_ju^j=\varphi(u)=\sum_bx_be^{y_bu}$.

\emph{Proof of (b).} Apply Lemma~\ref{lem:LB}(b) with $H(s)=e^{as}$, so $H'(s)\varphi(s)^N=\sum_{B^N}\prod x_{b_i}\,a\,e^{s(a+S)}$, convergent
coefficientwise in $s$. Taking the coefficient of $s^{N-1}$ gives, for $N\ge1$,
\[
[t^N]\,e^{a\,u(t)}=\frac1{N!}\sum_{(b_1,\dots,b_N)\in B^N}\Big(\prod_ix_{b_i}\Big)\,a\,(a+S)^{N-1}\ \in I^N .
\]
So $e^{au(t)}\in\mathcal A\langle t\rangle$, and $\mathrm{ev}_1(e^{au(t)})$ is the left side of (b). On the other hand
$e^{au(t)}=\sum_na^nu(t)^n/n!$ with $u(t)^n\in\mathcal{A}^n\langle t\rangle$, so $\mathrm{ev}_1(e^{au(t)})=\sum_na^nu^n/n!=e^{au}$.
\end{proof}

\subsection{Proof of Theorem \ref{thm:J}}\label{sec:proof}

Now we are ready to prove our first main result. Let us apply the fixed point localization formula (see Proposition \ref{prop:loc}) to
\[
  J_p(\tau,z)=-z+\tau|_p+\sum_{\substack{(N,\beta)\ne(0,0),(1,0)}}
  \frac{q^\beta}{N!}\,
  \int_{[\Mbar_{0,N+1}(X,\beta)]^{\mathrm{vir}} }
  \frac{\ev_0^*\phi_p}{-z-\psi_0} \cup
  \ev_1^*\tau\cup\dots \cup \ev_N^*\tau
  \ .
\]
We obtain a sum over decorated graphs $\Gamma$ (see Section
\ref{sec:graphs}):
\[
 J_p(\tau,z)=
  \sum_{\Gamma} \frac{1}{|\operatorname{Aut}(\Gamma)| N!} \,
  J_\Gamma(\tau,z)
\]
Each graph $\Gamma$ is a rooted tree. The idea is to decompose $\Gamma$ into subtrees by removing the root and to check that each summand factorizes compatibly. 
\paragraph{The branch sequence.}
The vertex $v_0=\mu(0)$ will be called the {\em root} of $\Gamma$. 
Note that $v_0$ lies over $p$. Let us denote by $e_1,\dots, e_n$ the edges of $\Gamma$ incident with $v_0$. Let $(C_i,k_i)$ be the label of $e_i$ and $v_i$ be the other vertex of $e_i$, i.e. $e_i$ is an edge connecting $v_0$ and $v_i$. 
If we cut the edge $e_i$, then the vertex $v_i$ becomes the root of a subtree $\Gamma_i$ of $\Gamma$. Note that the decoration of $\Gamma$ induces a decoration of $\Gamma_i$. Indeed, this is clear for the decorations of vertices and edges because $\Gamma_i$ is a subgraph of $\Gamma$ while for the marking $\mu_i:\{0,1,\dots,N_i\}\to \{\text{Vertices of } \Gamma_i\}$, we define $\mu_i(0):=v_i$ and we let $\mu_i(1),\dots,\mu_i(N_i)$ to be the subsequence obtained from $\mu(1),\mu(2),\dots,\mu(N)$ by removing the elements that are not vertices of $\Gamma_i$. In particular, $N_i$ is the number of marked points on the vertices of
$\Gamma_i$ and we have $N=m+N_1+\dots+N_n$. Therefore, we have constructed a map 
\begin{equation}\label{eq:branch_map}
\Gamma\to [(C_1,k_1,\Gamma_1),\dots , (C_n,k_n,\Gamma_n) ],
\end{equation}
where each triple $(C_i,k_i,\Gamma_i)$ consists of a $T$-curve $C_i$, an integer $k_i\geq 1$, and a decorated rooted tree $\Gamma_i$ such that $C_i$ connects $p$ and $p_i=p(v_i)$. 
Each such triple (more precisely its isomorphism class) is called a \emph{branch} of $\Gamma$. The square bracket denotes the \emph{permutation} equivalence class of the sequence, i.e. two sequences are said to be equivalent if there is a permutation that transforms one into the other.  Let us point out the properties of the map \eqref{eq:branch_map}. 
\begin{enumerate}
    \item[(i)]
    The map is well defined, i.e. it is independent of the enumeration that we chose for $e_1,\dots,e_n$
    \item[(ii)]
    Two graphs $\Gamma$ and $\Gamma'$ determine the same branch sequence if and only if they differ only by the choice of the marking.
    \item[(iii)]
    Each sequence $[(C_1,k_1,\Gamma_1),\dots , (C_n,k_n,\Gamma_n)]$ of branches is induced by precisely $\tfrac{N!}{m!N_1!\dots N_n!}$ rooted trees $\Gamma$ whose root carries $m+1$ marked points.
    \item[(iv)]
    We have the following short exact sequence of groups:
    \begin{equation}\label{eq:autseq}
1\longrightarrow\prod_{i=1}^n\Aut(\Gamma_i)\longrightarrow\Aut(\Gamma)\longrightarrow
\Aut\big((C_1,k_1,\Gamma_1),\dots,(C_n,k_n,\Gamma_n)\big)\longrightarrow1,
\end{equation}
where the last group is the subgroup of the symmetric group $S_n$ of permutations $\sigma$ with
$(C_{\sigma(i)},k_{\sigma(i)},\Gamma_{\sigma(i)})\cong(C_i,k_i,\Gamma_i)$ for all $i$.
\end{enumerate}

\paragraph{Factorization of $J_\Gamma$.} The term $J_\Gamma(\tau,z)$, in the summand corresponding to $\Gamma$, is a product of numbers associated with the vertices and the edges of $\Gamma$. The factors attached to the root, i.e. the vertex $v_0$ are:
\begin{equation}\label{eq:root_int}
  e_p^n \, (\tau|_p)^m\, \int_{\Mbar_{0,1+n+m}}
  \frac{1}{(-z-\psi_0)(\omega_1-\psi_1)\dots (\omega_n-\psi_n)},
\end{equation}
where
\begin{itemize}
\item the vertex factor $e_p^{\,n-1}$ and the factor $e_p$ from
  $\phi_p|_{p}=e_p$ combine into $e_p^{\,n}$ in total;
\item $(\tau|_p)^m$ comes from the insertions at the marked points of
  $v_0$.  
\item the integral is the vertex integral of $v_0$ with flag weights
  $\omega_i=\varepsilon_{C_i}(p)/k_i=z_{b_i}$ where
  $b_i=(p,C_i,k_i)\in B$.
\end{itemize}
The remaining factors in $J_\Gamma(\tau,z)$ can be identified with the
following product
\[
  \prod_{i=1}^n \frac{q^{k_i [C_i]}}{k_i E_{k_i}(C_i)}
  J_{\Gamma_i}(\tau,\omega_i),
\]
where each edge $e_i$ contributes its edge factor
$1/(k_iE_{k_i}(C_i))$ and the factor $q^{k_i[C_i]}$ coming from the
Novikov variable $q^\beta =\prod_e q^{k_e C(e)}$, and the factors in
$J_\Gamma$ corresponding to vertices and edges in $\Gamma_i$ form a
product that coincides with
$J_{\Gamma_i}(\tau,\omega_i)$. Indeed, recall that the root of
$\Gamma_i$ is the vertex $v_i$ of $\Gamma$ connected to $v_0$ via the
edge $e_i$. Put $p_i=p(v_i)$ and $F'_i=(p_i,C_i)$. Then let us
compare the common vertex and edge factors in $J_\Gamma$ and
$J_{\Gamma_i}$.  
\begin{itemize}
\item At $v_i$, the vertex factor (in $J_\Gamma(\tau,z)$) is
  $e_{p_i}^{\val(v_i)-1}=e_{p_i}^{\val_{\Gamma_i}(v_i)}$, while
  the corresponding vertex factor in $J|_{\Gamma_i}$ is
$e_{p_i}^{\val_{\Gamma_i}(v_i)-1}$ times $e_{p_i}$ from $\phi_{p_i}|_{p_i}$. These agree.
\item The integrand at $v_i$ contains $1/(\omega_{F'_i}-\psi_{F'_i})$ in place of $1/(-z-\psi_0)$; these agree at $z=\omega_i=\omega_{F_i}=-\omega_{F'_i}$ where $F_i:=(p,C_i)$ -- here $F_i$ and $F'_i$ are the two flags of $e_i$.
\item The factors of $J_\Gamma(\tau,z)$ and
  $J_{\Gamma_i}(\tau,\omega_i)$ at the all other vertices and
  all edges of $\Gamma_i$ are identical.
\end{itemize}
We get the following formula 
\begin{equation}\label{eq:J_Gamma}
J_\Gamma(\tau,z)= (\tau|_p)^m\, 
\frac{1}{(-z)\,
\omega_1\cdots\omega_n}\Big(\frac{1}{-z}+\frac1{\omega_1}+\dots+\frac1{\omega_n}\Big)^{n+m-2}
\prod_{i=1}^n \frac{e_p q^{k_i [C_i]}}{k_i E_{k_i}(C_i)}
J_{\Gamma_i}(\tau,\omega_i),
\end{equation}
where the factor $e_p^n$ in the root contribution \eqref{eq:root_int} was distributed into the product of branch contributions and the integral in \eqref{eq:root_int} was evaluated using the formula in Lemma \ref{lem:vertex}.
\paragraph{From graphs to branch sequences.} Since formula \eqref{eq:J_Gamma} expresses $J_\Gamma$ in terms of quantities that depend only on the branches of $\Gamma$, we can transform the sum in the localization formula as follows:
\[
\sum_{\Gamma}\frac{1}{|\operatorname{Aut}(\Gamma)| N!} = 
\sum_{m=0}^\infty 
\frac{1}{m!}
\sum_{[(C_1,k_1,\Gamma_1),\dots,(C_n,k_n,\Gamma_n)]}
\frac{1}{
|\operatorname{Aut}((C_1,k_1,\Gamma_1),\dots,(C_n,k_n,\Gamma_n))|
}
\prod_{i=1}^n\frac{1}{|\operatorname{Aut}(\Gamma_i)| N_i!},
\]
where we used properties (iii) and (iv) of the branch map \eqref{eq:branch_map}.
Since the sum over permutation equivalent classes of sequences weighted by $1/|\operatorname{Aut}((C_1,k_1,\Gamma_1),\dots,(C_n,k_n,\Gamma_n))|$ is the same as the sum over all sequences weighted by $1/n!$, the graph summation takes the following form:
\begin{equation}\label{eq:graph_sum}
\sum_{\Gamma}\frac{1}{|\operatorname{Aut}(\Gamma)| N!} = 
\sum_{n,m=0}^\infty 
\frac{1}{m!n!}
\sum_{((C_1,k_1,\Gamma_1),\dots,(C_n,k_n,\Gamma_n))}
\prod_{i=1}^n\frac{1}{|\operatorname{Aut}(\Gamma_i)| N_i!},
\end{equation}
Let us apply the sum \eqref{eq:graph_sum} to \eqref{eq:J_Gamma}. The sum over all branches splits naturally into sum over all edge markings $(C_1,k_1),\dots,(C_n,k_n)$ and a sum over all decorated rooted trees $\Gamma_1,\dots,\Gamma_n$ such that if $C_i$ is a $T$-curve between $p$ and $p_i$, then the root of $\Gamma_i$ is labeled by $p_i$. Since
\[
\sum_{\Gamma_i} 
\frac{1}{|\operatorname{Aut}(\Gamma_i)| N_i!}
\frac{e_p q^{k_i [C_i]}}{k_i E_{k_i}(C_i)}
J_{\Gamma_i}(\tau,\omega_i) = 
\frac{e_p q^{k_i [C_i]}}{k_i E_{k_i}(C_i)} \, J_{p_i}(\tau,z_{b_i}) =
R_{b_i}(\tau),
\]
where $b_i=(p,C_i,k_i)\in B$ and $z_{b_i} = \varepsilon(C_i)/k_i=\omega_i$, we get 
\begin{equation}\label{eq:rootsum}
J_p(\tau,z)=\sum_{n,m\ge0}\frac{1}{n!}\sum_{(C_1,k_1),\dots,(C_n,k_n)}
\frac{1}{(-z)\,\omega_1\cdots\omega_n}\Big(\frac1{-z}+\frac1{\omega_1}+\dots+\frac1{\omega_n}\Big)^{n+m-2}
\frac{(\tau|_p)^m}{m!}\;\prod_{i=1}^nR_{b_i}(\tau),
\end{equation}
where the inner sum runs over ordered $n$-tuples of pairs ($T$-curve through $p$, multiplicity). The above formula allows us to compute the poles of the $J$-function in the $z$-plane. Poles different from $z=0$ could come only from terms with $n+m\le2$. These are precisely the terms that reproduce the unstable conventions: $(n,m)=(0,0)$ gives $-z$, $(0,1)$ gives $\tau|_p$, and $(1,0)$ gives the term $R_{b}(\tau)/(-z+z_b)$. Therefore, $J_p(\tau,z)$ has a pole at $z=z_b$ with residue $-R_{b}(\tau)$ for all $b\in B(p)=\{b\in B\ |\ p\in b\}$.   

\paragraph{Rewriting \eqref{eq:rootsum}.} 
Put $a=-1/z$, $x_b:=R_{b}(\tau)/z_b$, and $y_b:=1/z_b$ where $b=(p,C,k)\in B(p)$ and $z_b= \varepsilon_C(p)/k$. Since $\omega_i=z_{b_i}$ where $b_i=(p,C_i,k_i)$,  \eqref{eq:rootsum} becomes
\[
J_p(\tau,z)=a\sum_{n,m\ge0}\frac{(\tau|_p)^m}{m!\,n!}
\sum_{(b_1,\dots,b_n)\in B(p)^n}\Big(\prod_ix_{b_i}\Big)
\big(a+S\big)^{n+m-2},
\]
where $S:=y_{b_1}+\dots+y_{b_n}$.
Now regard each insertion at the root as one more branch type $b_0$, with $x_{b_0}=\tau|_p$ and $y_{b_0}=0$. Choosing which $m$ of $N=n+m$ positions carry $b_0$ turns $\frac{1}{m!\,n!}$ into $\frac1{N!}$, and
\[
\frac{J_p(\tau,z)}{-z}=a^2\sum_{N\ge0}\frac1{N!}\sum_{b_1,\dots,b_N}\Big(\prod_ix_{b_i}\Big)(a+S)^{N-2},
\]
with $b_i$ running over $\{b_0\}\cup B(p)$. Here $x_b\to0$ $I$-adically, so Lemma~\ref{lem:abel} applies in $\mathcal A=\mathcal R$. The terms with
$N\le1$ involve $(a+S)^{-1}$ and $(a+S)^{-2}$, which are invertible in $\mathcal R$ since $a+y_b=-(1-zy_b)/z$.

\paragraph{Resummation via the Abel's identity.} Write $a^2(a+S)^{N-2}=a(a+S)^{N-1}-aS(a+S)^{N-2}$. By Lemma~\ref{lem:abel}(b), the first part sums to 
$e^{au}$, where $u\in I$ solves $u=\sum_bx_be^{y_bu}$. Since $x_{b_0}=\tau|_p$ and 
$x_b=R_{b}(\tau)/z_b$ for $b=(p,C,k)\in B(p)$, this is the equation defining $u_p$. In the second
part, $S=\sum_iy_{b_i}$, and the summand is symmetric in the positions. Singling out the last branch $b_N=b$ and writing $S'$ for
the sum over the others, the $N$ choices of the singled-out position cancel the factor $N$ in $N!$, and we get
\[
\sum_bx_by_b\sum_{M\ge0}\frac1{M!}\sum_{b_1,\dots,b_M}\Big(\prod x_{b_i}\Big)\,a\,(a+y_b+S')^{M-1}
=\sum_bx_by_b\,\frac{a}{a+y_b}\,e^{(a+y_b)u}.
\]
Here we used Lemma~\ref{lem:abel}(b) with $a$ replaced by $a+y_b$, multiplied by $a/(a+y_b)$. The branch $b_0$ has $y_{b_0}=0$ and drops out. Hence
\[
\frac{J_p(\tau,z)}{-z}=e^{-u/z}\Big(1-\sum_{b=(F,k)}x_{b}y_{b}\,\frac{e^{y_{b}u}}{1-zy_{b}}\Big)=
e^{-u/z}\Big(1+\sum_{b=(F,k)}\frac{x_b\, e^{y_{b}u}}{z-1/y_{b}}\Big),
\]
which is the formula stated in Theorem \ref{thm:J}, since $x_{b}=R_{b}(\tau)/z_b$ and 
$y_{b}=1/z_b$.

\paragraph{Canonical coordinates.} It remains only to prove that $u_p$ ($p\in X^T$) are canonical coordinates. 
\begin{proof}
Let us denote by $QH^*_T(X):=H^*_T(X)\otimes_S K[[q,t]]$ the big quantum cohomology of $X$. 
In this proof we use the fixed-point basis $\phi_a=i_{a*}(1)$ ($a\in X^T$) of $H^*_T(X)\otimes_SK$, and we denote by
$\partial_a:=\partial/\partial t_a$, where $t_a$ are the coordinates with respect to this basis, so that $\tau|_p=e_p\,t_p$.
Note that the Jacobi matrix $(\partial_a(u_p))$ ($p,a\in X^T$) is non-degenerate because its restriction to $q=0$ is an invertible diagonal matrix: the diagonal entries are $e_p=e^T(T_pX)$ which are non-zero thanks to the GKM condition (G1). Therefore, $u_p$ ($p\in X^T$) are formal coordinates. We have to prove that 
$\partial/\partial u_p := \sum_a \partial_{u_p} (t_a) \phi_a$ are idempotents for the quantum product:
\[
\frac{\partial}{\partial u_p}\bullet 
\frac{\partial}{\partial u_{p'}}=\delta_{p,p'} 
\frac{\partial}{\partial u_p}\qquad
\Leftrightarrow \qquad
\sum_{a,b}
\frac{\partial t_a}{\partial u_p}\,
\frac{\partial t_b}{\partial u_{p'}} \, \phi_a\bullet \phi_b =
\sum_c 
\frac{\partial t_c}{\partial u_p} \delta_{p,p'} \phi_c\, .
\]
Since the Jacobi matrix is invertible, this is equivalent to proving that the $K[[q,t]]$-linear maps 
\[
\chi_p: QH_T^*(X)\to K[[q,t]],\qquad 
\chi_p(\phi_a):= \partial_a(u_p)
\]
are ring homomorphisms for the quantum product, i.e. $\chi_p(\phi_a\bullet \phi_b)=\chi_p(\phi_a)\chi_p(\phi_b).$
This however follows immediately from the fact that the $J$-function satisfies the following differential equations:
\[
\Big(z\frac{\partial^2}{\partial t_a \partial t_b} + 
\sum_{c\in X^T} A_{ab}^c(q,t)\, 
\frac{\partial}{\partial t_c}\Big)\, J_p(\tau,z) = 0,
\]
where $A_{ab}^c(q,t)$ are the structure constants of the quantum product:
$\phi_a\bullet\phi_b=\sum_c A_{ab}^c(q,t)\phi_c$.
Indeed, we already proved that $J_p(\tau,z)=-z e^{-u_p/z} g_p(q,t,z)$ where 
$g_p(q,t,z)\in K[[z]][[q,t]]$ and $g_p\equiv1$ modulo $I$. The above differential equation is an identity in
$K(z)[[q,t]]$, because the coefficients of $J_p$ are rational functions of $z$ (Lemma~\ref{lem:poles}). Hence it remains
valid after expanding these rational functions in Laurent series at $z=0$; in this expansion $e^{-u_p/z}\in K[z^{-1}][[q,t]]$
and the factorization $J_p=-z e^{-u_p/z}g_p$ holds. Substituting it into the differential equation and canceling
$e^{-u_p/z}$ we get 
\[
\Big(
\sum_c A_{ab}^c(q,t) \partial_c(u_p)-
\partial_a(u_p)\partial_b(u_p)\Big)g_p(q,t,0) + z(\cdots) =0,
\]
where the dots involve terms with only non-negative powers of $z$. 
The vanishing of the coefficient in front of $z^{0}$, since $g_p(q,t,0)\equiv1$ modulo $I$ is invertible, coincides with the condition that $\chi_p$ is a ring homomorphism. 
This completes the proof of Theorem \ref{thm:J}.
\end{proof}

\subsection{Comparison with Brown's characterization}\label{sec:brown}
In this section we will give an alternative proof of Theorem \ref{thm:J} by using Brown's characterization of the Givental's cone.
\paragraph{The cone and Brown's conditions.}
Let $\mathcal H:=H^*_T(X)\otimes K((z^{-1}))[[q]]$ with the symplectic form $\Omega(f,g)=\Res_{z=0}\big(f(-z),g(z)\big)dz$, and let
$\mathcal L\subset\mathcal H$ be Givental's Lagrangian cone of the equivariant genus-$0$ theory \cite{Givental2004,CoatesGivental2007}: its
points are
\[
f=-z+t(z)+\sum_{n,\beta}\frac{q^\beta}{n!}\sum_a\Big\langle t(\psi),\dots,t(\psi),\frac{\phi_a}{-z-\psi}\Big\rangle_{0,n+1,\beta}\phi^a,
\qquad t(z)\in H^*_T(X)\otimes K[z][[q,t]]
\]
(in a formal neighbourhood of $-z$). With the normalization \eqref{eq:Jdef}, $J(\tau,\cdot)\in\mathcal L$ is the point with
$t(z)=\tau$. In the toric case Brown \cite[Thm.~2]{Brown2014} and Coates--Corti--Iritani--Tseng \cite[\S\S3--4]{CCIT2015}
characterised the points of $\mathcal L$ by their restrictions to fixed points; their proof uses only the localization of
\S\ref{sec:loc} and the structure of $T$-curves, and applies verbatim to GKM manifolds (cf.\ also \cite{Givental1996,Liu2013}).
In our normalization the characterization reads as follows. Let $f\in\mathcal H$ with $f\equiv-z$ modulo $I$.

\begin{itemize}
\item[(C1)] For every $p$, coefficientwise in $(q,t)$, $f|_p$ is a rational function of $z$ whose poles away from $z=0,\infty$
are simple, lie at the points $z_b=e^T(T_pC)/k$ where $b=(p,C,k)\in B$, and satisfy
\[
\Res_{z=z_b}f|_p(z)=-q^{k[C]}\,\frac{e_p}{k\,E_{k}(C)}\;f|_{p'}(z_b),
\]
where $p'\in C^T\setminus\{p\}$ and $E_k(C)$ is the edge factor \eqref{eq:edge}. 
\item[(C2)] For every $p$, the Laurent expansion of $f|_p$ at $z=0$ is a point of the Lagrangian cone $\mathcal L_{\mathrm{pt}}$ of a point,
i.e.\ $f|_p=-z\,e^{-c_p/z}g(z)$ with $c_p\in I$ and $g\in K[[z]][[q,t]]$.
\end{itemize}
The Brown's characterization of $\mathcal L$ can be stated as follows (see \cite{Brown2014,CCIT2015}).
\begin{proposition}[Brown's Theorem]\label{thm:brown}
$f\in\mathcal L$ if and only if $f$ satisfies \textup{(C1)} and \textup{(C2)}.
\end{proposition}
\begin{remark}
In \cite{Brown2014,CCIT2015} condition (C2) is formulated with the cone of a point twisted by $e^{-1}_T(T_pX)$. In genus $0$ this
twist is harmless: the twisting class on $\Mbar_{0,n}(\mathrm{pt})$ is $e^{-1}_T(R\pi_*\ev^*T_pX)=1/e_p$, a constant, which rescales the
potential and the pairing by the same factor; hence the twisted cone equals $\mathcal L_{\mathrm{pt}}=\bigcup_cz\,e^{-c/z}K[[z]]$.
\end{remark}

We will use Proposition \ref{thm:brown} to give an alternative proof of Theorem \ref{thm:J}.  
Put $f(q,t,z)=J(\tau,z)$. By definition $f\in \mathcal L$ and therefore $f$ satisfies conditions \textup{(C1)} and 
\textup{(C2)}. Put
\[
R_{b}:=q^{k[C]}\frac{e_p}{kE_{k}(C)}f|_{p'}(z_b),\qquad
b:=(p,C,k)\in B,
\]
then by condition \textup{(C1)} we have $\Res_{z=z_b}f|_p(z)=-R_{b}$. We first prove that 
\[
f|_p(z)=-z\,e^{-c_p/z}\Big(
1+\sum_{b\in B(p)}\frac{\rho_b}{z-z_b}\Big),\qquad 
c_p=\tau|_p+\sum_{b\in B(p)} \rho_{b},
\]
where $\rho_{b}:=R_{b}e^{c_p/z_b}/z_b$ and $B(p)=\{b\in B\ |\ p\in b\}$.  
Coefficientwise, $h(z):=f|_p(z)\,e^{c_p/z}/(-z)$ is rational ($c_p\in I$, so 
$e^{c_p/z}$ is a polynomial in $z^{-1}$ in each degree). 
By (C2) $h=g$ near $0$, so $h$ is regular at $0$. 
By (C1) its other finite poles are simple, they are at $z=z_b$, and the corresponding residues
$-R_{b}e^{c_p/z_b}/(-z_b)=\rho_{b}$. At infinity, since 
$f|_p=-z+\tau|_p+O(z^{-1})$, we have
$h=(1-\tau|_p/z+O(z^{-2}))(1+c_p/z+O(z^{-2}))\to1$. By partial fractions
$h=1+\sum\rho_{b}/(z-z_b)$, and the coefficient of $z^{-1}$ gives 
$c_p-\tau|_p=\sum_{b\in B(p)} \rho_{b}$. Finally, it remains only to prove that $c_p$ ($p\in X^T$) 
are canonical coordinates but this is done in the same way as in the last step of the proof
in Section \ref{sec:proof}.


\section{Analyticity of the big J-function}\label{sec:analytic}

In this section we show that Theorem~\ref{thm:J} does more than determine the big $J$-function as a formal power series:
read as a fixed-point problem for the residues, it is a contraction in a suitable Banach space of sequences. This gives the convergence of the big equivariant $J$-function near the large radius limit and the meromorphy of $J_p(\tau,z)$ in $z\in\C^\times$.

\subsection{The fixed-point problem}
We will view equations \eqref{eq:up} and \eqref{eq:Jp} as equations for
\begin{equation}\label{eq:x-var}
x=(u,R)\ \in\ \C^{X^T}\times \C^B,\qquad 
u=(u_p)_{p\in X^T},\qquad R=(R_b)_{b\in B}
\end{equation}
depending on the complex parameters
\[
\tau|_p\ (p\in X^T),\qquad 
q^{[C]}\ (C \text{ a $T$-curve}),\qquad
z\in \C_B^*:=\C^*\setminus \{ z_b\ |\ b\in B\}.
\]
The remaining ingredients of the formulas are either explicitly known: 
$e_p$, $z_b$, and $E_k(C)$, or we can express them in terms of the variables 
\eqref{eq:x-var}, that is $J_p(\tau,z)= \mathcal J_p(x;z)$ where
\begin{equation}\label{eq:ker-J}
\mathcal J_p(x;z):=-z\,e^{-u_p/z}\Big(1-\sum_{b\in B(p)}\frac{R_{b}}{z_b^{2}}\,\frac{e^{u_p/z_b}}{1-z/z_b}\Big).
\end{equation}
Put $\Phi(x):=\Phi(x;q,\tau)$ for the following vector of formal power series:
\begin{equation}\label{eq:Phi}
\Phi_p(x):=\tau|_p+\sum_{b\in B(p)}\frac{R_{b}}{z_b}\,e^{u_p/z_b},\qquad
\Phi_b(x):=\big(q^{[C]}\big)^k\,\frac{e_p}{k\,E_{k}(C)}\;
\mathcal J_{p'}\big(x;z_b\big),
\end{equation}
where $p\in X^T$ and $b\in B.$ We would like to prove that the fixed point problem $\Phi(x)=x$ has an analytic solution $x=(u(q,t), R(q,t))$ such that the series \eqref{eq:ker-J} is uniformly convergent on all compact subsets of 
\[
|q^{[C]}|<r\ (\forall \text{$T$-curves $C$}),\qquad 
|\tau|_p|<r\ (p\in X^T),\qquad z\in \C^*_B
\]
where $r$ is a sufficiently small positive number to be specified later on. 

To begin with, let us point out that according to Theorem \ref{thm:J}, the fixed-point problem has a unique formal solution such that $u_p, R_b\in I$ where $I$ is the ideal $(q,t)$ in $K[[q,t]]$ where $K=\operatorname{Frac} \C[\ft]$. Note however that the coefficients of the formal power series $u_p$ and $R_b$ are rational functions on $\ft$ with possible poles only at $\lambda\in \ft$ such that $u(\lambda)=0$ for some weight $u$ of $T_pX$, or $u(\lambda)/u'(\lambda)\in \mathbb Q$ for some pair $u\neq u'$ of weights of $T_pX $, for some fixed point $p\in X^T$ -- this 
is because the coefficients produced by the recursion
involve only the denominators $e_p$, $E_{k}(C)$, $z_b$ and $1-z_b/z_{b'}$.

\begin{definition}\label{def:admissible}
A point $\lambda\in\ft$ is called \emph{admissible} if for every $p\in X^T$ we have
\begin{enumerate}
    \item[(i)] $u(\lambda)\neq 0$ for all weights $u$ of $T_pX$.
    \item[(ii)] $u(\lambda)/u'(\lambda)\notin\mathbb R$ for all 
    distinct weights $u\ne u'$ of $T_pX$.
\end{enumerate}
\end{definition}

Since holomorphic functions $f:U\to \C$ defined on some open subset $U\subset \ft$ can not take only real values unless they are constant, by using the GKM condition \textup{(G2)}, we get that admissible points form a non-empty dense open subset of $\ft$. Moreover, the set of admissible points is stable under $\lambda\mapsto\zeta\lambda$, $\zeta\in\C^\times$. 

From now on we fix an admissible $\lambda\in \ft$. All equivariant parameters, viewed as elements in $H^*_T({\rm pt})=\C[\ft]$, will be evaluated at $\lambda$. Finally, we also specify the Novikov variables to be complex numbers, so that $q^{[C]}\in \C$. Since there are finitely many $T$-curves, by choosing the Novikov variables small enough, we can arrange that $|q^{[C]}|<r$ for any fixed $r>0$. Therefore, in the fixed-point problem, all known parameters and the unknown variables will be complex numbers while the formal solution has components in $\C[[q,t]]$.

For admissible $\lambda$ all weights are non-zero, and we define the following constants that depend only on the GKM manifold $X$ and $\lambda$:
\[
\varepsilon_{\min}:=\min_{p,\,u}|u(\lambda)|,\qquad \varepsilon_{\max}:=\max_{p,\,u}|u(\lambda)|,\qquad
s_0:=\min_{p}\ \min_{u\ne u'}\ \min_{t\in\mathbb R}\frac{|u(\lambda)-t\,u'(\lambda)|}{|u(\lambda)|}\in(0,1],
\]
the minima over weights $u,u'$ of $T_pX$. Thus $s_0$ is the sine of the smallest angle between two tangent weights at a fixed point, viewed as vectors in $\C\cong\mathbb R^2$. Furthermore, we introduce the notation
\[
n:=\dim X,\qquad \gamma:=1/\varepsilon_{\min},\qquad
|b|=k\qquad \text{ for } b=(p,C,k)\in B
\]
and we write $p'=p'(b)$ for the other fixed point of $C$ where $b=(p,C,k)\in B$.
Finally, the following estimates for the poles of the $J$-function are direct consequence of the above definitions:
\[
\varepsilon_{\min}/k\leq |z_b| \le\varepsilon_{\max},
\qquad
|1/z_b|\le\gamma|b|.
\]

\subsection{Analyticity of the fixed-point problem}
Our fixed-point problem was defined only formally. On the other hand, the vector space $\mathfrak X:=\C^{X^T}\times \C^B$ has a structure of Banach vector space given by the following norm: 
\[
\|x\|:=\max\Big(
\frac{\max_p|u_p|}{\eta},\ 
\frac{\sup_{b\in B}|R_{b}|\,\rho^{-|b|}}{M}\Big),
\]
where the constants $\eta,\rho,$ and $M$ will be specified later on. 
Let $\Omega:=\{\|x\|<1\}$ be the open unit ball and 
$|\tau|_p|<r, |q^{[C]}|<r$. We would like to prove that the constants $\eta,\rho,M$, and $r$ can be chosen in such a way that $\Phi:\Omega\to \mathfrak X$ is a holomorphic map depending holomorphically on the parameters $\tau$ and $q$. 
The proof is based on the Weierstrass $M$-test which is known to hold for Banach spaces -- see Proposition \ref{prop:M-test} in the Appendix for a precise statement and proof. 
The strategy of the proof is the following. First, note that the components of the map $\Phi$ ( 
see \eqref{eq:ker-J} and \eqref{eq:Phi}) are infinite series whose summand terms are holomorphic functions on $\mathfrak X$: each summand depends polynomially on finitely many of the variables $R_b$ and $e^{u_p/z_b}$ ($p\in X^T$, $b\in B$). Therefore, our main task is to 
find an upper bound for each term in the summand which is uniform in $x\in \Omega$, $|\tau|_p|<r, |q^{[C]}|<r$ and such that we can apply Proposition \ref{prop:M-test}.

\begin{lemma}[Edge factors grow at most exponentially]\label{lem:edgebound}
There are constants $A,\Theta >0$, depending only on the GKM data at $\lambda$, such that
\[
\Big|\frac{e_p}{k\,E_{k}(C)}\Big|\le A\,\Theta^k\qquad\text{for all }
b=(p,C,k)\in B.
\]
\end{lemma}

\begin{proof}
Recall the product formula of Lemma~\ref{lem:edgegen}(b), with $u_1=\varepsilon=\varepsilon_C(p)$. The tangent part has absolute value
\[
\prod_{j=-k,\,j\ne0}^{k}|j\varepsilon/k|=(k!)^2(|\varepsilon|/k)^{2k}\ge(|\varepsilon|/e)^{2k},
\]
since $k!\ge(k/e)^k$. For $i\ge2$ the weight $u_i$ is distinct from
$\varepsilon$, so every factor $u_i-\frac mk\varepsilon$ with $a_i\ge0$ has absolute value $\ge s_0|u_i|\ge s_0\varepsilon_{\min}$; there are $ka_i+1$ of them.
Every $H^1$-factor ($a_i<0$) has absolute value $\le|u_i|+|a_i||\varepsilon|\le(1+|a_i|)\varepsilon_{\max}$; there are $k|a_i|-1$ of them. Finally
$|e_p|\le\varepsilon_{\max}^n$. Since there are finitely many $T$-curves, the claim follows, with
\[
\Theta=\max_{p,C}\big(e/\varepsilon_{\min}\big)^2\,\max(1,(s_0\varepsilon_{\min})^{-1})^{\sum_{a_i>0}a_i}\,\max(1,\varepsilon_{\max})^{\sum_{a_i<0}|a_i|}\prod_{a_i<0}(1+|a_i|)^{|a_i|}.
\]
\end{proof}

\begin{lemma}[No small divisors]\label{lem:nosmall}
Let $b\in B(p)$ with other end $p'$, and let $b'\in B(p')$. Then $|1-z_b/z_{b'}|\ge s_0$.
\end{lemma}

\begin{proof}
Write $b=(C,k)$, $b'=(C'',j)$. Then $z_b=\varepsilon_C(p)/k=-\varepsilon_C(p')/k$ and $z_{b'}=\varepsilon_{C''}(p')/j$, so
$1-z_b/z_{b'}=\big(\varepsilon_{C''}(p')+\tfrac jk\varepsilon_C(p')\big)/\varepsilon_{C''}(p')$. If $C''=C$ this equals $1+j/k\ge1$. Otherwise $\varepsilon_C(p')$ and
$\varepsilon_{C''}(p')$ are distinct weights of $T_{p'}X$, and the numerator has absolute value $\ge s_0|\varepsilon_{C''}(p')|$.
\end{proof}

\begin{remark}
Lemmas~\ref{lem:edgebound} and~\ref{lem:nosmall} are the places where admissibility is essential: both use the bound
$|u-t\,u'|\geq s_0|u|$ for $t\in\mathbb R$. For real generic $\lambda$ the ratios of distinct weights at a fixed point are real
and irrational, and both the factors $u_i-\frac mk\varepsilon$ of the edge factors and the numbers $1-z_b/z_{b'}$ can be
arbitrarily small: one faces a small-divisor problem.
\end{remark}

Put 
$\eta:=\varepsilon_{\min}=1/\gamma$, 
$M:=16A\varepsilon_{\max}$, and write 
$\Phi(x)=(\tau,0)+\Psi(x)$.

\begin{lemma}\label{lem:psibound}
 There are $\rho\in(0,\frac12]$ and $r>0$ such that, for $|q^{[C]}|\le r$ for all $T$-curves
$C$ and all $\tau$, the map $\Psi$ is defined on $\Omega$, each of its coordinates is a holomorphic function of $(x,q,\tau)$, and
$\|\Psi(x)\|\le\frac18$ for $x\in\Omega$.
\end{lemma}

\begin{proof}
\emph{The $u$-component.} Since 
$x\in\Omega$, we have 
\[
|u_p|<\eta,\qquad 
|R_{b}|<M\rho^{|b|},\qquad
|e^{\pm u_p/z_b}|\le e^{\gamma\eta|b|}=e^{|b|},
\]
where we used that $|1/z_b|\le\gamma|b|$. Suppose that $\kappa:=e\rho\leq \tfrac12$, then 
\[
|\Psi_p(x)|\leq \sum_{b\in B(p)} |R_b e^{u_p/z_b}/z_b|\leq 
\sum_{b\in B(p)} M\gamma |b| (e\rho)^{|b|} \leq 
n M\gamma \sum_{k=1}^\infty k (e\rho)^k \leq 
4nM\gamma e\rho,
\]
where we used that for $b=(p,C,k)\in B$ with $p$ fixed there are at most $n$ choices for the $T$-curve $C$ passing through $p$  and the inequality $\sum_{k=1}^\infty k \kappa^k\leq 4\kappa$ for $\kappa=e\rho\leq \tfrac12$. In order to achieve $|\Psi_p(x)|\leq \eta/8$, it is sufficient to require that 
\[
\rho\leq \operatorname{min}\Big(\frac1{2e}, \frac\eta{32nM\gamma e}\Big).
\]
\emph{The values $\mathcal J_{p'}(x;z_b)$.} Let $b=(p,C,k)\in B$ and $z=z_b$. Then $|z|\le\varepsilon_{\max}$, $|e^{-u_{p'}/z}|\le e^{\gamma\eta k}=e^k$, and by
Lemma~\ref{lem:nosmall} the bracket is bounded by 
\[
1+\sum_{b'\in B(p')}M\rho^{|b'|}\gamma^2|b'|^2e^{\gamma\eta|b'|}/s_0\le
1+\frac{nM\gamma^2}{s_0} \sum_{k=1}^\infty k^2 (e\rho)^k\le
1+12nM\gamma^2\kappa/s_0 ,
\]
where we used $\sum_kk^2\kappa^k\le12\kappa$ for $\kappa=e\rho\le\frac12$.
Therefore, the term in the bracket is $\le2$ if $\kappa\le s_0/(12nM\gamma^2)$. 
Hence $|\mathcal J_{p'}(x;z_b)|\le2\varepsilon_{\max}e^{k}$ for 
\[
\rho\leq \operatorname{min}\Big(\frac1{2e}, \frac{s_0}{12nM\gamma^2 e}\Big).
\]
\emph{The $R$-component.} By Lemma~\ref{lem:edgebound},
\[
|\Psi_b(x)|\le r^kA\Theta^k\cdot2\varepsilon_{\max}e^{ k}=\frac M8\,(r\Theta e)^k\le\frac M8\rho^k
\]
if $r\le\rho/(\Theta e)$. 

So choose $\rho$ with $e\rho\le\min\big(\tfrac12,\eta/(32nM\gamma),s_0/(12nM\gamma^2)\big)$ and then $r:=\rho/(e\Theta)$. All series involved converge
uniformly on $\Omega\times\{|q^{[C]}|\le r\}$, so every coordinate of $\Psi$ is
independent of $\tau$ and holomorphic in $(x,q)$.
\end{proof}

\begin{proposition}\label{prop:contraction}
With $\rho,r$ as in Lemma~\ref{lem:psibound}, let $|q^{[C]}|\le r$ for all $C$ and $|\tau|_p|\le\eta/8$ for all $p$. Then $\Phi$ maps the closed
ball $\overline\Omega_{1/4}:=\{\|x\|\le\frac14\}$ into itself and is a contraction there with constant $\frac13$. Its unique fixed point
$x^*(q,\tau)$ depends holomorphically on $(q,\tau)$, coordinate by coordinate, and its Taylor expansion at $(q,\tau)=0$ is the formal
solution given by the big $J$-function.
\end{proposition}

\begin{proof}
For $\|x\|\le\frac14$, $\|\Phi(x)\|\le\frac18+\frac18=\frac14$. Let $x,x'\in\overline\Omega_{1/4}$ with $x\ne x'$, put $h:=x'-x$ (so $\|h\|\le\frac12$) and
$e:=h/\|h\|$, and fix one coordinate $\Psi_i$ together with its weight $w_i$ in the norm ($w_i=1/\eta$ or $\rho^{-|b|}/M$). The function
$g(\zeta):=w_i\Psi_i(x+\zeta e)$ is holomorphic for $|\zeta|<\frac34$, because then $x+\zeta e\in\Omega$, and $|g|\le\frac18$ by Lemma~\ref{lem:psibound}.
By the Schwarz lemma, $|g(\zeta)-g(0)|\le2\cdot\frac18\cdot|\zeta|/\frac34=|\zeta|/3$. Taking $\zeta=\|h\|$ gives $w_i|\Psi_i(x')-\Psi_i(x)|\le\frac13\|x'-x\|$, and
taking the supremum over $i$, $\|\Phi(x')-\Phi(x)\|=\|\Psi(x')-\Psi(x)\|\le\frac13\|x'-x\|$.

The iterates $x^{(0)}:=(\tau,0)$, $x^{(m+1)}:=\Phi(x^{(m)})$ have coordinates holomorphic in $(q,\tau)$ and converge to $x^*$ uniformly in
$(q,\tau)$, so the coordinates of $x^*$ are holomorphic. Expanding the identity $x^*=\Phi(x^*)$ in Taylor series in $(q,\tau)$ shows that
the Taylor coefficients of $x^*$ form a formal solution in the ideal $(q,\tau)$; by uniqueness it is the one given by the big $J$-function.
\end{proof}

\paragraph{Proof of Theorem~\ref{thm:J-convergence}.}
Let us choose $0<r<\eta/8$ that is smaller than the radius of Proposition~\ref{prop:contraction}. 
Let $x^*(q,\tau)$ be the solution of the fixed-point problem constructed in Proposition~\ref{prop:contraction}. By substituting the Taylor series expansion at $q=0,\tau=0$ of $x^*(q,\tau)$ in  $\mathcal J_p(x^*;z)$  we obtain a formal series which must coincide with 
big J-function $J_p(\tau,z)$ because they both are solutions to the equations \eqref{eq:up}--\eqref{eq:J-res}. Therefore, we need only to prove that $\mathcal J_p(x^*,z)$ is convergent uniformly in $(q,\tau)$ on the compact subsets $K\subset \C^\times$ that are disjoint from the poles $z_b$ ($b\in B$).
Since $z_b\to0$, there is $\delta_K>0$ with $|1-z/z_b|\ge\delta_K$ for all $z\in K$ and all $b\in B$. 
The estimates that we need are essentially the same as in the proof of Lemma \ref{lem:psibound}. Indeed, let us recall the definition \eqref{eq:ker-J} and focus on the estimate of the summand on the RHS. Using that $x^*\in \Omega$, we have $|u_p|<\eta$ and $|R_b|<M \rho^{|b|}$. We also have $|1/z_b|\leq \gamma |b|$
and therefore
\[
\left|
\frac{R_b}{z_b^2}\, \frac{e^{u_p/z_b}}{1-z/z_b}
\right| \leq 
\frac{|R_b| e^{|u_p|/|z_b|}}{|z_b|}\, \frac{1}{|z_b|\,|1-z/z_b|} \leq M \gamma^2 \delta_K^{-1} |b|^2 (e\rho)^{|b|}.
\]
Since we fixed $\rho$ such that $e\rho\leq 1/2$ and at most $n$ $T$-curves pass through $p$, we have
\[
\sum_{b\in B(p)} |b|^2 (e\rho)^{|b|} \leq 12 n e\rho,
\]
and therefore, by recalling the Weierstrass $M$-test, we get that the series $\mathcal J_p(x^*,z)$ converges uniformly on $K$. Near a given pole $z_{b_0}$ only the term $b_0$ is singular, and it has a
simple pole; hence the sum is meromorphic on $\C^\times$ with simple poles contained in $\{z_b\}_{b\in B(p)}$. Finally, for
$|z|\geq 2\varepsilon_{\max}$ we have $|z_b/z|\leq\frac12$ for all $b$, and
\[
\frac{R_b}{z_b^2}\,\frac{1}{1-z/z_b}=-\frac{1}{z}\,\frac{R_b}{z_b}\,\frac{1}{1-z_b/z},
\]
so the same estimates show that $\mathcal J_p(x^*;z)/(-z)$ is holomorphic in $z^{-1}$ near $z^{-1}=0$ and equal to $1$ there. Since
$-z\,e^{-u_p/z}=-z+u_p+O(z^{-1})$, the function $J_p(\tau,z)+z$ is holomorphic at $z=\infty$. 

\paragraph{Proof of Corollary \ref{cor:none-J}.}
Fix an admissible $\lambda\in\ft$ as in Theorem~\ref{thm:J-convergence}. In Steps~1 and~2 all equivariant
quantities are specialized at $\lambda$. Put $D_i:=c_1(L_i)$, where $L_1,\dots,L_r$ are the ample line bundles
used to define the Novikov variables, so that $q^\beta=q_1^{\langle D_1,\beta\rangle}\cdots q_r^{\langle D_r,\beta\rangle}$
and $\langle D_i,\beta\rangle\geq1$ for every $i$ and every effective $\beta\neq0$. In particular, $|q^{[C]}|$ is small
for all $T$-curves $C$ if all $|q_i|$ are small.

\emph{Step 1: equivariant lifts.}
Let $\{\varphi_a\}_{a=1}^N$ be a homogeneous basis of $H^*(X;\C)$. Fixing a fiber of the principal $T$-bundle
$ET\to BT$ gives an inclusion $T\hookrightarrow ET$, hence $T\times X\hookrightarrow ET\times X$ and, after taking the
quotient by $T$, an inclusion $i\colon X\hookrightarrow ET\times_T X$. The induced restriction map
\[
i^*\colon H^*_T(X)=H^*\big((ET\times X)/T\big)\longrightarrow H^*\big((T\times X)/T\big)=H^*(X)
\]
is surjective, because $X$ is equivariantly formal. Moreover, by the Leray--Hirsch theorem, any homogeneous
lifts $\widetilde\varphi_a\in H^*_T(X)$ with $i^*\widetilde\varphi_a=\varphi_a$ form a basis of $H^*_T(X)$ over
$H^*_T({\rm pt})\cong\C[\ft]$; see \cite{GKM1998}. In other words, $H^*_T(X)\cong H^*(X)[\ft]$ as $\C[\ft]$-modules
(but in general not as rings). We write $\tau=\sum_at_a\widetilde\varphi_a$ in the equivariant theory and
$\tau=\sum_at_a\varphi_a$ in the non-equivariant one, so that the two are related by $i^*$. Below, $\deg$ denotes the
complex degree and $n=\dim_\C X$.

\emph{Step 2: convergence of the genus-$0$ potential.}
For $p\in X^T$ we have $\tau|_p=\sum_at_a\,\widetilde\varphi_a|_p(\lambda)$, so $|\tau|_p|<r$ for all $p$ if all $|t_a|$
are small. Hence, by Theorem~\ref{thm:J-convergence}, there is a polydisc $\Delta=\{|q_i|<r',\ |t_a|<r'\}$ on which
every $J_p(\tau,z)$ converges. Its poles lie in the set $\{z_b\}_{b\in B(p)}$, which is bounded since
$|z_b|\leq|\varepsilon_C(p)|$; choose $R$ larger than all $|z_b|$. By the proof of Theorem~\ref{thm:J-convergence}, the
convergence is uniform for $z$ on the circle $|z|=R$. By the localization formula,
\[
\big(J(\tau,z),\widetilde\varphi_a\big):=\int_X J(\tau,z)\cup\widetilde\varphi_a
=\sum_{p\in X^T}\frac{J_p(\tau,z)\,\widetilde\varphi_a|_p}{e_p}.
\]
The coefficient of $(-z)^{-1}$ in this expression is
\[
\sum_{\beta\in\operatorname{NE}(X)}\sum_{m\geq0}\frac{q^\beta}{m!}
\langle \widetilde\varphi_a,\tau,\dots,\tau\rangle^T_{0,m+1,\beta}
=\frac{\partial F^{(0)}}{\partial t_a}(q,\tau),
\]
where the terms with $\beta=0$ and $m<2$ are omitted. This coefficient equals the Cauchy integral
$-\frac{1}{2\pi i}\oint_{|z|=R}(J,\widetilde\varphi_a)\,dz$, so it is holomorphic on $\Delta$. The
$\tau$-independent part can also be expressed as a Cauchy integral thanks to the divisor equation
\[ 
q_1\partial_{q_1} F^{(0)}(q,0)=
\sum_{\beta\neq0}\langle D_1,\beta\rangle\; q^\beta\langle\ \rangle^T_{0,0,\beta}=
-\frac{1}{2\pi i}\oint_{|z|=R}(J(0,z),\widetilde D_1)\,dz,
\]
where $\widetilde D_1$ is an equivariant lift of the divisor class $D_1$. 
Hence $F^{(0)}$ is holomorphic on $\Delta$. The radius of convergence of a Taylor series does not change when we differentiate it term-wise so $F^{(0)}(q,0)$ must be convergent too. 

Write
\[
F^{(0)}(q,\tau)=\sum_{\beta}\sum_{\alpha\in\mathbb{N}^N}F^{(0)}_{\beta,\alpha}\,q^\beta t^\alpha,
\qquad
F^{(0)}_{\beta,\alpha}=\frac{1}{\alpha!}\big\langle\widetilde\varphi_1^{\,\alpha_1},\dots,
\widetilde\varphi_N^{\,\alpha_N}\big\rangle^T_{0,|\alpha|,\beta}\in\C[\ft],
\]
where $\widetilde\varphi_a^{\,\alpha_a}$ means that $\widetilde\varphi_a$ is inserted $\alpha_a$ times. By the dimension
constraint, $F^{(0)}_{\beta,\alpha}$ is a homogeneous polynomial of degree
\[
d(\beta,\alpha)=\sum_{a=1}^N\alpha_a\deg\varphi_a+3-|\alpha|-n-\int_\beta c_1(TX),
\]
and it vanishes if $d(\beta,\alpha)<0$. If $d(\beta,\alpha)=0$, then $F^{(0)}_{\beta,\alpha}$ is a constant, equal to
its value at $0\in\ft$, i.e., to the corresponding non-equivariant invariant
\[
\frac{1}{\alpha!}\big\langle\varphi_1^{\,\alpha_1},\dots,\varphi_N^{\,\alpha_N}\big\rangle_{0,|\alpha|,\beta}.
\]
If $d(\beta,\alpha)>0$, the corresponding non-equivariant invariant vanishes for dimensional reasons. Therefore the
non-equivariant genus-$0$ potential is
\[
F^{(0)}_X(q,\tau)=\sum_{d(\beta,\alpha)=0}F^{(0)}_{\beta,\alpha}(\lambda)\,q^\beta t^\alpha .
\]
The Taylor series of the holomorphic function $F^{(0)}$ converges absolutely on compact subsets of $\Delta$, and so
does any subseries of it. Hence $F^{(0)}_X$ converges on $\Delta$. In particular, the big quantum product
$\bullet_{q,\tau}$ of $X$, whose structure constants are the third derivatives of $F^{(0)}_X$, is holomorphic on
$\Delta$.

\emph{Step 3: convergence of the non-equivariant $J$-function.}
From now on we work in non-equivariant cohomology, and we put $w:=z^{-1}$. We have
$J(\tau,z)=-z\,S(q,\tau,z)^{-1}1$, where
\[
S(q,\tau,z)=1+\sum_{k\geq1}S_k(q,\tau)\,z^{-k},\qquad S_k(q,\tau)\in\operatorname{End}(H^*(X))[[q,t]],
\]
is the \emph{calibration} of $X$ \cite{Givental2001b}. It satisfies
\begin{equation}\label{eq:S-tau}
z\,\partial_{t_a}S=\varphi_a\bullet_{q,\tau}S,\qquad 1\leq a\leq N,
\end{equation}
and, by the divisor equation,
\begin{equation}\label{eq:S-q}
z\,q_i\partial_{q_i}S=D_i\bullet_{q,\tau}S-S\,(D_i\cup),\qquad 1\leq i\leq r.
\end{equation}
Equivalently, the series $P:=S\,q_1^{D_1\cup/z}\cdots q_r^{D_r\cup/z}$ satisfies
\[
z\,q_i\partial_{q_i}P=D_i\bullet_{q,\tau}P,\qquad 1\leq i\leq r.
\]
Moreover, $S(q,\tau,-z)^*\,S(q,\tau,z)=1$, where ${}^*$ denotes the adjoint with respect to the Poincar\'e pairing
\cite{Givental2001b}. We show that $S$ converges, and is holomorphic, for $(q,\tau)$ in a polydisc
$\Delta'\subset\Delta$ and $w\in\C$.

\smallskip
(a) \emph{The case $\tau=0$.} Put $S^0(q,w):=S(q,0,z)$, and write $q=(q_1,q')$ with $q'=(q_2,\dots,q_r)$. Since $D_1$ is
ample, every quantum correction contains a positive power of $q_1$. Hence
\[
S^0|_{q_1=0}=1,\qquad D_1\bullet_{q,0}\big|_{q_1=0}=D_1\cup .
\]
Equation \eqref{eq:S-q} with $i=1$ and $\tau=0$ is the linear system
\begin{equation}\label{eq:fuchs}
q_1\frac{\partial S^0}{\partial q_1}=\mathcal{A}(q_1)\,S^0,\qquad
\mathcal{A}(q_1)Y:=w\big(D_1\bullet_{q,0}Y-Y\,D_1\cup\big),
\end{equation}
for an $\operatorname{End}(H^*(X))$-valued function $Y$. By Step~2, $\mathcal{A}$ is holomorphic in $q_1$, and in the
parameters $(q',w)$, for $q$ in a polydisc and $w\in\C$. Thus \eqref{eq:fuchs} has a regular singular point at
$q_1=0$, with residue $\mathcal{A}(0)=\operatorname{ad}(w\,D_1\cup)$, which is nilpotent for all values of the
parameters. Write $S^0=\sum_{k\geq0}S^0_k(q',w)\,q_1^k$ and $\mathcal{A}(q_1)=\sum_{j\geq0}\mathcal{A}_j\,q_1^j$. Then
\eqref{eq:fuchs} is equivalent to the recursion
\[
\big(k-\operatorname{ad}(w\,D_1\cup)\big)S^0_k=\sum_{j=1}^{k}\mathcal{A}_j\,S^0_{k-j},\qquad S^0_0=1 .
\]
Since $\operatorname{ad}(w\,D_1\cup)$ is nilpotent, the operator $k-\operatorname{ad}(w\,D_1\cup)$ is invertible for
$k\geq1$. Hence \eqref{eq:fuchs} has a unique formal solution with $S^0_0=1$, and this solution is the Gromov--Witten
series $S^0$. Note that, since $D_1$ is ample, each $S^0_k$ is a polynomial in $(q',w)$. By the classical theorem on
systems with a regular singular point (see, e.g., \cite{Bolibruch2009}), this formal solution is convergent, and
$P_1:=S^0\,q_1^{w\,D_1\cup}$ is the fundamental solution of \eqref{eq:fuchs} normalized by $S^0|_{q_1=0}=1$. Since the
residue is non-resonant for all values of the parameters, the radius of convergence and the bounds in the proof are
uniform for $(q',w)$ in compact sets. Hence $S^0$ is holomorphic in $(q,w)$ on a 
$\{\text{polydisc in $q$}\}\times \C$.

\smallskip
(b) \emph{General $\tau$.} Restricting \eqref{eq:S-tau} to the ray $s\mapsto s\tau$, $0\leq s\leq1$, gives the linear
ODE
\[
\frac{d}{ds}S(q,s\tau,z)=w\Big(\sum_at_a\,\varphi_a\bullet_{q,s\tau}\Big)S(q,s\tau,z),\qquad S(q,0,z)=S^0(q,w),
\]
whose coefficients are holomorphic in $(s,q,t,w)$ for $(q,\tau)\in\Delta'$ and $w\in\C$. Its solution is therefore
holomorphic in $(q,t,w)$. Its Taylor expansion in $(q,t)$ satisfies \eqref{eq:S-tau} with the same value at $\tau=0$,
so it coincides with the formal series $S$.

\smallskip
(c) \emph{Conclusion.} By the unitarity relation, $S(q,\tau,z)^{-1}=S(q,\tau,-z)^*$ is holomorphic as well. Hence
$J(\tau,z)=-z\,S(q,\tau,z)^{-1}1$ converges for $(q,\tau)\in\Delta'$ and every $z\in\C^\times$, and $J/z$ is
holomorphic in $w=z^{-1}\in\C$, i.e., on $\PP^1\setminus\{0\}$. \qed

\subsection{Equivariant central charge}\label{sec:central-charge}

We specialize the equivariant parameters to an admissible point $\lambda\in\ft$ in the sense of Definition \ref{def:admissible} and write $u$ for $u(\lambda)\in\C$. In particular the points 
$z_b$ ($b=(p,C,k)\in B$) are pairwise distinct and non-zero.

Let $\Gh^*_X:=\Gh_{T^*X}=\prod_x\Gamma(1-x)$, the products over the equivariant Chern roots $x$ of
$TX$, be the dual Gamma class \cite{Libgober1999,Iritani2009,KKP2008,GGI2016}. Rescaling the degree-$2j$
part by $z^{-j}$ and restricting to $p$ we get the meromorphic functions of $z\in\C^\times$
\[
\Ghs_p(z):=\prod_{u\in\wts(T_pX)}\Gamma\Big(1-\frac uz\Big).
\]
$\Gamma(1-u/z)$ has simple poles exactly at $z=u/k$, $k\ge1$, and no zeros. So $1/\Ghs_p$ is holomorphic on $\PP^1\setminus\{0\}$, equal to $1$ at $\infty$, and its zeros are simple and are exactly the points $z_b$ with $b\in B(p)$, i.e.\ the possible poles of $J_p(\tau,z)$ (Lemma~\ref{lem:poles}). Therefore, it is natural to introduce the ratio 
\[
G_p(\tau,z)=\frac{J_p(\tau,z)}{-z\Ghs_p (z)},\qquad p\in X^T.
\]
For $q_i\neq0$ and fixed branches of $\log q_i$ and $\log z$, we also put
\[
H_p(\tau,z):=G_p(\tau,z)\,\exp\Big(\frac{c_1^T(TX)|_p\,(\log z-\pi i)}{z}-\sum_{i=1}^r\frac{P_i|_p\log q_i}{z}\Big),
\qquad p\in X^T,
\]
where $P_i=c_1^T(L_i)$ is as in the introduction.

\begin{lemma}\label{lem:gammanorm}
Let $X$ be a GKM manifold, $\lambda$ admissible, and $r=r(\lambda)>0$ as in Theorem~\ref{thm:J-convergence}.
\begin{enumerate}
\item[(a)] The ratio $G_p(\tau,z)$ is holomorphic for all $(q,\tau,z)$ such that
\[
|q^{[C]}|<r\ (\text{all $T$-curves $C$}),\qquad |\tau|_p|<r,\qquad z\in \PP^1\setminus \{0\}.
\]
\item[(b)] For every $b=(p,C,k)\in B$, where $p'\in X^T$ is the fixed point connected to $p$ via the $T$-curve $C$,
\[
G_p(z_b)=q^{k[C]}(-z_b)^{-kc_1(C)}\,G_{p'}(z_b).
\]
\item[(c)] For every $b=(p,C,k)\in B$ we have $H_p(z_b)=H_{p'}(z_b)$, for any fixed choice of the branches of $\log z$ and
$\log q_i$.
\end{enumerate}
\end{lemma}
\begin{proof}
(a) Since the gamma function in the denominator cancels the poles of the numerator, Theorem~\ref{thm:J-convergence}
implies that $G_p(\tau,z)$ is analytic in the stated domain.

(b) Since we can 
work with the Taylor expansions in $(q,t)$, it is sufficient to prove that
coefficient-wise in $q$ and $t$:
\[
\;G_p(\tau,z_b)=q^{k[C]}\,\big(
-z_b\big)^{-k\,c_1(C)}\;
G_{p'}(\tau,z_b), 
\]
where $b=(p,C,k)$. By definition $z_b:=\varepsilon/k$ where $\varepsilon:=e^T(T_pC) =kz_b$. 
Using the splitting of Lemma~\ref{lem:edgegen}(a): weights $u_1=\varepsilon,u_2,\dots$ at
$p$, $u_i-a_i\varepsilon=u_i-ka_iz_b$ at $p'$.

\emph{Residue of the Gamma factor.} With $s=1-\varepsilon/z$ we have $s(z_b)=1-k$ and $ds/dz=\varepsilon/z_b^2=k/z_b$, so
\[
\Res_{z=z_b}\Gamma\Big(1-\frac\varepsilon z\Big)=\frac{(-1)^{k-1}}{(k-1)!}\cdot\frac{z_b}{k}=\frac{(-1)^{k-1}z_b}{k!},
\]
and the other factors $\Gamma(1-u_i/z_b)$, $i\ge2$, are finite and non-zero.

\emph{Using Theorem~\ref{thm:J}.} Since both $J_p(\tau,z)$ and $\Ghs_p(z)$ have simple poles at $z=z_b$, we have $G_p(z_b)=(\Res_{z=z_b}J_p(\tau,z))\big/\big(-z_b\Res_{z=z_b}\Ghs_p(z)\big)$, and $\Res_{z=z_b}J_p(\tau,z)=-R_{b}$ with
\[
R_{b}=q^{k[C]}\frac{e_p}{kE_{k}(C)}\,J_{p'}(\tau,z_b)= 
q^{k[C]}\frac{e_p}{kE_{k}(C)}\,
(-z_b)\,\Ghs_{p'}(z_b)\,G_{p'}(z_b).
\]
Therefore
\[
\frac{G_p(z_b)}{q^{k[C]}G_{p'}(z_b)}=\frac{(-1)^{k}\,k!\;e_p\;\Ghs_{p'}(z_b)}{k\,z_b\,E_{k}(C)\,\prod_{i\ge2}\Gamma(1-u_i/z_b)} .
\]

\emph{The factors.} We use $e_p=\prod_iu_i$, Lemma~\ref{lem:edgegen}(b), and $u_i-\frac mk\varepsilon=u_i-mz_b$.
\begin{itemize}
\item $i=1$: the weight of $T_{p'}C$ is $-\varepsilon$ and $\Gamma(1+\varepsilon/z_b)=k!$; the tangent part of $E_{k}(C)$ is
$\prod_{j=-k,\,j\ne0}^{k}(jz_b)=(-1)^k(k!)^2z_b^{2k}$; and $e_p$ contributes $\varepsilon=kz_b$. The factor is
$\frac{(-1)^kk!\cdot kz_b\cdot k!}{k\,z_b\,(-1)^k(k!)^2z_b^{2k}}=z_b^{-2k}$.
\item $i\ge2$, $a_i\ge0$: $\dfrac{\Gamma(1-u_i/z_b+ka_i)}{\Gamma(1-u_i/z_b)}=\prod_{m=1}^{ka_i}\Big(m-\frac{u_i}{z_b}\Big)=(-z_b)^{-ka_i}\prod_{m=1}^{ka_i}(u_i-mz_b)$, the
corresponding part of $E_k(C)$ is $\prod_{m=0}^{ka_i}(u_i-mz_b)$, and $e_p$ contributes $u_i$: the factor is $(-z_b)^{-ka_i}$.
\item $i\ge2$, $a_i<0$: the Gamma ratio is $(-z_b)^{-ka_i}\prod_{m=ka_i+1}^{0}(u_i-mz_b)^{-1}$, the $H^1$-part contributes
$\prod_{m=ka_i+1}^{-1}(u_i-mz_b)$ to $1/E_k(C)$, and $e_p$ contributes $u_i$: the factor is again $(-z_b)^{-ka_i}$.
\end{itemize}
The product is $z_b^{-2k}(-z_b)^{-k\sum_{i\ge2}a_i}=(-z_b)^{-kc_1(C)}$.

(c) Let $L$ be a $T$-equivariant line bundle on $X$. By Proposition~\ref{prop:eqBG}, $L|_C\cong L(a,w)$ with $a=\langle c_1(L),[C]\rangle$,
so the fiber weights of $L$ at $p$ and $p'$ satisfy $c_1^T(L)|_p-c_1^T(L)|_{p'}=a\,\varepsilon=k\,a\,z_b$. Applied to
$L=L_i$ and to $L=\det TX$, this gives
\[
P_i|_p-P_i|_{p'}=k\,\langle P_i,[C]\rangle\,z_b,\qquad c_1^T(TX)|_p-c_1^T(TX)|_{p'}=k\,c_1(C)\,z_b .
\]
Hence
\[
\frac{H_p(z_b)}{H_{p'}(z_b)}=\frac{G_p(z_b)}{G_{p'}(z_b)}\,e^{k\,c_1(C)(\log z_b-\pi i)}\prod_{i=1}^r e^{-k\langle P_i,[C]\rangle\log q_i}
=\frac{G_p(z_b)}{G_{p'}(z_b)}\,(-z_b)^{k\,c_1(C)}\,q^{-k[C]}=1,
\]
where we used that $k\,c_1(C)$ and $k\langle P_i,[C]\rangle$ are integers, so the result does not depend on the branches,
and part (b).
\end{proof}

\paragraph{Proof of Corollary~\ref{cor:central-charge}.}
Let $r$ be as in Lemma~\ref{lem:gammanorm} and put $r_0:=\min(r,1)$. Since every $L_i$ is ample,
$\langle c_1(L_i),[C]\rangle\geq1$ for all $T$-curves $C$, so $|q_i|<r_0$ for all $i$ implies $|q^{[C]}|<r$ for all $C$.
The central charge is additive in $E$, and $K^0_T(X)$ is generated by $T$-equivariant vector bundles, so it suffices
to consider a $T$-equivariant vector bundle $E$. By Theorem~\ref{thm:J-convergence}, every term of \eqref{eq:Z-loc} is
holomorphic on the universal cover of the domain in the corollary, except for poles at the points $z_b$.

\emph{Step 1: a second formula for $Z(E)$.} Write $J_p(\tau,z)/(-z)=\Ghs_p(z)\,G_p(\tau,z)$. By the reflection formula
$\Gamma(1-x)\Gamma(1+x)=\pi x/\sin(\pi x)$ and the identity $1/\sin(\pi x)=2i\,e^{-\pi i x}/(1-e^{-2\pi i x})$,
\[
\frac{1}{e_p}\,\Ghs_p(z)\prod_{u\in\wts(T_pX)}\Gamma\Big(1+\frac uz\Big)
=\frac1{e_p}\prod_{u}\frac{\pi u/z}{\sin(\pi u/z)}
=\frac{(2\pi i)^n}{z^n}\,\frac{e^{-\pi i\,c_1^T(TX)|_p/z}}{\prod_u\big(1-e^{-2\pi iu/z}\big)},
\]
since $e_p=\prod_uu$ and $c_1^T(TX)|_p=\sum_uu$. Substituting in \eqref{eq:Z-loc}, we get
\begin{equation}\label{eq:Z-H}
Z(E)(\tau,z)=i^n\Big(\frac{2\pi}{z}\Big)^{\frac n2}\sum_{p\in X^T}H_p(\tau,z)\,\Phi_p(z),\qquad
\Phi_p(z):=\frac{\operatorname{ch}_T(E_p)(2\pi i/z)}{\prod_{u\in\wts(T_pX)}\big(1-e^{-2\pi iu/z}\big)},
\end{equation}
where $z^{-n/2}:=e^{-\frac n2\log z}$. By Lemma~\ref{lem:gammanorm}(a), $H_p$ is holomorphic on the universal cover of the
domain. The function $\Phi_p$ is meromorphic in $z\in\C^\times$.

\emph{Step 2: proof of (a).} The poles of $\Phi_p$ are among the points $z_*=u/k$ with $u\in\wts(T_pX)$ and
$k\in\Z\setminus\{0\}$. Fix such a point $z_*$. If two distinct weights $u,u'$ of $T_pX$ satisfied $u/z_*\in\Z$ and
$u'/z_*\in\Z$, then $u/u'$ would be real, contradicting admissibility (ii). Hence $\Phi_p$ has at most a simple pole at $z_*$,
coming from the unique weight $u$ of $T_pX$ with $u/z_*\in\Z\setminus\{0\}$. Let $C$ be the $T$-curve through $p$ with
$\varepsilon_C(p)=u$, and let $p'$ be its other fixed point. Then $\varepsilon_C(p')=-u$, so $\Phi_{p'}$ also has a pole at
$z_*$, and $C$ is the unique $T$-curve through $p'$ responsible for it. Therefore the fixed points $p$ for which $\Phi_p$ is
singular at $z_*$ split into disjoint pairs $\{p,p'\}$, the two ends of $T$-curves $C$, and it suffices to prove that
$H_p\Phi_p+H_{p'}\Phi_{p'}$ has a removable singularity at $z_*$. Exchanging $p$ and $p'$ if necessary, we may assume that
$z_*=\varepsilon/k$ with $\varepsilon:=\varepsilon_C(p)$ and $k\geq1$, i.e., $z_*=z_b$ with $b=(p,C,k)\in B$.

Write $\Phi_p(z)=\Phi^\flat_p(z)/(1-e^{-2\pi i\varepsilon/z})$ and $\Phi_{p'}(z)=\Phi^\flat_{p'}(z)/(1-e^{2\pi i\varepsilon/z})$, where
$\Phi^\flat_p$ and $\Phi^\flat_{p'}$ are holomorphic near $z_b$. By Lemma~\ref{lem:edgegen}(a), the remaining weights of
$T_{p'}X$ are $u_i-a_i\varepsilon$ ($i\geq2$), and by Proposition~\ref{prop:eqBG}, $E|_C\cong\bigoplus_jL(c_j,w_j)$, so the weights of
$E_{p'}$ are $w_j-c_j\varepsilon$. Since $e^{-2\pi i(u_i-a_i\varepsilon)/z_b}=e^{-2\pi iu_i/z_b}e^{2\pi i a_ik}=e^{-2\pi iu_i/z_b}$ and
similarly for the weights of $E$, we get $\Phi^\flat_p(z_b)=\Phi^\flat_{p'}(z_b)$. A direct computation gives
\[
\Res_{z=z_b}\frac{dz}{1-e^{-2\pi i\varepsilon/z}}=-\frac{z_b}{2\pi ik},\qquad
\Res_{z=z_b}\frac{dz}{1-e^{2\pi i\varepsilon/z}}=\frac{z_b}{2\pi ik}.
\]
Therefore
\[
\Res_{z=z_b}\big(H_p\Phi_p+H_{p'}\Phi_{p'}\big)dz=\frac{z_b}{2\pi ik}\,\Phi^\flat_p(z_b)\,\big(H_{p'}(z_b)-H_p(z_b)\big)=0
\]
by Lemma~\ref{lem:gammanorm}(c). Since the pole is simple, the singularity is removable. This proves that $Z(E)$ is holomorphic
on the universal cover of the domain, and in particular (a).

\emph{Step 3: proof of (b).} In \eqref{eq:Z-loc}, the functions $J_p(\tau,z)$, $\Gamma(1+u/z)$ and $\operatorname{ch}_T(E_p)(2\pi i/z)$
are single-valued in $z$, and only $z^{n/2}$ and $e^{c_1^T(TX)|_p\log z/z}$ depend on the branch of $\log z$. Along the loop
$z\mapsto e^{2\pi is}z$ they are multiplied by $(-1)^n$ and by $e^{2\pi i\,c_1^T(TX)|_p/z}$. Since $c_1^T(TX)|_p$ is the weight of
$K_X^{-1}|_p=\det T_pX$, we have $e^{2\pi i\,c_1^T(TX)|_p/z}\operatorname{ch}_T(E_p)(2\pi i/z)=\operatorname{ch}_T\big((E\otimes K_X^{-1})_p\big)(2\pi i/z)$.
Hence $Z(E)$ is transformed into $(-1)^nZ(E\otimes K_X^{-1})$. Similarly, only $e^{-P_i|_p\log q_i/z}$ depends on the branch of
$\log q_i$; along the loop $q_i\mapsto e^{2\pi is}q_i$ it is multiplied by $e^{-2\pi iP_i|_p/z}$, and $P_i|_p$ is the weight of
$L_i|_p$. Hence $Z(E)$ is transformed into $Z(E\otimes L_i^{-1})$.

\emph{Step 4: proof of (c).} Let $\delta>0$ be the smallest angle between $S$ and the excluded rays. We use \eqref{eq:Z-loc}.

First, by Theorem~\ref{thm:J-convergence} and its proof, $J_p(\tau,z)/(-z)=e^{-u_p/z}g_p(z)$ with
\[
g_p(z):=1+\sum_{b\in B(p)}\frac{R_b}{z_b(z-z_b)}\,e^{u_p/z_b},\qquad
g_p(z)-g_p(0)=z\sum_{b\in B(p)}\frac{R_b}{z_b^2(z-z_b)}\,e^{u_p/z_b}.
\]
For $b=(p,C,k)$, the point $z_b$ lies on the ray $\arg z=\arg\varepsilon_C(p)$, which is excluded. Hence for $z\in S$ we have
$|z-z_b|\geq s_\delta|z_b|$, where $s_\delta:=\sin\min(\delta,\pi/2)$. Using $|R_b|<M\rho^k$, $|u_p|<\eta$ and $|1/z_b|\leq\gamma k$ as in
the proof of Theorem~\ref{thm:J-convergence}, we get
\[
|g_p(z)-g_p(0)|\leq\frac{|z|}{s_\delta}\sum_{b\in B(p)}|R_b|\,e^{|u_p|/|z_b|}\,|z_b|^{-3}
\leq\frac{|z|}{s_\delta}\,nM\gamma^3\sum_{k\geq1}k^3(e\rho)^k ,
\]
so $g_p(z)=g_p(0)+O(z)$ for $z\in S$.

Second, for $z\in S$ and every weight $u$ of $T_pX$ we have $|\arg(u/z)|\leq\pi-\delta$. By Stirling's formula
\cite{WhittakerWatson}, $\Gamma(1+x)=\sqrt{2\pi}\,e^{(x+\frac12)\operatorname{Log}x-x}\,(1+O(1/x))$ uniformly for $|\arg x|\leq\pi-\delta$
as $|x|\to\infty$. With $x=u/z$ and $\operatorname{Log}(u/z)=\ell_u-\log z$, and using $c_1^T(TX)|_p=\sum_uu$, we get
\[
z^{\frac n2}\prod_{u\in\wts(T_pX)}\Gamma\Big(1+\frac uz\Big)\;e^{c_1^T(TX)|_p\log z/z}
=(2\pi)^{\frac n2}\prod_{u}e^{\ell_u/2}\;\exp\Big(\frac1z\sum_uu\,(\ell_u-1)\Big)\,\big(1+O(z)\big).
\]
Finally, $\operatorname{ch}_T(E_p)(2\pi i/z)=\sum_wm_w(E_p)\,e^{2\pi iw/z}$. Substituting these three expansions in \eqref{eq:Z-loc}
gives the formula in (c). Since $e^{\ell_u}=u$, we have $\big(\prod_ue^{\ell_u/2}\big)^2=e_p$, so $\kappa_p^2=e_p^{-1}$. \qed

\subsection{Evidence for the equivariant mirror symmetry conjecture}\label{sec:mirror-conj}

\begin{remark}\label{rem:mirror-consistency}
Conjecture~\ref{conj:mirror} is consistent with Corollary~\ref{cor:central-charge}.
\begin{enumerate}
\item[(i)] An oscillatory integral over a rapid-decay cycle has no poles in $z\in\C^\times$, in agreement with part (a).
\item[(ii)] Changing the branch of $\log g$ by $2\pi i\mu$ multiplies the integrand by $e^{2\pi i\mu(\lambda)/z}$. By
\eqref{eq:Z-loc} the same factor appears in $Z(E\otimes\C_\mu)=e^{2\pi i\mu(\lambda)/z}Z(E)$, where $\C_\mu$ is the trivial line
bundle with fiber the character $\mu$. This is the compatibility of the $\Z[X^*(T)]$-module structures, and it matches the
terms $-2\pi iw$ in the exponents $F_{p,w}$ of part (c). Similarly, the terms $P_i|_p\log q_i$ in $F_{p,w}$ and the
monodromy in part (b) indicate that the Novikov variables enter $F$ through $\log q_i$.
\item[(iii)] By the stationary phase method, the asymptotics of the integral as $z\to0$ in a sector is a sum of
contributions $e^{-F(y)/z}\big(\operatorname{Hess}_\omega F(y)\big)^{-1/2}(1+O(z))$ over critical points $y$ of $F|_{Y_s}$. Part (c)
therefore predicts that the critical points correspond to the fixed points $p\in X^T$, that the critical values are
the numbers $F_{p,w}$ (on the various branches of $\log g$), and that $\operatorname{Hess}_\omega F=e_p/g_p(0)^2$ at the critical
point corresponding to $p$. The term $-\sum_uu(\ell_u-1)$ in $F_{p,w}$ is the critical value of
$\sum_{u\in\wts(T_pX)}(e^{y_u}-u\,y_u)$, the equivariant mirror of $T_pX\cong\C^n$ \cite{Givental1998}, and with
$\omega=dy_1\wedge\dots\wedge dy_n$ its Hessian at the critical point is $\prod_uu=e_p$. This agrees with the prediction at
$q=0$, where $g_p(0)=1$, and it also fixes the normalization $(2\pi z)^{-n/2}$ in Conjecture~\ref{conj:mirror}.
\end{enumerate}
\end{remark}

\paragraph{Partial flag manifolds.} We follow the notation of Chow \cite{Chow2025}. Let $G^\vee$ be a simply connected
complex semisimple group, $P^\vee\subset G^\vee$ a parabolic subgroup and $X=G^\vee/P^\vee$, with the action of a maximal torus
$T^\vee\subset P^\vee$ (so the roles of $T$ and $T^\vee$ are interchanged compared with Conjecture~\ref{conj:mirror}). Let $G$ be the Langlands dual group, which is of adjoint type, and $T\subset G$ the
maximal torus dual to $T^\vee$. Let $B=TU$ and $B^-=TU^-$ be opposite Borel subgroups, $\{\alpha_i\}_{i\in I}$ the simple roots,
and fix $e_i\in\mathfrak g_{\alpha_i}$, $f_i\in\mathfrak g_{-\alpha_i}$ with $[e_i,f_i]=\alpha_i^\vee$. Let $x_i:\C\to U$ and $y_i:\C\to U^-$
be the homomorphisms $x_i(a)=\exp(ae_i)$ and $y_i(b)=\exp(bf_i)$, and let $\chi_i:U\to\C$ be the homomorphism of algebraic
groups determined by $\chi_i\circ x_j=\delta_{ij}\,{\rm id}_\C$; put $\chi:=\sum_{i\in I}\chi_i$. For $w\in W$ we denote by
$\overline w\in N_G(T)$ the Tits lift of $w$, i.e.\ $\overline w=\overline{s}_{i_1}\cdots\overline{s}_{i_m}$ for a reduced
expression $w=s_{i_1}\cdots s_{i_m}$, where $\overline{s}_i:=x_i(-1)\,y_i(1)\,x_i(-1)$; it does not depend on the reduced
expression.

Let $P\subset G$ be the parabolic subgroup corresponding to $P^\vee$, with Levi subgroup $L_P$, let $I_P\subset I$ be the set of
simple roots of $L_P$ and $W_P=\langle s_i\mid i\in I_P\rangle$. Let $w_0$ and $w_0^P$ be the longest elements of $W$ and $W_P$,
and put $w_P:=w_0^Pw_0$; then $\ell(w_P)=\ell(w_0)-\ell(w_0^P)=\dim X$. The center $Z(L_P)=\{t\in T\mid\alpha_i(t)=1,\ i\in I_P\}$
is a torus, and $t\mapsto(\alpha_i(t))_{i\in I\setminus I_P}$ identifies it with the torus of quantum parameters of $X$, one for
each generator of ${\rm Pic}(X)$. Rietsch's mirror family \cite{Rietsch2008} is
\[
\pi_P:X_P:=B^-\cap U\,Z(L_P)\,\overline{w}_P\,U\longrightarrow Z(L_P).
\]
Every $x\in X_P$ admits two factorizations
\[
x=u_0\,t_0=u_1\,t\,\overline{w}_P\,u_2,\qquad u_0\in U^-,\quad t_0\in T,\quad t\in Z(L_P),\quad u_1,u_2\in U,
\]
where $t_0$ and $t$ are uniquely determined by $x$, and we put
\[
\pi_P(x):=t,\qquad \gamma_P(x):=t_0\in T,\qquad f_P(x):=\chi(u_1)+\chi(u_2).
\]
If $P\neq B$, then $u_1$ and $u_2$ are not unique, but $f_P(x)$ does not depend on their choice. The function $f_P$ is the
decoration of Berenstein--Kazhdan's geometric crystal, and it agrees with Rietsch's superpotential; the observation that
Rietsch's mirror can be obtained from geometric crystals is due to Lam \cite{Lam2016} and Chhaibi \cite{Chhaibi}, and the
mirror is described in this form in \cite{LamTemplier}. Finally, the fibers of $\pi_P$ carry a distinguished volume form
$\omega$ without zeros and with simple poles along the boundary divisor, normalized as in \cite{Chow2025}. Thus $Y=X_P$, $\mathcal B=Z(L_P)$, $f=f_P$, $g=\gamma_P$ and
$\omega$ have the form required in Conjecture~\ref{conj:mirror}.

There are many partial results towards Conjecture~\ref{conj:mirror} for partial flag manifolds; see
\cite{Rietsch2008,Lam2016,Chhaibi,LamTemplier,Chow2023,LLYZ2026,HLYZZ2025,Chow2025}. In
particular, mirror symmetry is known at the level of $D$-modules, for the small quantum cohomology by Chow
\cite{Chow2023} and for the big quantum cohomology by Hinault--Li--Yu--Zhang--Zhang \cite{HLYZZ2025}. However, the existence of a primitive form and the integral
structure (the cycles $C(E)$) are open in general; the only result in this direction is due to Chow \cite{Chow2025}: the oscillatory integral of $e^{-f/z}g^{\lambda/z}\,\omega$ over the
totally positive part of the fiber coincides with the Gamma-class central charge of the structure sheaf, for $\tau\in H^2$
and $\lambda$ small compared with $z$; up to conventions, this is Conjecture~\ref{conj:mirror} for $E=\mathcal O_X$ and
$\tau\in H^2$.

\paragraph{Toric manifolds.} We use Givental's notation \cite{Givental1996,Givental1998}. Let $X$ be a smooth projective
toric variety, obtained as the GIT quotient $X=\C^N/\!/K$ of $\C^N$ by a subtorus $K\subset(\C^*)^N$ of dimension $r=N-n$.
After choosing an isomorphism $K\cong(\C^*)^r$, the embedding is given by an integer matrix $(m_{ij})$, $1\leq i\leq r$,
$1\leq j\leq N$,
\[
K\longrightarrow(\C^*)^N,\qquad s\mapsto\Big(\prod_{i=1}^rs_i^{m_{i1}},\dots,\prod_{i=1}^rs_i^{m_{iN}}\Big),
\]
and $T:=(\C^*)^N/K$. Let $K^\vee:=(\C^*)^r$ and let $\pi:(\C^*)^N\to K^\vee$ be the transposed map,
$\pi(u_1,\dots,u_N)=\big(\prod_ju_j^{m_{1j}},\dots,\prod_ju_j^{m_{rj}}\big)$. Since $X$ is smooth, $\Z^N/m^{t}\Z^r$ is torsion free;
it is the lattice $\ft_\Z$ of cocharacters of $T$, and therefore $\operatorname{Ker}\pi\cong{\rm Hom}(\ft_\Z,\C^*)=T^\vee$.
Givental's equivariant mirror of $X$ is the family
\[
\pi:Y:=(\C^*)^N\longrightarrow\mathcal B:=K^\vee,\qquad f(u)=u_1+\cdots+u_N,
\]
whose fibers $Y_q=\{u\in(\C^*)^N\mid\prod_ju_j^{m_{ij}}=q_i,\ 1\leq i\leq r\}$ are $T^\vee$-torsors; the coordinates
$q_i$ of $K^\vee$ are the Novikov variables of the basis of $H^2(X;\Z)\cong{\rm Hom}(K,\C^*)$ given by the characters $s\mapsto
s_i$. The volume form on the fibers is the residue $\omega:=\bigwedge_j\frac{du_j}{u_j}\big/\pi^*\bigwedge_i\frac{dq_i}{q_i}$. To define $g:Y\to T^\vee$ we choose a homomorphism
$\sigma:K^\vee\to(\C^*)^N$ with $\pi\circ\sigma={\rm id}$ and put
$g(u):=u\,\sigma(\pi(u))^{-1}$. If $(\lambda_1,\dots,\lambda_N)\in\C^N={\rm Lie}((\C^*)^N)$ is the lift of
$\lambda\in\ft=\C^N/{\rm Lie}(K)$ determined by $\sigma$, then $\langle\lambda,\log g\rangle=\sum_j\lambda_j\log u_j$ up to a
function of $\log q$. A different choice of
$\sigma$ changes $\langle\lambda,\log g\rangle$ by a linear function of $\log q_1,\dots,\log q_r$, which is constant on the
fibers; this has the same effect on the integral as changing the $T$-linearizations of the line bundles $L_i$ in the
definition of $\Psi(E)$.
Fang \cite{Fang2020} proved, for semi-positive toric orbifolds and $\tau\in H^2$, that (up to normalization conventions) the equivariant central charges of all
equivariant coherent sheaves are given by oscillatory integrals of $e^{-f/z}g^{\lambda/z}\omega$ over integral cycles, namely
the characteristic cycles of the corresponding constructible sheaves (homotopic to the SYZ T-dual branes). In the non-equivariant case, the corresponding statement for the Gamma
integral structure is due to Iritani \cite{Iritani2009}. For big $\tau$, the function $f$ has to be replaced by an unfolding
and $\omega$ by the corresponding primitive form; for toric stacks such primitive forms were constructed in \cite{CCIT2020}.

\appendix

%

\section{The Weierstrass M-test for holomorphic maps}\label{app:M-test}

Let $E$ and $F$ be complex Banach spaces and let $U\subset E$ be an open
subset. We denote by $B(a,r)\subset E$ the open ball
of radius $r$ centered at $a$ and by $\mathcal{L}(E,F)$ the Banach space of bounded linear
maps $E\to F$ with the operator norm. We write $F^*\coloneqq\mathcal{L}(F,\C)$ 
for the continuous dual of $F$.
Recall that a map $g\colon U\to F$ is said to be
\emph{holomorphic} if it is Fr\'echet differentiable at every point $x\in U$, 
i.e. for every $x\in U$ there exists a $\C$-linear map
$Dg(x)\in\mathcal{L}(E,F)$ satisfying the following condition: for every $\epsilon>0$ there exists $\delta>0$ such that 
\[
\|g(x+h)-g(x)-Dg(x)h\|\leq \epsilon \|h\| \qquad 
\forall h\in E: \|h\|<\delta.
\]
The goal of this appendix is to prove the following statement.

\begin{proposition}[Weierstrass M-test]\label{prop:M-test}
Let $f_n\colon U\to F$, $n\geq 1$, be holomorphic maps. Suppose that every
point $a\in U$ has a ball $B(a,2r)\subset U$ and a sequence of non-negative
numbers $(M_n)_{n\geq 1}$, both depending on $a$, such that
\[
\sup_{x\in B(a,2r)}\|f_n(x)\|\leq M_n\quad (n\geq 1),
\qquad
\sum_{n=1}^\infty M_n<\infty .
\]
Then the series $S(x)\coloneqq\sum_{n=1}^\infty f_n(x)$ converges uniformly on
each such ball $B(a,2r)$, the map $S\colon U\to F$ is holomorphic, and
\[
DS(x)=\sum_{n=1}^\infty Df_n(x),\qquad x\in U,
\]
where the series converges absolutely in $\mathcal{L}(E,F)$.
\end{proposition}

The proof uses only the following standard facts from functional analysis
(see, e.g., \cite{Yosida1995}):
\begin{enumerate}
\item[\normalfont(A)] $F$ and $\mathcal{L}(E,F)$ are complete, the latter
  because $F$ is complete.
\item[\normalfont(B)] (Hahn--Banach) If $y\in F$ and $\varphi(y)=0$ for all
  $\varphi\in F^*$, then $y=0$.
\item[\normalfont(C)] If $\gamma$ is a circle in $\C$ and
  $u\colon\gamma\to F$ is continuous, then the Riemann integral
  $\int_\gamma u(\lambda)\,d\lambda\in F$ exists, satisfies
  $\big\|\int_\gamma u(\lambda)\,d\lambda\big\|\leq
  \operatorname{length}(\gamma)\,\sup_\gamma\|u\|$, and
  $\varphi\big(\int_\gamma u(\lambda)\,d\lambda\big)=
  \int_\gamma\varphi(u(\lambda))\,d\lambda$ for every $\varphi\in F^*$.
\item[\normalfont(D)] The Cauchy integral formula for scalar holomorphic
  functions of one complex variable.
\end{enumerate}
For a systematic treatment of holomorphic maps between Banach spaces we refer
to \cite{Mujica1986}.

\begin{lemma}[Vector-valued Cauchy formula]\label{lem:vector-Cauchy}
Let $R>0$ and let $\psi$ be an $F$-valued function, complex differentiable
on a neighborhood of the closed disc $\{|\lambda|\leq R\}\subset\C$. Then for
$|\mu|<R$
\[
\psi(\mu)=\frac{1}{2\pi i}\oint_{|\lambda|=R}\frac{\psi(\lambda)}{\lambda-\mu}\,d\lambda,
\qquad
\psi'(0)=\frac{1}{2\pi i}\oint_{|\lambda|=R}\frac{\psi(\lambda)}{\lambda^2}\,d\lambda .
\]
\end{lemma}
\begin{proof}
The integrands are continuous, so both integrals exist by (C). For every
$\varphi\in F^*$ the function $\varphi\circ\psi$ is a scalar holomorphic
function, and $(\varphi\circ\psi)'(0)=\varphi(\psi'(0))$. Hence, by (D), both
formulas hold after applying $\varphi$ to both sides, where we used (C) to
interchange $\varphi$ with the integrals. Since $\varphi\in F^*$ is arbitrary,
the claim follows from (B).
\end{proof}

\begin{lemma}[Uniform Taylor estimate]\label{lem:Taylor-estimate}
Let $g\colon U\to F$ be holomorphic, let $B(a,2r)\subset U$, and suppose that
$\|g(x)\|\leq M$ for all $x\in B(a,2r)$. Put $\rho\coloneqq r/2$. Then for all
$x\in B(a,r)$:
\begin{enumerate}
\item[\normalfont(i)] $\|Dg(x)\|\leq M/\rho$;
\item[\normalfont(ii)] $\|g(x+h)-g(x)-Dg(x)h\|\leq 2M\rho^{-2}\,\|h\|^2$ for all
  $h\in E$ with $\|h\|\leq\rho/2$.
\end{enumerate}
\end{lemma}
\begin{proof}
Fix $x\in B(a,r)$ and $h\in E$ with $h\neq 0$; put $R\coloneqq\rho/\|h\|$ and
$\psi(\lambda)\coloneqq g(x+\lambda h)$. If $|\lambda|<r/\|h\|$, then
$\|\lambda h\|<r$, so $x+\lambda h\in B(a,2r)$. Since $R<r/\|h\|$, the
function $\psi$ is defined on a neighborhood of the closed disc
$\{|\lambda|\leq R\}$, and $\|\psi\|\leq M$ there. By the chain rule, $\psi$ is
complex differentiable with $\psi'(\lambda)=Dg(x+\lambda h)h$; in particular
$\psi'(0)=Dg(x)h$.

(i) By Lemma~\ref{lem:vector-Cauchy} and (C),
\[
\|Dg(x)h\|=\|\psi'(0)\|\leq\frac{1}{2\pi}\cdot 2\pi R\cdot\frac{M}{R^2}
=\frac{M}{R}=\frac{M}{\rho}\,\|h\| .
\]

(ii) Suppose $\|h\|\leq\rho/2$, i.e., $R\geq 2$. Using the identity
\[
\frac{1}{\lambda-1}-\frac{1}{\lambda}-\frac{1}{\lambda^2}
=\frac{1}{\lambda^2(\lambda-1)}
\]
and Lemma~\ref{lem:vector-Cauchy} with $\mu=1$ and $\mu=0$, we get
\[
g(x+h)-g(x)-Dg(x)h=\psi(1)-\psi(0)-\psi'(0)
=\frac{1}{2\pi i}\oint_{|\lambda|=R}\frac{\psi(\lambda)}{\lambda^2(\lambda-1)}\,d\lambda .
\]
Since $|\lambda-1|\geq R-1\geq R/2$ on the circle $|\lambda|=R$, the norm of
the right-hand side is at most
\[
\frac{1}{2\pi}\cdot 2\pi R\cdot\frac{M}{R^2\,(R/2)}=\frac{2M}{R^2}
=\frac{2M}{\rho^2}\,\|h\|^2 .\qedhere
\]
\end{proof}

\begin{proof}[Proof of Proposition~\ref{prop:M-test}]
Fix $a\in U$, together with $r$ and $(M_n)$ as in the hypothesis, and put
$\rho\coloneqq r/2$. Let $S_N\coloneqq\sum_{n=1}^N f_n$.

\emph{Convergence.} For $N'\geq N$ and $x\in B(a,2r)$ we have
$\|S_{N'}(x)-S_{N}(x)\|\leq\sum_{n>N}M_n$, which tends to $0$ as
$N\to\infty$ uniformly in $x$. By completeness of $F$, the partial sums
converge uniformly on $B(a,2r)$ to a map $S$.

\emph{Convergence of the derivatives.} By
Lemma~\ref{lem:Taylor-estimate}(i) applied to $g=f_n$ and $M=M_n$, we have
$\|Df_n(x)\|\leq M_n/\rho$ for $x\in B(a,r)$. Hence the series
\[
L(x)\coloneqq\sum_{n=1}^\infty Df_n(x)
\]
converges absolutely in $\mathcal{L}(E,F)$, which is complete by (A). Each
$Df_n(x)$ is $\C$-linear, and the $\C$-linear maps form a closed subspace of
$\mathcal{L}(E,F)$; therefore $L(x)$ is $\C$-linear.

\emph{Differentiability.} Let $x\in B(a,r)$ and $h\in E$ with
$\|h\|\leq\rho/2$. All three series $\sum f_n(x+h)$, $\sum f_n(x)$ and
$\sum Df_n(x)h$ converge, so we may combine them termwise. Applying
Lemma~\ref{lem:Taylor-estimate}(ii) to each $f_n$, we obtain
\[
\|S(x+h)-S(x)-L(x)h\|
\leq\sum_{n=1}^\infty\|f_n(x+h)-f_n(x)-Df_n(x)h\|
\leq\frac{2}{\rho^2}\Big(\sum_{n=1}^\infty M_n\Big)\|h\|^2 .
\]
The right-hand side is $o(\|h\|)$ as $h\to 0$. Thus $S$ is Fr\'echet
differentiable at $x$ with $DS(x)=L(x)$, which is $\C$-linear. In particular
this holds at $x=a$. Since $a\in U$ was arbitrary, $S$ is holomorphic on $U$
and $DS=\sum_n Df_n$.
\end{proof}

\bibliographystyle{amsplain}
\bibliography{gkm-bib}

\end{document}